\documentclass[12pt,a4paper]{article}
\usepackage{calc}
\usepackage[all]{xy}
\usepackage[centertags]{amsmath}
\usepackage{latexsym}
\usepackage{amsfonts}
\usepackage{graphicx}
\usepackage{tikz}
\usepackage{cases}
\usepackage{amssymb}
\usepackage{amsthm}
\usepackage{mathptmx}
\usepackage{color}
\usepackage[colorlinks,citecolor=blue,linkcolor=blue]{hyperref}
\usepackage{fancyhdr}
\usepackage{epsfig,float,epstopdf,array,bm,cite}
\usepackage{newlfont}
\usepackage[latin1]{inputenc}
\usepackage[latin1]{inputenc}
\usepackage[english]{babel}
\usepackage{tikz}
\usepackage{graphicx}
\usepackage{t1enc}
\usepackage[english]{babel}
\usepackage{fancybox}
\usepackage{graphicx}
\usepackage{t1enc}
\usepackage[french2]{minitoc}
\usepackage{amsfonts}

\usepackage{natbib}

\usepackage{url}

\newtheorem{theo}{Theorem}[section]
\newtheorem{example}{Example}[section]

\newtheorem{lem}{Lemma}[section]
\newtheorem{remark}{Remark}[section]
\newtheorem{propo}{Proposition}[section]

\newtheorem{defi}{Definition}[section]

\newcommand{\beaa}{\begin{eqnarray*}}
\newcommand{\eeaa}{\end{eqnarray*}}
\numberwithin{equation}{section}
\newcounter{mnotecount}[section]

\renewcommand{\themnotecount}{\thesection.\arabic{mnotecount}}

\newcommand{\N}{\mathbb N}

\newcommand{\R}{\mathbb R}

\newcommand{\vp}{\varphi}

\newcommand{\dsp}{\displaystyle}

\newcommand{\no}{\nonumber}

\newcommand{\warn}[1]
{\protect{\stepcounter{mnotecount}}$^{\mbox{\footnotesize
			$
			\bullet$\themnotecount}}$ \marginpar{
		\raggedright\tiny\em $\!\!\!\!\!\!\,\bullet$\themnotecount: {\bf
			Warning:} #1} }

\author{\bf Frank Kemayou $^{1 }$ \quad Roger Tagne Wafo $^{1\dag}$ \quad  Samuel Bowong$^{1,2}$ \\
{\small $^1$ Laboratory of   Mathematics,
Department of Mathematics and Computer Science,}
\\
{\small  Faculty of Science, University of Douala, PO Box 24157  Douala, Cameroon}
\\
\small{$^2$ {\small IRD,  Sorbonne Université, UMMISCO, F-93143, Bobby, France}}\\
{\small $^\dag$ Corresponding author,  email: rtagnewafo@yahoo.com }
}
\title{Global stability analysis of an age-structured model assessing the impact of \textit{Radopholus similis} on banana-plantain production}
\begin{document}
\date{}
\maketitle

 \section*{Abstract}
 In this paper, we develop and analyse a mathematical model to investigate the interactions between banana and plantain plants and the nematode \textit{Radopholus similis}, a pest species occurring in banana plantations worldwide, with particularly high prevalence in Central Africa.  The model incorporates root infection and mortality rates as functions of root age, providing a more realistic representation of the infection dynamics. We prove that the model is well-defined by establishing the existence and uniqueness of a mild solution using the theory of semi-groups for nonlinear evolutionary systems. We further show that the solution is positive and bounded. Under additional assumptions on the regularity of the parameters and the data, we prove that the mild solution is indeed a classical solution. We then study the asymptotic behaviour of the solution by deriving a threshold parameter 
\(\mathcal{N}\), which determines the stability of the disease-free equilibrium. Finally, we perform a numerical analysis of the proposed model using the semi-implicit Euler method. The biological consistency of the numerical solutions is established, and simulations are carried out to illustrate the theoretical results and estimate yield losses caused by nematodes. We conclude the numerical analysis by implementing an impulsive control strategy, which confirms that the use of nematicides - whether chemical or biological - helps to mitigate the devastating effects of nematodes and enhances crop yield.

 Keywords : \textit{Radopholus similis}, Age-structured, Semi-group theory, Stability, Mild solution,  Numerical analysis, Numerical simulations, Impulsive control.
 
\section{Introduction}
 The production of banana and banana  plantain has played an important role in the ongoing development of many countries and continues to be an essential crop  on which the economy of many countries is based \cite{a14}. Banana and/or banana plantain are cultivated both by small-scale farmers as a household cash crop and by large agricultural companies for local consumption, but especially for export. The primary banana plantain production zones are concentrated in Central and West Africa, followed by Latin America.
 Central Africa, in particular, is recognized as the region with the highest levels of banana and plantain production and consumption, making it a critical area for the crop's economic and food security relevance.
 Cameroon is currently ranked among the world's top producers of sweet banana and banana plantain, with an estimated 2024 annual output of over 5.4 millions tonnes. The government aims to increase this to 7.5 millions tonnes by 2025 and 10 millions tonnes by 2030,~\cite{2024banana}.
 Unfortunately, banana and plantain cultivation in all these regions is increasingly threatened by pests and pathogens, which often have a significant and detrimental impact on the yields of cultivated areas.
 The main pests affecting the banana plantain production  are the banana weevil (\textit{Cosmopolites sordidus}) and nematodes. They primarily attack the roots of the banana plants, penetrating them to feed and reproduce.

  The most harmful nematode in Cameroon is \textit{Radopholus similis} \cite{a13} also known as \textit{burrowing nematodes}, \cite{BN}. It is responsible for significant damage to banana cultivation worldwide, especially in developing countries. Its presence leads to substantial losses in banana production and, consequently, impacts trade (the annual loss of yield due to \textit{Radopholus similis} in Cameroon can reach 20 to 60\%, see~\cite{bridge1995, fogain1998, fogain2000,loubana2007}). 
 This nematode, originating from Africa, is becoming an invasive pest in almost all banana-growing countries worldwide. It feeds on and spends most of its life inside the roots of banana plants \cite{a12}.
 This particular aspect of the nematode's biology makes it difficult to study its behaviour and challenging to control. The perforation of immature roots by \textit{Radopholus similis} causes them to fall prematurely during the fruiting period, leading to significant yield reductions. Therefore, its control is a major issue, especially in Cameroon, where agriculture plays a vital role in income generation and employment.

The existing control methods against the nematode \textit{R. similis}  are of several types. (i) The Chemical control. This involves the application of nematicides to soil to reduce nematode populations. It is generally effective and provides rapid suppression of\textit{ R. similis}, but it is expensive, potentially toxic to humans and non-target organisms, polluting to the environment, and strictly regulated in most countries. 
(ii) 
The biological and natural control. 
The Biological control relies on the use of natural enemies or organic materials to suppress \textit{R. similis} populations. Beneficial microorganisms such as \textit{Paecilomyces lilacinus}, \textit{Pochonia chlamydosporia}, and \textit{Bacillus spp.} have shown nematicidal activity by parasitizing eggs or producing toxic metabolites, \cite{Gowen2005}. Natural products like plant extracts, composts, and agricultural residues have also demonstrated nematicidal effects and can improve soil health and plant resistance, \cite{loubana2007, fogain1998}. Their efficacy can vary with environmental conditions and soil microbiota, \cite{Sarah1996}.
(iii) Agroecologycal control. Agroecological management focuses on creating cropping systems that reduce nematode pressure through ecological balance and soil health improvement. Practices such as crop rotation, intercropping, and the use of cover or trap crops help interrupt the life cycle of \textit{R. similis} by depriving it of a suitable host, \cite{Gowen2005}. Although agroecological practices are environmentally friendly and sustainable, they have some limitations. Their effectiveness is often slower and less predictable compared to chemical control, as results depend on soil conditions, climate, and the presence of beneficial microorganisms.
(iv)
Mathematical models. Mathematical modeling provides a powerful tool for understanding and predicting the population dynamics of \textit{Radopholus similis} and its interactions with host plants and the environment.
They allow to:
\begin{enumerate}
\item[(a)] Simulate nematode population dynamics,
\item[(b)] Virtually evaluate the effectiveness of different control strategies, and
\item[(c)] Predict yield losses and assess associated costs.
\end{enumerate}

So far, several mathematical models have been proposed and studied (see \cite{tankam2020, fotso2022, tankam2021} and the references therein) to assess the impact of \textit{R. similis} on banana and plantain production, as well as to propose various strategies to curb the rapid spread of this infesting pest. In particular, \textit{I. Tankam Chedjou}, \cite{tankam2021},  devoted an entire chapter of his thesis (Chapter III) to the review of existing mathematical models on the subject. All these models have the common feature that they are systems of Ordinary Differential Equations (ODE) which assume that the state variables depend only on time. In this case, it implies that the plant population (or their roots) is considered homogeneous, meaning that each individual experiences the same rates of death and/or infection  at a given time $t$, regardless of differences in age or other individual traits. 
However, this homogeneity assumption is often unrealistic in biological systems. In the case of banana and plantain crops, the roots exhibit age-dependent physiological properties. When modelling the infestation by \textit{Radopholus similis}, it is therefore important to account for the age structure of the root system, since roots of different ages may not be equally vulnerable to nematode attack and may respond differently to infestation dynamics.

In this paper, our main objective is to develop and analyze a mathematical model that describes the interactions between \textit{Radopholus similis} nematodes and banana-plantain plants, while explicitly accounting for the age structure of the root system.

This paper is organized as follows. Section 2 is devoted to the formulation of the model based on the nematode life cycle. In Section 3, we present the mathematical analysis of the model. More precisely, we first address its well-posedness using perturbation semigroup theory, which provides tools to investigate the solvability of our model by exploiting the operators involved in the associated Cauchy problem (Laroche and Perasso \cite{a3}, 2015; Pazy \cite{a4}, 1985; Perasso and Laroche \cite{a10}, 2008). Secondly, we prove the positivity and boundedness of solutions. Thirdly, we perform an asymptotic analysis by computing the disease-free equilibrium and the basic reproduction number, which determine the stability of the system. In Section 4, we carry out the numerical analysis of the model. A semi-implicit Euler scheme is implemented, and the positivity, boundedness, and consistency of the numerical approximations are established. We then present numerical simulations to illustrate the theoretical results, including the application of an impulsive control strategy. Finally, additional simulations of the cumulative and daily production are performed to estimate the percentage of yield losses caused by nematode damage. The last section of the paper is devoted to concluding remarks and discussion.

 \section{Model formulation  }
 \subsection{Biological background of the \textit{Radopholus similis} }
 \textit{Radopholus similis} is a burrowing nematode, meaning it burrows into the roots of its host plants. They are also migratory endoparasites, meaning they penetrate roots and are able to move throughout the host. 
In the population of  \textit{Radopholus similis}, infection of the root tissues is caused by the mobile juveniles and adult females, whereas the males, lacking a functional mouthpart, are not infectious.
 They are able to reproduce both sexually and asexually and therefore exist in female, hermaphroditic, and male forms. At all developmental stages, this nematode penetrates living, young, and tender roots. It lives and reproduces within root tissues, causing lesions and necrosis. Over time, the root progressively deteriorates and  when it no longer provides suitable living and feeding conditions, the nematode leaves and moves through the soil in search of a new root to penetrate. They often inhabit the parenchyma, and females lay about 4-5 eggs per day within the tissues they occupy, over a period of two weeks. Once laid, eggs usually take about 5-10 days to hatch, 10-13 days to develop into adults, and approximately 2 days to become gravid. This corresponds to a life cycle of 20-25 days from egg to gravid adult. See~\cite{fotso2022}  and~\cite{a1}    and the references therein for more details on the biology of \textit{R. Similis}.

\subsection{Model formulation}
In this section, based on the life cycle of the nematode \textit{Radopholus Similis}, we suggest an age-structured model describing
 the interaction between banana-plantain plants and the \textit{Radopholus similis}. The model takes into account the age of banana-plantain plant,
 free nematode in the soil and the infecting nematode in the roots of banana-plantain plants. We use the compartmental modeling approach with a progression of
stages. The banana-plantain plants in the field are subdivided into two compartments according to their infestation status. We consider the compartment of healthy (\textit{resp.} infected) banana-plantain plants. We denote by $S(a,t)$ (resp. $I(a,t)$) the density of healthy (\textit{resp.} infected) banana-plantain plants at time $t$ and of age $a$.
The \textit{Radopholus similis} population is also subdivided into two compartments at time $t$: the compartment of free nematodes in the soil, whose number is denoted by $N_{F}(t)$, and the compartment of infesting nematodes within the roots, whose number is denoted by $N_{I}(t)$.
 The mathematical formulation of the model is based on the following biological assumptions:
 \begin{enumerate}
  \item[(A1)]  We assume that the healthy banana-plantain are produced at time dependent rate denoted $m(t)$.
  \item[(A2)] The infection is initiated by free nematodes in the soil and the force of infection is modeled by $\dsp \beta(a)\mathcal{W}(P,N_{F})$, where
  $\beta(.)$ represents the age infection rate per unit of time and the interaction between free nematodes and banana-plantain plants is given by $\mathcal{W}$ which is defined
  by
 $$ \mathcal{W}(P,N_{F})=\frac{N_{F}}{P} \quad \mbox{with}\quad    P(t)=\int_{0}^{a_{\dagger}}
 (S(a,t)+I(a,t))da\,,$$
 where $a_{\dagger}$ is the maximum age of an individual which in our case is the banana-plantain plant. After infection, we have immediately two different losses in the compartments.
 \begin{itemize}
  \item[i)]   The lost of healthy banana plantain plants that become infected plants modeled by~$-\beta(a)\mathcal{W}(P,N_{F})S(a,t)\,.$
  \item[ii)]   The lost in the compartment of free nematodes in the soil,  becoming infecting nematodes at the rate 
     $\dsp -\alpha \mathcal{W}(P,N_{F}) \int_0^{a_{\dagger}}\beta(a) S(a, t)d\,a\,$,  $\alpha$ being 
  the conversion rate of the biomass of newly infected banana-plantain roots into infecting nematodes per unit of plant.   Their natural mortality is denoted by $\mu_F$. 
   \end{itemize}
 \item[(A3)] Infecting pests $ N_{I}  $ feed on plant roots with a Holling-type-II functional response of the form $\dsp  -  \frac{d(a)N_I\, I(a,t)}{K_d+\int_0^{a_{\dagger}}I(a, t)d\,a} \, $ in the compartment of infected plants. Their natural mortality rate,  $\mu_{I}$  is different from the mortality rate of free nematodes in the soil because the environments are different and the root, which serves as both host and food for the nematode, is more favorable to its survival than the soil: $\mu_{I} < \mu_{F}$ . Here, $K_d$ is known as the half saturation constant and $d$ the consumption rate.
 \item[(A4)] Infecting pests reproduce inside the roots following the logistic function defined by \break  $\dsp \rho \frac{\int_0^{a_{\dagger}}d(a )I(a, t)d\,a\,}{K_d+\int_0^{a_{\dagger}}I(a, t)d\,a\,}N_{I}\left(1- \frac{N_I}{K}\right)$, where $\rho$ and $K$ are respectively the conversion rate of ingested roots into pests and the maximum reception capacity of the environment.
 \item[(A5)] When the plants (roots) no longer present favorable living and feeding conditions, infected plants release nematodes at the rate $\gamma$.
 \item[(A6)] Free nematodes also infest the infected plants (reinfection) at rate $\beta(a)$  and immediately, we have the lost of free nematodes that become infecting nematodes denoted by $\dsp -e\alpha \int_0^{a_{\dagger}}\beta(a)\mathcal{W}(P,N_{F})I(a, t)d\,a\, $, where $0< e < 1$ is the probability of reinfection.
 \item[(A7)] All plants undergo natural mortality at rate $\mu(a)$. 
  \end{enumerate}
Based on these biological assumptions, Figure~\ref{a} shows the flowchart representing the interactions between the different variables.
		\begin{small}
\begin{center}
	\begin{figure}[H]
		\begin{center}
			\begin{tikzpicture}[x=0.85cm, y=0.85cm]
				\node[very thick, black, circle,fill=blue!30, inner sep=0.1cm, draw]  (a) at (-3,0){$S(a,t)$};
				
				\node[very thick, black, circle,fill=red!30, inner sep=0.1cm, draw] (b) at (3,0){$I(a,t)$};
				
				\node[very thick, black, rectangle,fill=blue!30, inner sep=0.2cm, draw]  (d) at (-1,-4){$N_{F}(t)$};
				
				\node[very thick, black, rectangle,fill=red!30, inner sep=0.2cm, draw] (e) at (6,-4){$N_{I}(t)$};
				
				\draw[->, very thick] (-5.5,0)--node[above]{$m(t)$}(a);
				
				\draw[ ->, dashed, very thick] (d)--node[pos=0.4, above]{}(-1,0.01);    
				
				\draw[ ->, dashed, very thick] (3,-4)--node[pos=0.4, above]{$\qquad (iii)$}(b); 
				
				\draw[ <-, dashed, very thick] 
			    (d) to[left=15] node[pos=0.4, above]{$\gamma$}(b);
				
				\draw[ ->, very thick] (a)--node[above]{(i)}(b);
				\draw[ ->, very thick] (d)--node[below]{(iv)}(e);
			
	            \draw[->, very thick, black, font=\large] (5.5,-4.5) .. controls +(0,-1) and +(0,-1) .. node[below] {($v$)} (6.5,-4.5); 
				
				\draw[->, very thick] (a)--node[right]{$\mu(a)$}(-3,1.5);
				
				\draw[->, very thick] (b)--node[left]{$\mu(a)$}(3,1.5);
				
				\draw[->, very thick] (b)--node[above]{(ii)}(5,0);
				
				\draw[->, very thick] (d)--node[left]{$\mu_{F}$}(-1,-5.5);
				
				\draw[->, very thick] (e)--node[above]{$\mu_{I}$}(8,-4);
			\end{tikzpicture}
			\caption{
			Diagram of the nematode propagation in the banana-plantain plantation. State variables: Healthy plants (S), Infected plants (I), Free nematodes $(N_{F})$ and Infesting nematodes $(N_{I}).$
			The logistic expressions of the model are given by:
						 $(i):\, \beta(a)\mathcal{W}(P,N_{F}), $
						  $(ii):\, \dfrac{d(a) N_{I}}{K_d+\int_a^{a_{\dagger}} I(a,t)da},$
						  $(iii):\, e\int_{0}^{a_{\dagger}}\beta(a)\mathcal{W}(P,N_F)da,$
					  $(iv):\,\alpha( \int_{0}^{a_{\dagger}}\beta(\cdot) S(\cdot\,,\,t)da + e \int_{0}^{a_{\dagger}}\beta(a)I(a,t))da) \mathcal{W}(P,N_{F}) $ and 
					 $(v) :\, \rho N_{I}\dfrac{\int_0^{a_{\dagger}} d(a)I(a,t)da}{K_d+\int_0^{a_{\dagger}}
					  I(a,t)da}\left(1-\dfrac{N_I}{K}\right).$
				}
						\label{a}
					\end{center}
				\end{figure}
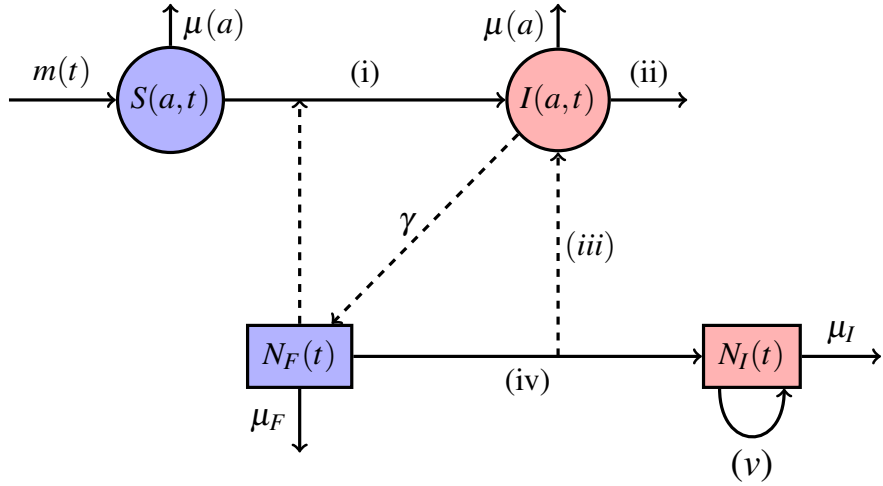
			\end{center}
 \end{small}
\begin{remark}\label{space}
Based on the definitions of $S$ and $I$, the integrals
$\displaystyle\int_0^{a_{\dagger}} S(a,t)\,da$ and
$\displaystyle\int_0^{a_{\dagger}} I(a,t)\,da$ represent, respectively, the total number of healthy and infected banana-plantain plants at time~$t$. These quantities are denoted by $\|S(\cdot, t)\|_{L^1}$ and $\|I(\cdot, t)\|_{L^1}$, which is consistent with the definition of the $L^1$-norm, since both $S$ and $I$ remain nonnegative, as will be shown later.
\end{remark}
From the diagram in Figure\ref{a}, we derive the following system of differential equations:
\begin{small}
  \begin{eqnarray}\label{Model diagram}
 \partial_{t}S(a,t)+\partial_{a}S(a,t) & = &  -\beta(a)\mathcal{W}(P,N_{F})S(a,t)-\mu(a)S(a,t) \nonumber \\
 \partial_{t}I(a,t)+\partial_{a}I(a,t) & = & \beta(a)\mathcal{W}(P,N_{F})S(a,t)-d(a)\dfrac{N_{I}(t)I(a,t)}{K_d+\lVert I(\cdot\,,\,t)\lVert}-\mu(a)I(a,t)\qquad \qquad\qquad \qquad \\
\dot{N_F}  & = & -\alpha(\lVert\beta(\cdot) S(\cdot\,,\,t)\lVert + e \lVert\beta(\cdot)I(\cdot\,,\,t))\lVert) \mathcal{W}(P,N_{F})+ \gamma\lVert I(\cdot\,,\,t) \lVert -\mu_{F}N_{F}  \nonumber  \\
\dot{N_I}  & = &  \alpha(\lVert\beta(\cdot) S(\cdot\,,\,t)\lVert + e \lVert\beta(\cdot)I(\cdot\,,\,t))\lVert) \mathcal{W}(P,N_{F}) + \rho N_{I}\dfrac{\lVert d(\cdot)I(\cdot\,,\,t)\lVert}{K_d+\lVert I(\cdot\,,\,t)\lVert}\left(1-\dfrac{N_I}{K}\right) -\mu_{I}N_{I}\nonumber
\end{eqnarray}
\end{small}
System \eqref{Model diagram} is supplemented with the following
 initial and boundary   conditions:
\begin{equation}\label{InCond}
 S(a,0)=S_{0}(a),\,\,\, 
I(a,0)=I_{0}(a), \,\,\, N_{F}(0)=N_{F0},\,\,\, N_{I}(0)=N_{i0}\,,
\end{equation}
and
\begin{equation}
S(0,t)=m(t),\,\,\,I(0,t)=0\,.
\end{equation}
\noindent
The cumulative production over the time period $[0, t]$ of healthy banana and plantain bunches, denoted by $P_c(t)$, and the daily production $P_d(D)$ at day $D$ are respectively given by:
\begin{equation}\label{cumprod}
P_c(t) = \int_{0}^{t} \int_{a^*}^{a_{\dagger}} \theta(a) S(a, u) \, da \, du, 
\qquad \text{where} \qquad 
\theta(a) = \frac{a}{a + a_0} \, \theta_{\max}\,,
\end{equation}
and
\begin{equation}
P_d(D) = P_c(D) - P_c(D - 1),
\qquad \text{for all } D \geq 1 \text{ and } t \geq 0.
\end{equation}
In~\eqref{cumprod}, $\theta_{\max}$ denotes the maximum weight of healthy banana or plantain bunches at the maturity age $a_{\dagger}$, while $a^*$ and $a_0$ are respectively the ages at which the bunch becomes mature enough to be harvested and the age at which the bunch reaches half of its maximum weight.  
For the purpose of our simulations, we choose the following parameter values:
\[
\theta_{\max} = 35\,\text{kg}, \qquad 
a^* = 240\,\text{days}, \qquad 
a_0 = 270\,\text{days}, \qquad 
a_{\dagger} = 300\,\text{days},
\]
since the cropping season of banana varies from 10 to 12 months~\cite{a3}, and we take the maximum $T = 550\,\text{days}$.
The variables and the parameters of the model and their biological
significance are summarized in Table \ref{t1}  and Table
\ref{t2} respectively.
\begin{table}[H]
 \caption{Variables of system \ref{Model diagram}.}
\begin{center}
\begin{tabular}{lll}
\hline \hline
\textbf{Variable} & \textbf{Definition} & \textbf{Units}\\
\hline \hline
$ S(a,t) $ & Healthy banana-plantain plants & Plants\\
& at time t and of age a .& \\
$ I(a,t) $ &Infected banana-plantain plants & Plants\\
& at time t and of age a. & \\
$ N_{F}(t) $ & Free nematodes  at time t & Nematodes.\\
$ N_{I}(t) $ & Infecting nematodes  at time t & Nematodes .\\
\hline
\end{tabular}
\label{t1}
\end{center}
\end{table}
\begin{center}
\begin{table}[H]
\caption{Biological significance of parameters}
\begin{center}
\begin{tabular}{lll}
\hline \hline
\textbf{Symbols} & \textbf{Meaning} & \textbf{Units}\\
\hline \hline
$m(t)$ & New healthy plants & $\text{plants}.\text{day}^{-1}$\\
$\beta(a)$ & Infection rate of $N_F$ & $\text{plants}.\text{nematode}^{-1}.\text{day}^{-1}$\\
$\alpha$ & Conversion of the number & $\text{plant}^{-1}$\\ 
         & of nematodes per plant unit & \\
$\mu(a)$ & Mortality rate of plant & $day^{-1}$   \\
$d(a)$  &  Consumption rate & $\text{plants}.\text{nematode}^{-1}.\text{day}^{-1}$ \\
$K_d$  & Half-saturation constant & plants \\
$\gamma$ & Production rate of free nematodes  & $ \text{nematodes}.\text{plant}^{-1}$\\
$e$ & Probability of reinfection  & \\
$\rho$ & Conversion rate of ingested roots & $ \text{nematodes}.\text{plant}^{-1}$\\
$K$ & Maximum capacity  & nematodes\\
$\mu_{F}$ & Mortality rate of $N_F$ & $day^{-1}$\\
$\mu_{I}$ & Mortality rate of $N_I$ & $day^{-1}$\\
$\theta$ &  Production rate &  $kg.\text{plant}^{-1}$ \\
 \hline
\end{tabular}
\label{t2}
\end{center}
\end{table}
\end{center}

\section{Mathematical analysis of the model} 

In this section, we carry out the mathematical analysis of the model system \eqref{Model diagram}, provided with its boundary and initial conditions.
\subsection{Functions spaces}
The mathematical analysis of the model system will be carried out within the following functional spaces. They are motivated by Remark\,\ref{space}. 
  Let $\Omega$ be an interval of the set of
real numbers $\mathbb{R}$. We denote by:
\begin{itemize}
\item $\dsp C^k(\Omega)\,,\; k\in\N\,,$ the set  of real-valued functions  which $k$-times continuously differentiable on  $\Omega$.
\item $\dsp C_b^k(\Omega)\,,\; k\in\N\,,$  the set of  functions in  $\dsp C^k(\Omega)$ which are bounded on $\Omega$ together with all their derivatives up to order~$k$ .
\item  $L_{+}^{\infty}(\Omega)$ the
set of measurable, essentially bounded and nonnegative functions
defined on $\Omega$ 
\item $L_{+}^{1}(\Omega) $   the
set of measurable functions that are nonnegative and
Lebesgue-integrable   over the set
$\Omega$. 
\item $W^{k,1}(\Omega)\,,\; k\in\N\,,$ denotes the Sobolev space of functions of
$L^{1}(\Omega)$ which are such that all their derivatives (weak) up to order $k$ are functions of  $L^{1}(\Omega)$. 
\end{itemize}

\subsection{Assumptions}\label{tata} \noindent
We make the following assumptions, which are motivated by biology.
We assume that:
\begin{enumerate}
\item[(A1)]  The parameters $K_d$, $\alpha$, $e$, $\gamma$, $\rho$, $K$, $b$, $\mu_{F}$, $\mu_{I}$ and the initial conditions $N_{F0}$ and $N_{I0}$ are nonnegative.
\item[(A2)]   $m \in L_{+}^{\infty}(0,\infty)\bigcap C^1([0,\infty))\,$ and $\beta \in
L_{+}^{\infty}(0,a_{\dagger})$
\item[(A3)] $\frac{S_0}{\pi}\,;\, \frac{I_0}{\pi}\;\in\;
W^{1,1}(0,a_{\dagger})\;,$ with $\,S_0(a), \; I_0(a) \geq 0\, $ a.e.
on $[0, a_{\dagger}]$.
\item[(A4)] There exists a constant $b> 0$ such that $P(t)\geq b, \, \forall t\geq 0$.
\item[(A5)]  $\mu\in  C([0,a_{\dagger}))$ and there exists a real number $\tilde{\mu}>0$
such that $\mu(a)\geqslant \tilde{\mu}$ for almost every a
$\in[0,a_{\dagger})$; $a_{\dagger}$ being the maximum age of an
individual which in our case is the banana-plantain plant. It is part of our hypotheses that $\mu $ is an increasing function   of the banana-plantain  age. 
\item[(A6)] $\;\dsp \int_{0}^{a_{\dagger}}\mu(a)da=+\infty\;$ and  $\;\dsp d  \in L_{+}^{\infty}(0,a_{\dagger})\,$.
\end{enumerate}
\begin{remark}\label{rmq2}
\begin{enumerate}
\item In the third assumption, $\pi$ is the
survival probability which will be shortly defined (see Equation \eqref{pi} bellow).
\item 
The sixth assumption is required for the survival probability \(\pi\) (defined below) to vanish at age \(a_{\dagger}\). 
An example of a function satisfying this condition is
\[
\mu(a) = \frac{\alpha_0}{(a_{\dagger} - a)^{m}}
\;\ge\;
\frac{\alpha_0}{a_{\dagger}^{m}},
\]
where \(\alpha_0 > 0\) and \(m > 0\). See, for example, \cite[p.~13]{a2}.

\item  Recall that the function $\mathcal{W}$ defined as $\displaystyle
\mathcal{W}(P,
  N_F)=\frac{N_F}{P}$ on the set $[b,\infty)\times[0,\infty)$ is of $C^1-$differentiability class and then is a locally
 Lipschitz-continuous function in both arguments i.e there exists a positive real number M such that:
\begin{equation}\label{LL}
\arrowvert\mathcal{W}(P^{1},N_{F}^{1})-\mathcal{W}(P^{2},N_{F}^{2})\arrowvert
\leq M(\arrowvert P^{1}-P^{2}  \arrowvert +\arrowvert
N_{F}^{1}-N_{F}^{2}\arrowvert)\;, \end{equation}
for all $(P^{1},N_{F}^{1}), \,(P^{2},N_{F}^{2})$ in any
 bounded subset of  $[b,\infty)\times[0,\infty)$.
\end{enumerate}
\end{remark}

\subsection{Existence and uniqueness of a mild solution}\label{se1}
In this section, we aim at proving existence and uniqueness of a mild
solution of the model system \eqref{Model diagram}. The preliminary step
consists in finding an equivalent system to system \eqref{Model diagram} with zero boundary condition. Using the ideas of Laroche
and Perasso (\cite{a3},\, 2015), we introduce the probability for a
healthy banana-plantain plant to survive to age $a$,
  defined by 
 \begin{align}\label{pi}
 \pi(a)=exp \left(-\displaystyle\int_{0}^{a}\mu(s)ds\right) \,,
 \end{align}  
    and for all $t \geq 0$, we define a new
  function $\phi$ as: 
 \begin{align}\label{rrr}
    \phi(t, \cdot ):a \mapsto   \phi(a,t): & = \left\{
    \begin{array}{ll}
   \pi(a)m(t-a) \hspace{1em} \mbox{if} \hspace{1em}  0\leq a \leq t,   \\
   0 \hspace{1em} \mbox{otherwise}.
  \end{array}
   \right.
 \end{align}
 We have the following
 \begin{propo}
 \label{prop_pi}
 From assumptions in Section \ref{tata}, functions
 $\pi$ and $\phi$ have the following properties:
\begin{itemize}
 \item   $\pi\in L_{+}^{1}(0,a_{\dagger})\cap C^1(0,a_{\dagger})$ and satisfies the differential equation
 $$  \pi\,'(a) \,=\, -\mu(a)\pi(a) \quad \mbox{for  every }\quad   a \in
 [0,a_{\dagger}] \;.$$
\item 
  For all $ a  \in [0,a_{\dagger}]$ and $t\geq 0$ such that $0\leq t  \leq
    a,$ 
    \begin{equation}\label{bn}
    0\leq \pi(a) \leq \exp(-\tilde{\mu}a)\quad \mbox{and}\quad 
   \dfrac{\pi(a)}{\pi(a-t)} \leq \pi(t)  
   \,.
   \end{equation} 
 \item  The function t $\mapsto \phi(\cdot\,,\,t)$ belongs to the space $\mathcal{C}(\mathbb{R_{+}},\,L^{1}_{+}(0,a_{\dagger}))$, the set of continuous functions defined on $\mathbb{R_{+}}$ with values in $L^{1}_{+}(0,a_{\dagger})$ .
\item 
  For all $t  \geqslant$ 0, $\phi(0,t)= m(t)$ and $\phi$ satisfies the following partial differential equation:
 \begin{align}
  \partial_{t}\phi(a,t)+\partial_{a}\phi(a,t)~=~-\mu(a)\phi(a,t).
  \end{align}
\end{itemize}
 \end{propo}
 \begin{proof}
  Assumption $(A2)$ ensures that $\phi$ is continuously differentiable and from assumption $(A4)$
   it is easy to see that $0 \leq \pi\leq 1 $ and then is integrable on $(0,
   a_{\dagger})$.
   By hypothesis $\mu(a)\geq \tilde{\mu}$ for all $a \in [0;a_{\dagger}]$, thus
  $0 \leq \pi(a)\leq \exp (-\tilde{\mu}a) $.
  Let us prove that $\frac{\pi(a)}{\pi(a-t)} \leq \pi(t)$.
  From the definition of $\pi$, we have
  $
  \frac{\pi(a)}{\pi(a-t)}= exp \left(\int_{a}^{a-t}\mu(x)dx\right)
$
 and since the function  $ t \mapsto e^{t} $ is increasing, it is sufficient to show that
$\int_{a-t}^{a}\mu(x)dx \geq \int_{0}^{t}\mu(x)dx\, $ for $\,0\leq
a\leq t$.
  Using the change of variable $ y = x-(a-t) $, we obtain
 $$
   \int_{a-t}^{a}\mu(x)dx  =  \int_{0}^{t}\mu(y+a-t)dy
 \geq  \int_{0}^{t}\mu(y)dy\,,
$$
 since $a - t \geq 0 $ and  $\mu$, the mortality rate is an increasing function   of the banana-plantain
 age.  Let  $t >0,$ for $\,0< a< t$, straightforward computations give
   $$
   \partial_{t}\phi(a,t)+ \partial_{a}\phi(a,t)  = -\mu (a)\phi(a,t)\;.
 $$
   Let us prove that $ t\longmapsto \phi(\cdot\,,\,t) \in \mathcal{C}(\mathbb{R}_{+}, \displaystyle L_{+}^{1}(0, a_{\dagger})) .$
   For this, let $ t_{n} > 0 $ such that $ \lim\limits_{n \rightarrow \infty} t_{n} =
   t_{0}.$
   $ (\phi(\cdot\,,\,t^n))_{n \in \mathbb{N}}$ is a sequence of measurable functions and for all $a\in(0,a_{\dagger}), $
   $\lim\limits_{n \rightarrow \infty} \phi(a\,,\,t_{n}) =\phi(a\,,\,t_{0})
   $. On the other hand, since $m\in L_{+}^{\infty} (0;+\infty)\cap
   C([0,\infty))$,
   $ \arrowvert \phi(a\,,\, t_{n}) - \phi(a\,,\, t_{0})\arrowvert \leq 2\Arrowvert m \Arrowvert_{\infty}\pi(a) $. Note that the mapping
    $ a \mapsto  2\Arrowvert m \Arrowvert_{\infty}\pi(a)$ belongs to $ L_{+}^{1}(0,
    a_{\dagger})$ and applying the
    dominated convergence theorem the result follows.
 \end{proof}
 In order to nullify the boundary conditions and to get rid of the terms $\mu S$ and $\mu I$ in the model system \eqref{Model diagram}, we introduce now new
 state variables  (see \cite{a2}, and \cite{a3}):
 \begin{equation}\label{nfonct}\hat{S}(a,t):=\frac{S(a,t)}{\pi(a)}-\phi(a,t) \quad \mbox{and}\quad  \hat{I}(a,t):=\frac{I(a,t)}{\pi(a)}\;.
 \end{equation}
  System \eqref{Model diagram} is then
 equivalent to the following system with zero boundary conditions:
 \begin{small}
  \begin{eqnarray}
  \label{K2}
  \partial_{t}\hat{S}(a,t)+\partial_{a}\hat{S}(a,t)
  &=&-\beta(a)\mathcal{W}(\hat{P},N_{F})[\hat{S}(a,t)+\phi(a,t)]+\mu(a)\phi(a,t)\,, \nonumber \\
  \partial_{t}\hat{I}(a,t)+\partial_{a}\hat{I}(a,t) &=&\beta(a)\mathcal{W}(\hat{P},N_{F})[\hat{S}(a,t)+\phi(a,t)] -d(a)N_{I}\dfrac{\hat{I}(a,t)}{K_d+\left\lVert \pi \hat{I}(\cdot\,,\,t)\right\lVert}\,, \nonumber\\
 \dot{N_F} &=&-\alpha \left( \left\|\beta\pi\left( \hat{S}(\cdot\,,\,t)+\phi(\cdot\,,\,t)\right)\right\lVert +e \left\lVert \beta\pi \hat{I}(\cdot\,,\,t)\right\lVert\right)\mathcal{W}(\hat{P},N_{F})+ \gamma\lVert \pi \hat{I}(\cdot\,,\,t) \lVert-\mu_{F}N_{F} \,,\qquad \qquad  \qquad \qquad \qquad  \\
 \dot{N_I} &=& \alpha \left( \left\|\beta\pi\left( \hat{S}(\cdot\,,\,t)+\phi(\cdot\,,\,t)\right)\right\lVert +e \left\lVert \beta\pi \hat{I}(\cdot\,,\,t)\right\lVert\right)\mathcal{W}(\hat{P},N_{F})-\mu_{I}N_{I}\,,\\ 
 & & + \rho N_{I}\dfrac{\lVert d \pi I(\cdot\,,\,t)\lVert}{K_d+\lVert \pi I(\cdot\,,\,t)\lVert}\left(1-\dfrac{N_I}{K}\right)\,,
 \nonumber
 \end{eqnarray}
  \end{small}
 where, $$\hat{P}(t)= \int_{0}^{a_{\dagger}} \pi(a)\left(\hat{S}(a,t)+\hat{I}(a,t)+\phi(a,t)\right)da.$$
The initial and the boundary conditions now read:
\begin{equation}\label{eq1}
\hat{S}(a,0)=\frac{S_{0}(a)}{\pi(a)}\,,\;\hat{I}(a,0)=\frac{I_{0}(a)}{\pi(a)}\,,\;N_{F}(0)=N_{F0}\,,\;N_{I}(0)=N_{i0}\,,
\end{equation}
and 
\begin{equation}\label{eq2}
 \hat{S}(0,t)=  \hat{I}(0,t)=0\,.
\end{equation}
It is straightforward to verify that the vector of state variables
\(\dsp 
(\hat{S}(a,t), \hat{I}(a,t), N_{F}(t), N_{I}(t))
\)
satisfies the model system~\eqref{K2} if and only if the original vector of state variables \break
\(\dsp 
(S(a,t), I(a,t), N_{F}(t), N_{I}(t))
\)
satisfies the initial model~\eqref{Model diagram}.
We shall therefore analyze system~\eqref{K2}, which can be formulated as a standard abstract Cauchy problem and is more convenient to study.

Let $\mathcal{A}_{1}$ and $\mathcal{A}_{2}$ be the
operators defined by: 
\begin{align*}
\mathcal{A}_{1}: D(\mathcal{A}_{1})\subset L^{1}([0,a_{\dagger}))&\rightarrow L^{1}([0,a_{\dagger}))\\
\varphi& \mapsto\mathcal{A}_{1}\varphi\,\,\, such \,\,\, that \,\,\,
\mathcal{A}_{1}\varphi(a)=-\varphi'(a)  \,\,\, for \,\,\, all \,\,\,
a \in [0,a_{\dagger}]
\end{align*}
and
 \begin{align*}
 \mathcal{A}_{2}: D(\mathcal{A}_{2})\subset \mathbb{R}^{2}&\rightarrow\mathbb{R}^{2} \\
 (\varphi_{1},\varphi_{2})^{T}& \mapsto\mathcal{A}_{2}(\varphi_{1},\varphi_{2})^{T}
 = (-\mu_{F}\varphi_{1},-\mu_{I} \varphi_{2})^{T},
 \end{align*}
with the subspaces $D(\mathcal{A}_{1})$ and $D(\mathcal{A}_{2})$
defined  by (see  \cite{a7} and \cite{a3})
$$D(\mathcal{A}_{1}):=\left\{\varphi\in
W^{1,1}(0,a_{\dagger}):\;\varphi(0)=0 \right\} \,\,\ \mbox{and}\,\,\
D(\mathcal{A}_{2})=\mathbb{R}^{2}\;. $$ As it is well known in the literature, the
linear operators $\mathcal{A}_{1}$ and
$\mathcal{A}_{2}$ represent the mortality process related to the
healthy banana-plantain plants  and the nematodes respectively.

Consider the product space
$ \mathbb{X}= L^{1}([0,a_{\dagger}))\times
L^{1}([0,a_{\dagger}))\times \mathbb{R}\times  \mathbb{R}\,.
$  Endowed with the product norm
$\displaystyle
 \Arrowvert x \Arrowvert_{\mathbb{X}}=\Arrowvert x_{1}  \Arrowvert_{L_{1}}+\Arrowvert x_{2} \Arrowvert_{L_{1}}+\lvert x_{3} \lvert+\lvert x_{4} \lvert 
$  
 ( for   $x =(x_{1} ,x_{2} ,x_{3} ,x_{4} ))$, \,  $ \mathbb{X}$ is a Banach space.
 The non-negative cone of the Banach space $\mathbb{X}$ is defined by
 \begin{equation*}
 \mathbb{X}_{+}= L^{1}_{+}([0,a_{\dagger}))\times L^{1}_{+}([0,a_{\dagger}))\times \mathbb{R}_{+} \times  \mathbb{R}_{+}.
 \end{equation*}
Set
$u(t)=(\hat{S}(\cdot\,,\,t),\hat{I}(\cdot\,,\,t),N_{F}(t),N_{I}(t))$,
the solution of the system \eqref{K2} corresponding to the initial
condition
$u(0)=(\hat{S}(\cdot, 0),\hat{I}(\cdot, 0),N_{F}(0),N_{I}(0))$,
and consider  the linear operator
\begin{equation*}
 \mathcal{A}:D(\mathcal{A}) \subset \mathbb{X} \rightarrow \mathbb{X}\,\,\
 \mbox{where}\,\,\ D(\mathcal{A})= D(\mathcal{A}_{1}) \times D(\mathcal{A}_{1})\times D(\mathcal{A}_{2}) \,\,\ \mbox{ such that}\,\,\
 \end{equation*}
 \begin{equation}
\mathcal{A}:= \mbox{diag}(\mathcal{A}_{1},\mathcal{A}_{1},\mathcal{A}_{2}).
 \end{equation}
 Define now the nonlinear perturbation map  
 \begin{eqnarray*}
  \mathcal{H}:\mathbb{R_{+}} \times D(\mathcal{A}) \subset \mathbb{R_{+}}\times\mathbb{X}&\rightarrow&\mathbb{X} \\
  (t, u(t))&\mapsto &\mathcal{H}(u(t)), \,\,\, such \,\,\, that \,\,\, for \,\,\, all \,\,\, a \,\,\, \in [0, a_{\dagger}],
  \end{eqnarray*}
  \begin{small}
  \begin{align}
  \label{nlt}
   \mathcal{H}(u(t))(a) & =
   \begin{pmatrix}
   -\beta(a)\mathcal{W}(\hat{P},N_{F})[\hat{S}(a,t)+\phi(a,t)]+\mu(a)\phi(a,t)\\
    \beta(a)\mathcal{W}(\hat{P},N_{F})[\hat{S}(a,t)+\phi(a,t)]-d(a)N_{I}\dfrac{\hat{I}(a,t)}{K_d+\left\lVert \pi \hat{I}(\cdot\,,\,t)\right\lVert}\\
  -\alpha \left( \left\|\beta\pi\left( \hat{S}(\cdot\,,\,t)+\phi(\cdot\,,\,t)\right)\right\lVert +e \left\lVert \beta\pi \hat{I}(\cdot\,,\,t)\right\lVert\right)\mathcal{W}(\hat{P},N_{F})+ \gamma\lVert \pi \hat{I}(\cdot\,,\,t) \lVert \\
    \alpha \left( \left\|\beta\pi\left( \hat{S}(\cdot\,,\,t)+\phi(\cdot\,,\,t)\right)\right\lVert +e \left\lVert \beta\pi \hat{I}(\cdot\,,\,t)\right\lVert\right)\mathcal{W}(\hat{P},N_{F})+\dfrac{\rho N_{I}\lVert d \pi I(\cdot\,,\,t)\lVert}{K_d+\lVert \pi I(\cdot\,,\,t)\lVert}\left(1-\dfrac{N_I}{K}\right)  
  \end{pmatrix}
  \,.
 \end{align}
 \end{small}
Then mapping $\mathcal{H}$ is well defined. To prove that, we only need to show that $\mu\phi \in L^1(0, a_{\dagger})$ for any fixed $t > 0$. This follows from the observation 
\begin{align*}
\int_{0}^{a_{\dagger}}\mu(a) \phi(a) da & = \int_{0}^{t}\mu(a)\pi(a) m(t-a) da\,,\\ 
& \leq \|m\|_{\infty}\int_{0}^{t}\mu(a)\exp{\left(-\int_{0}^{a}\mu(s)ds\right)}da \,,\\
& \leq \|m\|_{\infty}(1-\pi(t)) < \infty \,.
\end{align*}
Then, using these operators, the model system \eqref{K2} with the
boundary and initial conditions take the following form of  abstract Cauchy
problem on the Banach space $\mathbb{X}$:
 \begin{equation}\label{K7}
 \left\{
  \begin{array}{ll}
 \dot{u}(t) ~=~\mathcal{A}u(t)+\mathcal{H}(u(t)),   \qquad \mbox{for} \qquad  t \geqslant 0, \\
 u(0) ~=~ u_{0} \in D(\mathcal{A}).
  \end{array}
   \right.
 \end{equation}
Note that the property $u_0 \in D(\mathcal{A})$ follows from assumption (A3) of Section \ref{tata} (See equations \eqref{eq1} and \eqref{eq2}). We study now the semi linear Cauchy problem (\ref{K2}). The idea is to use the theory of Lipschitz perturbations of linear evolution equations of A. Pazy \cite{a4}. See also \cite{a5,a6,a7,a8} and the references therein. We will then
show that the linear operator $\mathcal{A}$ is the infinitesimal generator of a strongly continuous semigroup $ \{T(t),\,t\geq 0 \}$ and that the nonlinear part of the system $\mathcal{H} : \mathbb{R_{+}}\times\mathbb{X} \rightarrow  \mathbb{X}$ is continuous in $t\geq 0$ and locally Lipschitz in $u$, uniformly in $t$ on bounded intervals. See A. Pazy (\cite{a4}, Theorem 1.4, page 185).
 For this purpose, denote by $\rho(\mathcal{A}_{1})$ the resolvent of the linear operator $\mathcal{A}_{1}$ i.e. the set of all complex numbers $\lambda$ for which ($\lambda  - \mathcal{A}_{1}$) is invertible and $(\lambda  -\mathcal{A}_{1})^{-1}$ is a bounded linear operator on $L^{1}(0,a_{\dagger})$.
 We have the following (see Laroche and Perasso \cite{a3}; 2015):
 \begin{propo}\label{fr} The linear operator
  $(\mathcal{A}_{1},D(\mathcal{A}_{1}))$ is such that:
  \begin{enumerate}
  \item  $D(\mathcal{A}_{1})$ is dense in $L^{1}(0,a_{\dagger})$.
  \item $\mathcal{A}_{1}$ is closed .
  \item The resolvent set $\rho(\mathcal{A}_{1}) \supset
  (0,\infty)$ and
  the family of bounded linear operators \break  $R_{\zeta}$ :=$(\zeta  -\mathcal{A}_{1})^{-1}$  called resolvent is given for each $\zeta \in (0,\infty)$ by
\begin{equation}
 R_{\zeta}(\varphi)(a)=\int_{0}^{a}\exp(\zeta(\eta-a)) \varphi(\eta)d\eta \;; \;  \forall
 \varphi \in L^{1}(0,a_{\dagger}).
\end{equation}
\noindent
Moreover, the norm of $R_{\zeta}$ satisfies the estimate $\Arrowvert
R_{\zeta} \Arrowvert \leqslant \dfrac{1}{\zeta}, \, \forall \zeta
\in (0,\infty).$
\end{enumerate}
\end{propo}
\begin{remark}
The proof of this proposition has to be compared with that of Lemma 2 of [\cite{a7}, 2022] where the authors consider that the function that defines the age-specific natural mortality rate in the compartment of healthy is an essentially bounded function. This less realistic
assumption is replaced in this work by assumptions (A5) and (A6) in Section \ref{tata}. See also Remark \ref{rmq2}.
\end{remark}

\proof
  Set
$E:=\left\{\varphi\in W^{1,1}(0,a_{\dagger}),\,\,\, supp(\varphi)
\,\, \,\, compact \,\,\, in \,\,\, [0, a_{\dagger}) \right\}$. Note
that, if $\Omega$ is an open subset of $\R^N, \, N\geq 1$ then the
set $C^{\infty}_c(\Omega)$ of smooth functions compactly supported
in $\Omega$ is dense in the space $L^p(\Omega),\, 1\leq p < +\infty$
 (see for example [\cite{a9}, 2011], Corollary 4.32 page 109). Since $C^{\infty}_c(0, a_{\dagger})\subset E\subset L^{1}(0,a_{\dagger})$
  it is clear that $E$ is dense in $L^{1}(0,a_{\dagger})$ and it will be sufficient to
show that $D(\mathcal{A}_{1})$ is dense in $E$ with the topology of
$L^{1}(0,a_{\dagger})$.
 Let $\phi
\in E$, we want to show that there exists a sequence of elements of $D(\mathcal{A}_{1})$ which converges to $\phi$ in $L^1(0, a_{\dagger})$.
  For   $0 < \eta < a_{\dagger}$ sufficiently small, we define a function denoted by $g_{\eta}$ by setting
 $$
               g_{\eta}(a)= \left\{
                \begin{array}{ll}
                \frac{a}{\eta}\phi(\eta) \quad\mbox{if}\quad  0\leq a\leq \eta  \\
                \phi(a) \quad\mbox{if}\quad  a\geq \eta
              \end{array}
              \;.
               \right.
$$
Obviously $g_{\eta}(0)=0, \,g_{\eta}\in L^1(0,a_{\dagger})$ and for
any $\psi\in C_c^{\infty}(0,a_{\dagger})$,
$$
\int_0^{a_{\dagger}} g_{\eta}(a)\psi'(a)da =
-\frac{1}{\eta}\phi(\eta)\int_0^{\eta}\psi(a)da
-\int_{\eta}^{a_{\dagger}}\phi'(a)\psi(a)\,,
$$
which means that the weak derivative of $g_{\eta}$ is defined as
$$ g_{\eta}'(a)=\left\{
\begin{array}
  {l}
  \frac{1}{\eta}\phi(\eta) \quad \mbox{if} \quad 0\leq a\leq \eta
  \\
  \phi'(a)\quad  \mbox{if}\quad  \eta\leq a_{\dagger}
\end{array}\,, \right.$$  which belongs to $L^1(0,a_{\dagger})$ since $\phi\in
W^{1,1}(0, a_{\dagger})$. We have thus prove that $g_{\eta}\in
W^{1,1}(0, a_{\dagger}).$
 Finally, we have
$$
  \|g_{\eta}-\phi\|_{L^1(0,a_{\dagger})}\leq \frac{1}{2}\eta|\phi(\eta)|+\int_{0}^{\eta}\left|\phi(\eta)-\phi(a)
  \right|da\,  \rightarrow  0,\; \eta   \rightarrow   0.
$$
Note that we have the embedding $W^{1,1}(0, a_{\dagger}) \hookrightarrow L^{\infty}(0, a_{\dagger})$ (see for example [\cite{a9}, 2011], Theorem 8.8 page 212), therefore $\phi$ is essential bounded on $(0, a_{\dagger})$ and the right hand side of the last inequality converges to $0$ as $\eta$ goes to $0$. This proves that $D(\mathcal{A}_{1})$ is dense in $E$ and thus in $L^{1}(0,a_{\dagger}).$

We prove now that the unbounded linear operator $\mathcal{A}_{1}$ is
closed. Let $(\vp_n)_{n\in\N}$ be a sequence of points in $
D(\mathcal{A}_{1})$ such that the sequence $(\vp_n,
\mathcal{A}_{1}\vp_n)_{n\in\N}$ converges to $(\vp,\psi)$ in $L^{1}(0,
a_{\dagger})\times L^{1}(0, a_{\dagger})$. We aim at proving that
the couple $(\vp,\psi)$ belongs to the graph of $\mathcal{A}_{1}$ i.e.
$\vp\in D(\mathcal{A}_{1}), \, \psi\in L^{1}(0, a_{\dagger})$ and $
\psi = \mathcal{A}_{1}\vp.$ By hypothesis, the sequences
$(\vp_n)_{n\in\N}$ and $(\mathcal{A}_{1}\vp_n)_{n\in\N}$ converge
respectively to $\vp$ and $\psi$ in $L^{1}(0, a_{\dagger})$ and
since $L^{1}(0, a_{\dagger})$ embeds continuously in
$\mathcal{D}'(0,a_{\dagger})$, the set of distributions on
$(0,a_{\dagger})$, these convergence also hold in
$\mathcal{D}'(0,a_{\dagger})$. Now the continuity of the derivation
operator on $\mathcal{D}'(0,a_{\dagger})$ implies that
$-\vp_n'=\mathcal{A}_{1}\vp_n \longrightarrow -\vp'$ in
$\mathcal{D}'(0,a_{\dagger})$ and since the topology of
$\mathcal{D}'(0,a_{\dagger})$ is Haussdorf, one has $-\vp'=\psi$;
but $\psi\in L^{1}(0, a_{\dagger})$ thus $\vp\in D(\mathcal{A}_{1})$
and $\psi= \mathcal{A}_{1} \vp$.

\noindent
Let $\varphi \in D(\mathcal{A}_{1})$. We have:
\begin{align*}
R_{\zeta}(\zeta  -\mathcal{A}_{1})\varphi(a) & = R_{\zeta}(\zeta\varphi+\varphi')(a)\,, \\
& =  \int_{0}^{a}\exp(\zeta(\eta-a)) (\zeta\varphi(\eta)+\varphi'(\eta)d\eta \,,\\
& = \int_{0}^{a}\exp(\zeta(\eta-a)) (\zeta\varphi(\eta) )d\eta +
\int_{0}^{a}\exp(\zeta(\eta-a)) \varphi'(\eta)d\eta \,,
\\ & =
  \varphi(a).
\end{align*}
In the same way, let $\varphi \in L^{1}(0,a_{\dagger})$. We have:
$$
(\zeta  -\mathcal{A}_{1})R_{\zeta}\varphi(a)
     = \zeta R_{\zeta}\varphi(a) + \dfrac{d}{da}(R_{\zeta}\varphi)(a)
   =  \zeta
R_{\zeta}\varphi(a) -\zeta R_{\zeta}\varphi(a)+\varphi(a)  =
\varphi(a).
$$
Finally,
\begin{align*}
\Arrowvert R_{\zeta}\varphi\Arrowvert_{L^{1}}
=\int_{0}^{a_{\dagger}}\arrowvert R_{\zeta}\varphi(a) \arrowvert da
&  \leq \int_{0}^{a_{\dagger}} \int_{0}^{a}\exp (\zeta(\eta-a))
\left|\varphi(\eta) \right| d\eta da  \,,\\
& \leq \int_{0}^{a_{\dagger}}\exp(-a \zeta )
\int_{0}^{a}\arrowvert\varphi(\eta)
\arrowvert\exp(\eta\zeta)d\eta da\,,\\
& \leq \left[-\dfrac{1}{\zeta}\exp(- \zeta
a)\int_{0}^{a}\arrowvert\varphi(\eta) \arrowvert\exp(\eta \zeta
)d\eta
\right]_{0}^{a_{\dagger}}+\dfrac{1}{\zeta}\int_{0}^{a_{\dagger}}\lvert
\varphi(a)\lvert da \,,\\
& \leq -\dfrac{1}{\zeta}\exp(- \zeta
a_{\dagger})\int_{0}^{a_{\dagger}}\arrowvert\varphi(\eta)
\arrowvert\exp(\eta \zeta
)d\eta+\dfrac{1}{\zeta}\int_{0}^{a_{\dagger}}\lvert
\varphi(\eta)\lvert d\eta\,, \\
& \leq \dfrac{1}{\zeta} \Arrowvert \varphi \Arrowvert_{L^{1}}
\,.
\end{align*}
Thus, $\displaystyle \Arrowvert R_{\zeta}\Arrowvert  \leq
\dfrac{1}{\zeta }$ and therefore, $\displaystyle
 R_{\zeta} =(\zeta  -\mathcal{A}_{1})^{-1}.$
\qed

\begin{remark}\label{rmq1}
Proposition \ref{fr} shows that the linear operator
$\mathcal{A}_{1}$ is the infinitesimal generator of $C_0$-semi group
of contraction. One can easily check that this semi group
$(T_1(t))_{t\geq 0}$ is defined by
\begin{equation}\label{tt}
T_1(t)\varphi(a)=\left\{
    \begin{array}{ll}
   \varphi(a-t) \qquad \mbox{if} \qquad  & t \leq a \leq a_{\dagger}  \\
   0 \qquad \mbox{otherwise} \qquad
  \end{array} 
   \right. , \;\; \mbox{forall} \;\varphi\in L^1(0,a_{\dagger})\,.
\end{equation}
On the other hand, the part $\mathcal{A}_{2}$ of the operator
$\mathcal{A}$is the infinitesimal generator of $C_0$-semi group of
contraction $\left(T_2(t)\right)_{t\geq 0}$ defined by
\begin{equation}
  \label{ti}
   T_{\mathcal{A}_{2}}(t)=\mbox{diag}\left(\exp(-\mu_{F}t),\exp(-\mu_{I}t)\right).
\end{equation}
\end{remark}
\noindent
We prove now that the nonlinear perturbation map
$\mathcal{H}(u)$ in the abstract system \eqref{K7} is continuous
in $t\geq 0$ and locally Lipschitz in $u$, uniformly in $t$.
  For every constant $\delta > 0$,
  denote by $\mathcal{B}_{\delta}$ the ball in the space $\mathbb{X}$ centered at $0$ and of radius $\delta$
  i.e.
 \begin{align*}
 \mathcal{B}_{\delta}:=\left\{u \in \mathbb{X}: \Arrowvert u \Arrowvert_{\mathbb{X}} \leq \delta \right\}.
 \end{align*}
 \begin{remark}\label{rreq}
 Note that for  $u \in \mathcal{B}_{\delta}$ and $P \in [b, +\infty)$ we have
 $\displaystyle\lvert\mathcal{W}(P,N_{F})\lvert< \frac{\delta}{b}=:\theta_1 $ and
 from assumption (A5) of Section \ref{tata}, on can easily see that   $\Arrowvert \phi(\cdot\,,\,t) \Arrowvert_{L_{1}} \leq \dfrac{1}{\tilde{\mu}}\Arrowvert m \Arrowvert_{\infty} =: \theta_2$.
 \end{remark}
\noindent We have the
following.
\begin{propo}\label{tup}
 The pertubation map $\mathcal{H}$ is a locally Lipschitz continuous function in $u$, uniformly in $t  \geq$ 0 i.e for all $\delta >$ 0,
 there exists a constant depending on $\delta$ denoted by $M_{\delta}$ such that
  \begin{equation*}
  \Arrowvert\mathcal{H}(u_{2})-\mathcal{H}(u_{1})\Arrowvert_{\mathbb{X}} \leq M_{\delta}\Arrowvert u_{2}-u_{1}\Arrowvert_{\mathbb{X}},
   \forall (u_{1},u_{2}) \in B_{\delta}\times B_{\delta}\,.
   \end{equation*}
   \end{propo}

\proof
In the sequel the dependence on $t$ and the space in which the norm is taken will be omitted. For example, we will write for $\Arrowvert S^{1}-S^{2}\Arrowvert$ for  $\Arrowvert S^{1}(\cdot\,,\,t)-S^{2}(\cdot\,,\,t)\Arrowvert_{L^{1}}$ and so on.
Let $\delta > 0$, for all $t \geq 0 $ consider 
$u_{1}=(S^{1},I^{1},N_{F}^{1},N_{I}^{1}) $ and 
$u_{2}=(S^{2},I^{2},N_{F}^{2},N_{I}^{2})$ two vectors of  $\mathcal{B}_{\delta}$ and denote by $\mathcal{W}^{\,i}, i=1,2$ the quantity $\mathcal{W}(P^{\,i}, N^{\,i}_F)$.
Set
\begin{align}\label{gh}
 \mathcal{H}(u_{1})-\mathcal{H}(u_{2}) & =
   \left(
   A(u),\;
   B(u) ,\;
   C(u) ,\;
   D(u)
  \right)^T
  \,,
 \end{align}
Here, \( A(u), B(u), C(u), \) and \( D(u) \) are the components of \( \mathcal{H}(u_{1}) - \mathcal{H}(u_{2}) \), as obtained from \eqref{nlt}. After straightforward manipulations, one has (   $\kappa(I^i) = K_d + \|\pi I^i\|$ ):
 \[
 A(u)= \beta\big( (\mathcal{W}^{2}-\mathcal{W}^{1})(S^{1}-\phi)+ \mathcal{W}^{2}(S^{2}-S^{1}) )\big)\,,  
 \]
\beaa
 B(u) 
  &=&
- A(u)-   \frac{d}{\kappa(I^1)\kappa(I^2)}(N_I^1 I^1 (\|\pi I^{2}\|-\|\pi I^{1}\|) + (K_d + \|\pi I^1\|)(N_I^1 (I^1 - I^2)
\\
&&
\qquad  + I^2(N_I^1 - N_I^2)))\,,
\eeaa
\beaa
C(u) &=&
\alpha\left(\mathcal{W}^2-\mathcal{W}^1\right)\int_0^{a_{\dagger}}\beta\pi\left( (S^2+\phi) +e I^{2} \right)da
+\alpha \mathcal{W}^1\int_0^{a_{\dagger}}\beta\pi\left( S^2-S^1\right)da\\
&&
+\alpha \mathcal{W}^1\int_{0}^{a_{\dagger}}\beta e\pi\left( I^{2}-I^{1} \right)da+\gamma(\|\pi I^{1}\|-\|\pi I^{2}\|)
\,,
\eeaa
\begin{small}
\[
D(u) = -C(u)+\gamma\left(\|\pi I^1\|-\|\pi I^2\| \right) +\rho N_{I}^1\dfrac{\lVert d \pi I^1\lVert}{K_d+\lVert \pi I^1\lVert}\left(1-\dfrac{N_I^1}{K}\right)
 -\rho N_{I}^2\dfrac{\lVert d \pi I^2\lVert}{K_d
 +\lVert \pi I^2\lVert}\left(1-\dfrac{N_I^2}{K}\right) 
\]
\end{small}
The last two terms of $D(u)$ can be written as 
\begin{small}
\begin{align*}
\rho N_{I}^1\dfrac{\lVert d \pi I^1\lVert}{\kappa(I^1)}\left(1-\dfrac{N_I^1}{K}\right)-\rho N_{I}^2\dfrac{\lVert d \pi I^2\lVert}{\kappa(I^2)}\left(1-\dfrac{N_I^2}{K}\right) =& \rho N_{I}^1\dfrac{\lVert d \pi I^1\lVert}{\kappa(I^1)}\left(\dfrac{N_I^2 - N_I^1}{K}\right)
 +\left(1-\dfrac{N_I^2}{K}\right)E(u)
 \,,
\end{align*}
\end{small}
where, 
\begin{align*}
E(u) = & \rho N_{I}^1\dfrac{\lVert d \pi I^1\lVert}{\kappa(I^1)} - \rho N_{I}^2\dfrac{\lVert d \pi I^2\lVert}{\kappa(I^2)} \\
= & \dfrac{\rho}{ \kappa(I^1) \kappa(I^2) }\left[(N_I^1)^2\|d\pi I^1\|(\|\pi I^2\|-\|\pi I^1\|) + (N_I^1)^2\|\pi I^1\|(\|d\pi I^1\|-\|d\pi I^2\|)\right]\\
& + \dfrac{\rho}{(K_d+\lVert \pi I^1\lVert)(K_d+\lVert \pi I^2\lVert)}K_d (N_I^1)^2(\|d\pi I^1\|-\|d\pi I^2\|)\\
& + \dfrac{\rho}{\kappa(I^1) \kappa(I^2)}(K_d +\|\pi I^1\|) \|d\pi I^1\|(N_I^1-N_I ^2)(N_I^1+N_I ^2)\,.
\end{align*}
Using these expressions, together with Inequality~\ref{LL}, Remark~\ref{rreq}, and the estimate \break 
\(\dsp 
|P^2 - P^1| \leq \|S^2 - S^1\| + \|I^1 - I^2\|,
\) 
we deduce that there exists a constant \( M_\delta > 0 \) such that  
\[
\|\mathcal{H}(u_1) - \mathcal{H}(u_2)\|_{\mathbb{X}} \leq M_\delta \|u_1 - u_2\|_{\mathbb{X}},
\]
where \( M_\delta \) is a polynomial function of the parameters  
$ 
\theta_1,$  $\theta_2,$  $\|\beta\|_{\infty}$, $\|d\|_{\infty}$,
 $e$, $\delta$, $\alpha$, $\rho$, $K$, $K_d$, $\dsp \frac{1}{K}$, $\frac{1}{K_d}$ and 
 $\frac{1}{K+K_d}$.   
This completes the proof. \qed

We now proceed to establish the existence of a unique mild solution to system~\eqref{K7}. Keeping in mind Remark~\ref{rmq1}, we begin by recalling the following.
   \begin{defi}
A mild solution of the abstract Cauchy problem \eqref{K7} is a
continuous function $u$ which satisfies the integral equation
$$
u(t) = T(t)u_0+\int_0^tT(t-s)\mathcal{H}(u(s))ds \;;
$$
where $\displaystyle T = \mbox{diag}(T_1, T_1, T_2)$ (see (\ref{tt}) and
(\ref{ti})).
   \end{defi}
   \noindent
We have the following
  \begin{theo}
   For any initial condition in $D(\mathcal{A})$, there exists a maximal interval of time $[0,t_{max})$ on which the abstract Cauchy problem \eqref{K7} has a unique mild solution.
 \end{theo}
 \begin{proof}
 The proof will be done using a fixed point argument by adapting the ideas of Pazy , (\cite{a4}; 1985), see also  Laroche and Perasso (\cite{a3}, 2015); Perasso and Laroche (\cite{a10}; 2008).
Let $\tau$ > 0 and $\hat{\delta}$ > 0 be such that
   \begin{align*}
 \hat{\delta}:= 3\Arrowvert u_{0} \Arrowvert_{\mathbb{X}}  \;, \;\;
 \tau : = \dfrac{\Arrowvert u_{0}\Arrowvert_{\mathbb{X}}}{\hat{\delta}M_{\hat{\delta}}+\|\mathcal{H}(u_0)\|_\mathbb{X}},
  \end{align*}
where the $M_{\hat{\delta}}$ is the local Lipschitz constant of the
map $\mathcal{H}$ defined in Proposition \ref{tup}. Consider the space
$\mathcal{C}([0,\tau],\mathbb{X})$ of continuous functions defined
on the time interval $[0,\tau]$ with values in $\mathbb{X}$. Endowed
with the supremum norm: $\displaystyle  \|u\|_{\infty}:=
\sup\limits_{t\in [0,\tau]}\|u(t)\|_{\mathbb{X}}$,
$\mathcal{C}([0,\tau],\mathbb{X})$ is a Banach  space. Next,
consider the ball in $\mathcal{C}([0,\tau],\mathbb{X})$ defined as
   \begin{align*}
   \tilde{\mathcal{B}}_{\hat{\delta}}:= \left\{u \in \mathcal{C}([0,\tau],\mathbb{X}) : \Arrowvert u -u_0 \Arrowvert_{\infty} \leq \hat{\delta}   \right\}\,.
   \end{align*}
   \noindent
Finally, let us consider the non linear mapping :
   \begin{align*}
   \mathcal{G}: \mathcal{C}([0,\tau],\mathbb{X}) &\rightarrow \mathcal{C}([0,\tau],\mathbb{X}) \\
   u(.)& \mapsto\mathcal{G}
(u)(t) = T_{\mathcal{A}}(t)u_{0} +
\int_{0}^{t}T_{\mathcal{A}}(t-s)\mathcal{H}(s , u(s))  ds.
\end{align*}
We prove that the mapping $\mathcal{G}$ preserves the ball $
\tilde{\mathcal{B}}_{\hat{\delta}}$. We have
\begin{eqnarray*}
  \|\mathcal{G}
(u)(t)-u_0\|_\mathbb{X}&\leq&  \|\mathcal{G} (u)(t)\|_\mathbb{X}+\|u_0\|_\mathbb{X}\,,\\
&\leq&
\|u_0\|_\mathbb{X}+\|T_{\mathcal{A}}(t)\|_{\mathcal{L}(\mathbb{X})}\|u_{0}\|_{\mathbb{X}}
+
\int_{0}^{t}\|T_{\mathcal{A}}(t-s)\|_{\mathcal{L}(\mathbb{X})}\|\mathcal{H}(u(s))\|_{\mathbb{X}}ds\,,
\\
&\leq& 2\|u_0\|_\mathbb{X}+\int_{0}^{t} \|\mathcal{H}(s ,
u(s))-\mathcal{H}(u_0)\|_{\mathbb{X}}ds+t\|\mathcal{H}(u_0)\|\,,
\\
&\leq& 2\|u_0\|_\mathbb{X}+M_{\hat{\delta}}\int_{0}^{t} \| u(s)- u_0
\|_{\mathbb{X}}ds+t\|\mathcal{H}(u_0)\|\,,\\
& \leq& 2\|u_0\|_\mathbb{X}+t\left(\hat{\delta}M_{\hat{\delta}} +
\|\mathcal{H}(u_0)\|\right)\,, \\
&\leq& 3\|u_0\|_\mathbb{X}\,.
\end{eqnarray*}
Next, we prove that $\mathcal{G}$ is a
 contraction mapping of
$ \tilde{\mathcal{B}}_{\hat{\delta}}$. We have
\begin{eqnarray*}
  \|\mathcal{G}
(u)(t)-\mathcal{G} (v)(t)\|_\mathbb{X}&\leq&
\int_{0}^{t}\|T_{\mathcal{A}}(t-s)\|_{\mathcal{L}(\mathbb{X})}\|\mathcal{H}(u(s))-\mathcal{H}(v(s))\|_{\mathbb{X}}ds\,,\\
&\leq&\tau M_{\hat{\delta}}\|u-v\|_{\mathbb{X}}\,,\\
&\leq& \frac13\|u-v\|_{\mathbb{X}}\;.
\end{eqnarray*}
The last inequality is a direct consequence of the definitions of
$\hat{\delta}$ and $\tau$. We therefore conclude using the Banach
fixed point theorem that the map $\mathcal{G}$ has a unique fixed
point  $u\in  \tilde{\mathcal{B}}_{\hat{\delta}}$. This fixed point
is thus the desired mild solution on the interval $[0, \tau]$ of the
evolution system \eqref{K7}.
   We repeat all the previous arguments,
   but now with the initial condition u($\tau$) instead $u(0)$,
   $\tau_{1}=\tau+\kappa$ with $\kappa > 0$, together with the functional space $\mathcal{C}([\tau,\tau_{1}],\mathbb{X})$ and the
   mapping $\mathcal{G}_{1}$ defines exactly as $\mathcal{G}$. The same computations now apply to $\mathcal{G}_{1}$ show
   also that $\mathcal{G}_{1}$ is again a strict contraction with Lipschitz constant $\dfrac{1}{3}$. This gives once more a unique fixed point
    which extends the previous solution on the interval [0,$\tau_{1}$]. By proceeding successively , we can extend the solution on a maximal time interval  $[0,t_{\max})$,
     so that u $\in \mathcal{C}([0,t_{\max}),\mathbb{X})$ is a mild solution of the system \eqref{K7}.

   As we have shown that the abstract Cauchy system \eqref{K7} has a unique mild solution \break  $ u \in \mathcal{C}([0,t_{\max}),\mathbb{X})$,
   then using the transformation $$u(t) \mapsto  \mbox{diag}(\pi, \pi, 1, 1)u(t)^{T}+ (\pi \phi(\cdot\,,\,t),0,0,0)=(S(a,t),I(a,t),N_{F},N_{I})\,,$$ yields the existence and uniqueness
   of a mild solution $(S, I, N_{F}$, $N_{I}) \in \mathcal{C}([0,t_{\max}],\mathbb{X})$ of the initial model \eqref{Model diagram}.
     \end{proof}

   \subsection{Positivity and boundedness of solutions}\label{se2}
 In this section we prove the positivity and the boundedness of solutions of our model. This
 is required from a biological point of view as in any theory of modeling of populations, all populations must be positive and bounded from above. As a byproduct we will obtain that our solutions are global (in time). This will be a direct consequence of the well known continuation criterion which says that, the breakdown of a solution u of a first order hyperbolic
 system of partial differential equations must involve a blow-up of the supremum $sup_{(a, t)} |u(a, t)|$. In other
 words, if the $C^0-$ norm of a solution $u(a, t)$ on an interval $[0,T)$ is uniformly bounded then this
 solution can be extended beyond T. We have the following results about the positivity of solutions.
  \begin{propo}
  For all $t \in [0, t_{\max}) $, the solution of system \eqref{Model diagram} remain nonnegative for nonnegative initial conditions.
  \end{propo}

\proof
	Let $n(a,t)= \mu(a)+\beta(a)\mathcal{W}(P(t),N_{F}(t))$.
	Straightforward computations (using the method of characteristics
	applied to the first equation of system \eqref{Model diagram} give the following formula for $S(a,t)$:
	$$
	S(a,t)   = \left\{
	\begin{array}{lll}
		S_{0}(a-t)\exp\left(-\int_0^tn(s+a-t,s)ds\right)   \quad &\mbox{if}& \quad   0\leq t\leq a \\
		m(t-a)\exp\left(-\int_0^an(s, s+t-a)ds \right)  \quad &\mbox{if}&
		\quad 0\leq a\leq t
	\end{array}
	\,.
	\right.
	$$
	From this formula, it's clear that $S(a,t)$ remains nonnegative for
	nonnegative data.
	In the same way, set  \[v(a)=\mu(a)+\dfrac{d(a) N_{I}}{K_d+\lVert I(\cdot\,,\,t)\lVert},\;
	\Pi(a)=\exp(-\int_{0}^{a}v(x)dx) \;\mbox{and}\;
 f(a,t) = \beta(a)\mathcal{W}(P(t),N_{F}(t))S(a,t)\,.\] 
 Using these new notations, the second equation of \eqref{Model diagram} reads
	\begin{equation}\label{eqI}
		\partial_tI(a,t)+\partial_aI(a,t)=-v(a)I(a,t)+ f(a,t)\,.
	\end{equation}
	Recall that this equation is supplemented with the initial and
	boundary conditions:
	\begin{equation}\label{ibc}
		I(0, t)= 0\qquad \mbox{and}\qquad I(a,0)=I_0(a) \;.
	\end{equation}
	Once more, using the method of characteristics one has
	\begin{align}\label{exp2}
		I(a,t)  = \left\{
		\begin{array}{lll}
			\frac{\Pi(a)}{\Pi(a-t)}I_0(a-t)+\displaystyle\int_0^t\frac{\Pi(a)}{\Pi(u+a-t)}f(u+a-t,u)du  & \mbox{if}&   0\leq t \leq a \\
			\displaystyle\int_{0}^{a}\frac{\Pi(a)}{\Pi(u)}f(u,u+t-a)du &\mbox{if}	
			& 0\leq a\leq t
		\end{array}
		\,.
		\right.
	\end{align}
	\noindent
	%
	Substituting these expressions in the two last equations of system \eqref{Model diagram}, we obtain 
	\begin{equation}\label{pox}
		\left\{
		\begin{array}{ll}
			\dot{N}_F(t) ~=~ g(N_F, N_I),\\
			\dot{N}_I(t) ~=~h(N_F, N_I).
		\end{array}
		\right.
	\end{equation}
Now we apply the \textit{barrier theorem for positivity of ODE solutions}.  
Let \( X = (N_F, N_I) \in \mathbb{R}^2 \) denote the state vector of system~\eqref{pox}, and suppose the initial condition \( X(0) \in \mathbb{R}_+^2 \).  
Assume that at some time \( t \geq 0 \), one component of \( X \), say \( X_i \), satisfies \( X_i(t) = 0 \), while the second component remain non-negative.  
If the time derivative \( \dot{X}_i(t) \geq 0 \), then \( X_i(t) \) cannot decrease below zero, and thus the solution remains in \( \mathbb{R}_+^2 \).  
This implies that the non-negative orthant is invariant under the flow of the system.
Consequently, we have three cases:
	\begin{itemize}
		\item If $N_F=0$ and  $N_I\geq 0$, one gets \[\dot{N_F}=g(0,N_I) =
		\left\{
		\begin{array}{l}
	  \gamma \int_{t}^{a_\dagger}\frac{\Pi(a)}{\Pi(a-t)}I_0(a-t)da   \quad \mbox{if}\quad 0\leq t \leq a_{\dagger}
	  \\
	  0\quad \mbox{if}\quad t> a_{\dagger}
		\end{array}
		\right.
		\geq 0 \,.
	\]	
		\item If $N_I=0$ and $N_F\geq 0$, we obtain
\(\dsp
\dot{N_I}   =    \alpha(\lVert\beta(\cdot) S(\cdot\,,\,t)\lVert + e \lVert\beta(\cdot)I(\cdot\,,\,t))\lVert) \mathcal{W}(P,N_{F})  
\).	From the expression of $I(a,t)$ given by (\ref{exp2}) and the hypothesis that $N_F\geq 0$ one easily obtains that $\dot{N_I}\geq 0.$ 
\item If $N_F=N_I=0$ at the same time, then 
		\begin{equation*}
			\left\{
			\begin{array}{ll}
				\dot{N}_F(t) ~=~ g(0, 0) ~=~\left\{
						\begin{array}{l}
					  \gamma \int_{t}^{a_\dagger}\frac{\Pi(a)}{\Pi(a-t)}I_0(a-t)da   \quad \mbox{if}\quad 0\leq t \leq a_{\dagger}
					  \\
					  0\quad \mbox{if}\quad t> a_{\dagger}
						\end{array}
						\right. \; \geq 0\,,\\
				\dot{N}_I(t) ~=~h(0, 0) ~=~0\,  \geq 0.
			\end{array}
			\right.
		\end{equation*}
	\end{itemize}
Hence, we conclude that both \( N_F(t) \) and \( N_I(t) \) remain non-negative for all \( t \in [0, t_{\max}) \), provided the initial conditions are non-negative.  
Furthermore, since \( S(a, t) \), \( N_F(t) \), and \( N_I(t) \) are non-negative, it follows that the expression for \( I(t, a) \) in equation~\eqref{exp2} is also non-negative.  
Therefore, for all non-negative initial conditions, the solutions of system~\eqref{Model diagram} remain non-negative on the existence interval \( [0, t_{\max}) \).
 \qed

We prove now that our solution is a global solution i.e
$t_{max}=+\infty$. According to the continuation principle, it will
be sufficient to prove that the solution is bounded (see the conclusion of Theorem~1.4,  p.185 of Pazy \cite{a4}). We have the
following result about the boundedness of trajectories.
   \begin{propo}
    Under the assumption of Section \ref{tata}, the solution of the model system (\ref{Model diagram}) is bounded on $[0, t_{max})$ and thus $t_{\max}=+\infty$. Moreover, the domain
   \begin{equation*}
    \varSigma := \{(S,I,N_{F},N_{I}) \in \mathbb{X}_{+} : \rVert S(\cdot\,,\,t)\rVert + \rVert I(\cdot\,,\,t)\rVert \leq \kappa_1, \,N_{F}(t)+ N_{I}(t) \leq  \kappa_2\}\,,
    \end{equation*}
    where
    \begin{align*}
  \kappa_1:=\max\left( \dfrac{\rVert m\rVert_{\infty}}{\tilde{\mu}}, \rVert S_{a0} \rVert + \rVert I_{a0} \rVert \right) \; \mbox{and} \; \kappa_2:=\max\left(N_{F0}+N_{I0}, \dfrac{4\gamma \kappa_1 + \rho K  \|d\|_{\infty}} {4\min\{\mu_{F}, \mu_{I}\}}\right)\,,
    \end{align*}
    is positively invariant under the flow of the system (\ref{Model diagram}).
     \end{propo}

   \proof
    Let $(S(\cdot\,,\,t),I(\cdot\,,\,t),N_{F}(t),N_{I}(t)) \in \mathbb{X}\,,$ be a solution of system \eqref{Model diagram}. Adding the integral (with respect to $a$) of
     the two first equations of system \eqref{Model diagram}, gives
    \begin{eqnarray*}
    \dfrac{d}{dt} (\rVert S(\cdot\,,\,t) \rVert_{L^{1}}+\rVert I(\cdot\,,\,t) \rVert_{L^{1}}) &=& -S(t,a_{\dagger})+S(0,t)-I(t,a_{\dagger})+I(0,t) \\
    &&- \int_{0}^{a_{\dagger}}
     \mu(a)(S(a,t)+I(a,t))da   -  \int_{0}^{a_{\dagger}}\dfrac{d(a)N_{I}I(a,t)}{K_d+\lVert I(\cdot\,,\,t)\lVert} da\,, \\
     & \leq & m(t)-\tilde{\mu}(\rVert S(\cdot\,,\,t)
      \rVert_{L^{1}}+\rVert I(\cdot\,,\,t) \rVert_{L^{1}})\,, \\
     & \leq & \rVert m\rVert_{\infty}-\tilde{\mu}(\rVert S(\cdot\,,\,t) \rVert_{L^{1}}+\rVert I(\cdot\,,\,t) \rVert_{L^{1}})\,.
    \end{eqnarray*}
    Hence using Gronwall's inequality, we have
   \begin{eqnarray*}
    \rVert S(\cdot\,,\,t) \rVert_{L^{1}}+\rVert I(\cdot\,,\,t) \rVert_{L^{1}} &\leq &\dfrac{\rVert m\rVert_{\infty}}{\tilde{\mu}}+( \rVert S_{a0}
    \rVert + \rVert I_{a0} \rVert - \dfrac{\rVert m\rVert_{\infty}}{\tilde{\mu}})\exp(-\tilde{\mu}t)\,,
     \\
      &\leq & \max\left( \dfrac{\rVert m\rVert_{\infty}}{\tilde{\mu}}, \rVert S_{a0} \rVert + \rVert I_{a0} \rVert \right) = \kappa_1 \,.
   \end{eqnarray*}
 Adding the third and the fourth equations of system \eqref{Model diagram} gives
$$\dot{N_{F}} +  \dot{N_{I}} =  \gamma \lVert I(\cdot\,,\,t)\lVert-\mu_{F}N_{F} -\mu_{I}N_{I} + \dfrac{\rho N_{I}\lVert d(\cdot)I(\cdot\,,\,t)\lVert}{K_d+\lVert I(\cdot\,,\,t)\lVert}\left(1-\dfrac{N_I}{K}\right).$$
Using the fact that \[\dfrac{\rho N_{I}\lVert d(\cdot)I(\cdot\,,\,t)\lVert}{K_d+\lVert I(\cdot\,,\,t)\lVert}\left(1-\dfrac{N_I}{K}\right) \leq \dfrac{\rho K}{4}\dfrac{\lVert d(\cdot)I(\cdot\,,\,t)\lVert}{K_d+\lVert I(\cdot\,,\,t)\lVert} \leq \dfrac{\rho K  \|d\|_{\infty} }{4}\, ,\]
 one has 
$$\dot{N_{F}} +  \dot{N_{I}} \leq \gamma D + \dfrac{\rho K  \|d\|_{\infty} }{4}-\min\{\mu_{F}, \mu_{I}\}\left(N_{F} +  N_{I}\right)$$
Thus for all $t \geq 0$, using again the Gronwall inequality, we have
\begin{align*}
N_{F} +  N_{I} & \leq \max\left(N_{F0}+N_{I0}, \dfrac{4\gamma \kappa_1 + \rho K  \|d\|_{\infty}} {4\min\{\mu_{F}, \mu_{I}\}}\right) = \kappa_2\,.
\end{align*}
Thus, the solution of the system \eqref{Model diagram} is bounded on $[0, t_{max})$ and consequently $t_{\max}=+\infty$. \qed
   
\subsection{Regularity of solutions}
Now that the existence and uniqueness of a global (in time), positive, and bounded mild solution have been established, we investigate under what additional conditions on the data and parameters System~\ref{Model diagram} admits a unique global, positive, and bounded classical (regular) solution $-$ that is, a unique solution defined on $[0, \infty)$, continuous on $[0, \infty)$, and continuously differentiable on $(0, \infty)$.
  The result here will be an application of Theorem~1.7, page~190 of \cite{a4}. To this end, we introduce the Banach space \( \mathbb{Y} = D(\mathcal{A}) \), endowed with the graph norm \( |\cdot|_{\mathbb{Y}} \); that is, for \( u = (S, I, N_F, N_I) \in D(\mathcal{A}) \), we define
\(
|u|_{\mathbb{Y}} = \|u\|_{\mathbb{X}} + \|\mathcal{A}u\|_{\mathbb{X}},
\)
(recall that,  \( {\mathbb{X}} \) is the ambient Banach space in which the operator \( \mathcal{A} \) is densely defined) i.e.
\begin{eqnarray} | u |_\mathbb{Y} &=&  \|u\|_\mathbb{X}+ \|\mathcal{A}(u)\|_\mathbb{X}
\,, \no\\
&
=
&
\|S\|_{W^{1,1}(0,a_{\dagger})}+\|I\|_{W^{1,1}(0,a_{\dagger})}+(1+\mu_F)|N_F|+(1+\mu_I)|N_I\,.
\no
 \end{eqnarray}
We have the following (from now on, we adopt the convention that the overdot $(\dot{\ })$ denotes differentiation with respect to $t$, while the prime $( ')$ denotes differentiation with respect to $a$). 
\begin{propo}\label{Reg}
In addition to the assumptions of Section~\ref{tata}, suppose that		
$ m \in C^{1}_b([0,\infty))$, $\beta, d \in C^{1}_b([0,a_{\dagger}))$, and $\pi \mu,\, \pi \mu',\, \pi\mu^2 \in L^1(0, a_{\dagger})$. Then the mapping $\mathcal{H}$ is defined on $\mathbb{Y}$ with values in $  \widetilde{\mathbb{Y}} := W^{1,1}(0, a_{\dagger}) \times W^{1,1}(0, a_{\dagger}) \times \mathbb{R} \times \mathbb{R}$, and is locally Lipschitz continuous in $u$, uniformly in \( t \geq 0 \); that is, for all \( \delta > 0 \), there exists a constant \( \kappa_\delta > 0 \) such that
\[
\|\mathcal{H}(u_{2}) - \mathcal{H}(u_{1})\|_{ \widetilde{\mathbb{Y}}} \leq \kappa_{\delta} \|u_{2} - u_{1}\|_{\mathbb{Y}}, \, \forall (u_{1}, u_{2}) \in \tilde{B}_{\delta} \times \tilde{B}_{\delta}\,,
\]
where \( \tilde{B}_{\delta} \) denotes the ball in \( \mathbb{Y} \) centered at the origin with radius \( \delta \).
\end{propo}
\proof 
First, we prove that for all
\(\displaystyle  
u  \in  \mathbb{Y} ,
\)
 \(\displaystyle \mathcal{H}(u) \in \widetilde{\mathbb{Y}}\).  
To establish this, it suffices to show that
\( \displaystyle 
  (\mu \phi)\,' \in L^1(0, a_{\dagger})\)
\,
for any fixed \(t > 0\) (see Proposition~\ref{prop_pi}).
We have 
\[
  (\mu \phi)\,'= \left\{ \begin{array}{ll}
    m(t-a)\big(\mu'(a)\pi(a) - \pi(a)\mu^2(a)\big) - m'(t-a)\pi(a)\mu(a) \hspace{1em} \mbox{if} \hspace{1em}  0\leq a \leq t    \\
     0 \hspace{1em} \mbox{otherwise} 
    \end{array}
     \right.
     \,,
\]
and from hypotheses, it easily follows that \( \displaystyle 
  (\mu \phi)\,' \in L^1(0, a_{\dagger})\). 
Secondly, we show that the perturbation mapping \(\mathcal{H}\) is locally Lipschitz continuous with respect to \(u\), uniformly for all \(t \geq 0\). 
The notation here are those of the proof of Proposition~\ref{tup}, see Equation~\ref{gh}. We have
\begin{align*}
  A\,'(t,u) & = \beta\,'  \left(-(\mathcal{W}^{1}-\mathcal{W}^{2})S^{1}- \mathcal{W}^{2}(S^{1}-S^{2})- \phi(\mathcal{W}^{1}-\mathcal{W}^{2})\right)\,,\\
	& + \beta  \left(-(\mathcal{W}^{1}-\mathcal{W}^{2}) S^{1}\,'- \mathcal{W}^{2}( S^{1}\,' S^{2}\,')-  \phi\,'(\mathcal{W}^{1}-\mathcal{W}^{2})\right)\,.
\end{align*}
Proceeding exactly as in the proof of Proposition \ref{tup}, we obtain that  
\begin{equation}\label{er11}
\|A'\|_{L^1(0, a_{\dagger})} \leq \kappa_1(\delta) \|u_2 - u_1\|_{\mathbb{Y}}\,,
\end{equation}
where \(\kappa_1(\delta)\) is a constant depending on \(\theta_1, \theta_2, \theta_3, M, \delta\), and on the \(L^\infty\)-norms of \(\beta\) and \(\beta'\); its exact value is irrelevant.
In the same way one has the  following straightforward computations 
\begin{small}
	\begin{align*}
	 B\,'(t,u) = & -\frac{d\,'(a)}{\left(K_d+\left\lVert \pi I^2\right\lVert\right)\left(K_d+\left\lVert \pi I^1\right\lVert\right)}\left[\left(K_d + \lVert \pi I^1 \lVert \right)\left(N_I^1(I^1-I^2)+ I^2(N_I^1-N_I^2) \right) \right]\\
		& -\frac{d'(a)}{\left(K_d+\left\lVert \pi I^2\right\lVert\right)\left(K_d+\left\lVert \pi I^1\right\lVert\right)}\left[N_I^1I^1\left(\left\lVert \pi I^2 \right\lVert - \left\lVert \pi I^2 \right\lVert\right)\right] \\
		& -\frac{d(a)}{\left(K_d+\left\lVert \pi I^2\right\lVert\right)\left(K_d+\left\lVert \pi I^1\right\lVert\right)}\left[\left(K_d + \lVert \pi I^1 \lVert \right)\left(N_I^1( I^1\,'- I^2\,')+  I^2\,'(N_I^1-N_I^2) \right) \right]\\
		& -\frac{d(a)}{\left(K_d+\left\lVert \pi I^2\right\lVert\right)\left(K_d+\left\lVert \pi I^1\right\lVert\right)}\left[N_I^1 I^1\,'\left(\left\lVert \pi I^2 \right\lVert - \left\lVert \pi I^2 \right\lVert\right)\right]- A\,'(t,u)\,,		
	\end{align*}
\end{small}
which leads, as in previous case to the estimate:
\begin{equation}\label{er12}
 \|B'\|_{L^1(0, a_{\dagger})} \leq \kappa_2(\delta) \|u_2 - u_1\|_{\mathbb{Y}}\,.
\end{equation}
Here again, $\kappa_2(\delta)$ is a generic constant as in (\ref{er11}) depending on 
 \(K_d, \theta_1, \theta_2, \theta_3, M, \delta\), and on the \(L^\infty\)-norms of \(\beta\,,\,\,\beta'\) and $d\,'$.
Combining Inequalities (\ref{er11}) and (\ref{er12}), we obtain
\begin{equation*}
\|\mathcal{H}(t,u_{2}) - \mathcal{H}(t,u_{1})\|_{ \widetilde{\mathbb{Y}}} \leq \kappa_{\delta} \|u_{2} - u_{1}\|_{\mathbb{Y}}
	 \,,
\end{equation*}
where $\displaystyle \kappa_{\delta}= \max(\kappa_1(\delta), \kappa_2(\delta))$.
The proof is then complete.
\qed
\begin{example}
The hypotheses of the previous Proposition are satisfied for the following choices of the rates $\mu$, $\beta$, and $d$:
\[
\mu(a)= \frac{\alpha_0}{(a_{\dagger}-a)^m}, \,  \;\; \beta(a)=\beta_{max}e^{-\frac{(a-a_{opt})^2}{2\sigma_p^2}}\;\;\mbox{and}\;\; d(a)= d_{0}e^{-\nu \,a}
\,.
\]
These choices are biologically motivated. See for example~\cite{a2} for different models of mortality rates.
\end{example}

\begin{theo}[Regularity] Under the hypotheses of Proposition~\ref{Reg}, the model system~\ref{Model diagram}
admits a unique, global, positive and bounded  \textbf{strong solution} which satisfies
$$
 u = (S, I, N_F, N_I) \in C([0, \infty); D(\mathcal{A})) \cap C^1([0, \infty); \mathbb{X}).
$$
\end{theo}
\proof
The theorem follows from the previous Proposition and Theorem~1.7 page 190 of Pazy,~\cite{a4}.
\qed 

\subsection{Asymptotic analysis for constant recruitment}\label{se3}
 \noindent

 In this section, we discuss the stability of the pest-free steady states or time independent solutions of the system \eqref{Model diagram}. We assume that the recruitment of new healthy banana plantain is constant, $\, m(t)=m \in \mathbb{R}_{+}$. We compute the disease-free steady state of system \eqref{Model diagram}. At this state, there is no infection ($I=0, \,  N_I= 0$) and
  the derivatives with respect to time are zero, we are thus led to a system with
 only one derivative with respect to age which has the unique solution (pest-free steady state) $\xi(a)=(S^{0}(a),0,0,0)$  where $S^{0}(a)=m \pi(a)$.
  Now, to compute the basic reproduction number,  $\mathcal{N}$, one can use the next generation operator approach (Diekman et al. (\cite{a11}, 1990), Inaba (\cite{a5}, 2017)) but instead, since the eigenvalue problem which arises from the linearized system is rather simple we proceed in a straightforward manner. We linearize the system \eqref{Model diagram} around the pest-free steady state $\xi(a)=(S^{0}(a),0,0,0)$.
   Set $$x(a,t)=S(a,t)-S^{0}(a), \, y(a,t)=I(a,t), \, v(t)=N_{F}(t) , \, w_{i}(t)=N_{I}(t)
  \,. 
    $$
%
 \noindent  
The linearized system  thus reads:
  \begin{equation}\label{c}
  \left\{\begin{array}
  {llll}
   x_t + x_a & = &-\ell(a)v-\mu x \\
   y_t + y_a & = &\ell(a)v-\mu y \\
   \dot{v} & = & \gamma \|y(a)\| - (B+\mu_F)v_0\\
  \dot{w} & = & \sigma v - \mu_I w
   \end{array}
   \,,
   \right.
   \end{equation}
  where $\dsp \ell(a) = \beta(a)\dfrac{S^0}{\|S^0\|} $ and $\dsp  \sigma=\alpha \dfrac{\|\beta S^0\|}{\|S^0\|}$. Next, we introduce the eigenvalues of system (\ref{c}) by looking for their solutions of the form 
  $$x(a,t)=e^{\lambda t} x(a), \, y(a,t)=e^{\lambda t} y(a), \, v(t)=e^{\lambda t} v_0 , \, w(t)=e^{\lambda t} w_0\,, $$
 which leads to the following eigenvalue problem:
 \begin{equation}\label{d}
  \left\{\begin{array}
  {llll}
\lambda x(a) + \dot{x}(a) & = &-\ell(a)v_0-\mu(a)x(a) \\
\lambda y(a) + \dot{y}(a) & = &\ell(a)v_0-\mu(a)y(a) \\
\lambda v_0 & = & \gamma \|y(a)\| - (B+\mu_F)v_0\\
\lambda w_0 & = & \sigma v_0 - \mu_I w_0
   \end{array}
     \right.
     \,,
     \end{equation}
with the initial conditions $x(0)=y(0)=0$. We need to determine the $\lambda$ that allow the system above to have a nonzero solution. It is easy to see that the unknowns $x(t)$ and $y(t)$ read:
\begin{equation}\label{e}
x(a) = - y(a) = -v_0\int_0^a \ell(\zeta)e^{-\int_{\zeta}^{a}(\lambda+\mu(r))dr}d\zeta.
\end{equation}
Substituting expressions \eqref{e} into the third equation of system \eqref{d} gives
$$v_0 = \dfrac{\gamma}{\sigma+\mu_F}v_0 \int_0^{a_{\dagger}} \int_0^a \ell(\zeta)e^{-\int_{\zeta}^{a}(\ell+\mu(r))dr}d\zeta da -\dfrac{\lambda}{\sigma+\mu_F}v_0
\,.
$$
Since we are looking for nonzero solutions ($v_0\neq 0 $) one obtains the following characteristic equation 
\begin{equation}\label{f}
1 = \dfrac{\gamma}{\sigma+\mu_F} \int_0^{a_{\dagger}} \int_0^a \ell(\zeta)e^{-\int_{\zeta}^{a}(\lambda+\mu(r))dr}d\zeta da -\dfrac{\lambda}{\sigma+\mu_F} =: \mathcal{K}(\lambda)\,.
\end{equation}

\noindent
 Set,
$$
\mathcal{N} = \mathcal{K}(0) = \dfrac{\gamma}{\sigma+\mu_F} \int_0^{a_{\dagger}} \int_0^a \ell(\zeta)e^{-\int_{\zeta}^{a}\mu(r)dr}d\zeta da
\,.
$$
We have the following 
\begin{theo}\label{thm1}
The hypotheses are those of Section~\ref{tata}. In addition, if we assume that the recruitment of new healthy banana plantain plants is constant ($m(t) = m > 0$), then the pest-free steady state $\xi(a)$ is locally asymptotically stable (l.a.s.) when  $\mathcal{N} < 1$
and unstable in the case where  $\mathcal{N} > 1$. 
\end{theo}
\proof
We show that Equation~\ref{f} admits a unique real solution $\lambda^{\star}$.
Viewed as a real-valued function of the real variable $\lambda$, 
the function $\mathcal{K}$ is decreasing and satisfies
\[
\lim_{\lambda \to +\infty} \mathcal{K}(\lambda) = -\infty
\quad\text{and}\quad
\lim_{\lambda \to -\infty} \mathcal{K}(\lambda) = +\infty .
\]
Therefore, Equation~\ref{f} has a unique real solution $\lambda^{\star}$ such that
\[
\lambda^{\star} < 0 \;\;\text{iff}\;\; \mathcal{N} < 1, \quad
\lambda^{\star} > 0 \;\;\text{iff}\;\; \mathcal{N} > 1  \quad\mbox{and}\quad
\lambda^{\star} = 0 \;\;\text{iff}\;\; \mathcal{N} = 1 .
\]
To complete the proof, we show that the real part of any complex solution of
Equation~\ref{f} is strictly less than $\lambda^{\star}$.
Let $\lambda =\lambda_1 +i\lambda_2$ be an arbitrary complex solution of Equation~\ref{f}, then we have (recall $\mathcal{R}e(z)$ denote the real part of $z$)
\begin{eqnarray}
\mathcal{K}(\lambda^{\star}) &=& \mathcal{K}(\lambda) \,,
\no\\
&=&\mathcal{R}e\big(\mathcal{K}(\lambda)\big)\,,
\no
\\
&=&
 \dfrac{\gamma}{\sigma+\mu_F} \int_0^{a_{\dagger}} \int_0^a \ell(\zeta)\cos(\lambda_2(a-\zeta))e^{-\int_{\zeta}^{a}(\lambda_1+\mu(r))dr}d\zeta da -\dfrac{\lambda_1}{\sigma+\mu_F} 
 \,,
\no
 \\
 &\leq&
 \mathcal{K}(\lambda_1)
 \,.
 \no      
\end{eqnarray}
This inequality and the properties of \(\mathcal{K}\) prove that $\dsp \lambda_{1}\leq \lambda^{\star}.$
  Hence $\xi(a)$ is l.a.s. if $\mathcal{N} < 1$ and unstable if $\mathcal{N} > 1$.
\qed
\begin{remark} $\mathcal{N}$, is the basic reproduction number and 
biologically, it is the average number of new infesting nematodes originated from a single free nematode during its lifespan (see \cite{a5,a11} for more details).
Mathematically $\mathcal{N}$ serves as threshold for the stability of the disease-free equilibrium. %
\end{remark}
\begin{remark}
Theorem \ref{thm1} means biologically that,  in the case where  $\mathcal{N}< 1$, the nematode population will disappear 
from the banana-plantain plantation if the initial sizes of
sub-population of susceptible are in the basin of attraction
 of the pest-free steady state. So, in this case it is not necessary to control the \textit{Radopholus similis}.
\end{remark}
We turn now to the proof of the global stability of the pest free steady state. We first introduce the following Lemma.
\begin{lem}\label{ll}
For any $ a \, \in  [0, a_{\dagger})$, if $\;S_0(a)+I_{0}(a) \leq
S^{0}(a)$ then  for all $t\geq 0$ ,  $ S (a,t) + I(a,t) \leq
S^{0}(a)$.
\end{lem}

\begin{proof}
 Set $V(a,t)=S(a,t)+I(a,t)-S^{0}(a)$. From the first and second equations of system 
\ref{Model diagram}, we have
\begin{equation}\label{ba1}
 \partial_{t}V(a,t)+\partial_{a}V(a,t)  =  -\mu(a)V(a,t)-z(a,t)
 \,,
\end{equation}
where   $\dsp z(a,t)=\dfrac{d(a)N_{I}I(a,t)}{K_d+\lVert I(\cdot\,,\,t)\lVert}\,.$ 
Note that equation \eqref{ba1} is supplemented with the following
boundary and initial conditions:
\begin{equation}\label{InCond2}
V(0,t)=0,\,\,\, V(a,0)=S_{0}(a)+I_{0}(a)-S^{0}(a)=V_0(a)
\,.
\end{equation}
Solving   \eqref{ba1}-\eqref{InCond2} along the characteristic line $t-a=c$,
gives the following representation formula for $V$:
$$ V(a,t)  = \left\{
         \begin{array}{l}
      V_0(a-t)e^{ -\int_{0}^{t}\mu(u+a-t)du  }  
       -\displaystyle\int_0^t z(r+a-t,r)e^{ -\int_{r}^{t}\mu(u+a-t)du  }dr \quad   \mbox{if}  \quad   0\leq t \leq a \\
       -\displaystyle\int_{0}^{a}z(r,r+t-a)e^{ -\int_{r}^{a}\mu(u)du  }dr \quad 
       \mbox{if} \quad 
         0\leq a\leq t
       \end{array}
        \right.
        \,.
$$
From the  hypothesis we readily see that $V(a,t) \leq
0$ . This completes the proof of the Lemma.
\end{proof}
\noindent 
Proceeding as in the proof of  Theorem~4.2 of \cite{castillo1998}, we can prove the following theorem
\begin{theo}
The hypotheses are those of Theorem~\ref{thm1}, with the additional condition that 
\begin{equation} \label{mbmu}
   \dfrac{m}{b}\int_{0}^{a_{\dagger}}\pi(a)da \geq 1
  \,.
\end{equation} 
Then , the  pest-free equilibrium $\xi(a)$ is globally asymptotically stable if 
$$\mathcal{N} <  \dfrac{(e\alpha\underline{\beta}+\mu_F)b}{(\sigma+\mu_F)\|S^0\|} =:\mathcal{N}_0\,,$$
where \( \dsp \underline{\beta} = \inf_{0 \leq a \leq a_{\dagger}} \beta(a)\,.\)
\end{theo}
\begin{remark}
The term \(~\dsp  \dfrac{m}{b}\int_{0}^{a_{\dagger}}\pi(a)da \) 
 represents the ratio between the equilibrium root biomass, maintained by a constant recruitment rate 
\(m\), and the minimal root banana plantain \(b\). It serves as an indicator of the banana plantain  system's renewal capacity: values above one indicate sufficient recruitment to sustain or increase banana plantain biomass, whereas values below one reflect insufficient renewal due to high mortality, such as that caused by Radopholus similis infestation.
\end{remark}

\proof
First, we observe that, from the additional hypothesis~\ref{mbmu}, it follows that $\mathcal{N}_0 \leq 1$. In fact, by the definition of  \( \dsp \underline{\beta}\), it follows that
 \begin{equation}\label{z1} 
\frac{\alpha \left( \lVert \beta S \rVert + e \lVert \beta I \rVert \right)}{P} \geq e \alpha \, \underline{\beta}
\quad \mbox{and} \quad  
\frac{e \alpha \, \underline{\beta} + \mu_F}{\sigma + \mu_F} \leq 1\,.
\end{equation}
We consider the long-term behavior of system~\ref{Model diagram}. Consequently, without loss of generality, we assume $t > a_{\dagger}$, the maximum age of a banana plantain plant. Solving the second equation in system \ref{Model diagram}, we obtain 
\begin{equation}\label{s4}
I(a,t)   = \int_{0}^{a} \omega(u,u+t-a) e^{-\int_{u}^{a}\nu(r)dr}du, 
\end{equation}
where we have set $\omega(a,t) = \beta(a)\mathcal{W}(P(t),N_{F}(t))S(a,t)$  and  $\nu(a)=\mu(a)+\dfrac{d(a)N_{I}}{K_d+\lVert I(\cdot\,,\,t)\lVert}\cdot$
Using assumption (A4) and  Lemma \ref{ll}, we obtain   
\begin{align}\label{s5}
I(a,t) & \leq \dfrac{1}{b}\displaystyle\int_{0}^{a} e^{-\int_{u}^{a}\mu(r)dr}\beta(u)N_F(u+t-a) S^0(u) du
\,.
\end{align}
On the other hand, the first inequality of \eqref{z1} and the   third equation of system~\eqref{Model diagram} yield
\begin{align}\label{s3}
\dot{N}_F  & \leq  -(e \alpha \, \underline{\beta}+\mu_F) N_F+\gamma\lVert I(\cdot\,,\,t) \lVert\,.
\end{align}
By Gronwall's Lemma one has  
\begin{align*}
N_F & \leq N_{F0}e^{-(e \alpha \, \underline{\beta}+\mu_F)t}+\displaystyle\int_{0}^{t}\gamma \lVert I(\cdot\,,\,s) \lVert e^{-(e \alpha \,  \underline{\beta}+\mu_F)(t-s)}ds\,,\\
& \leq N_{F0}e^{-(e \alpha \, \underline{\beta}+\mu_F)t}+\displaystyle\int_{0}^{t}\gamma \lVert I(\cdot\,,\,t-x) \lVert e^{-(e \alpha \,  \underline{\beta}+\mu_F)x}dx\,,\quad  (x = t-s)\,,\\
& \leq N_{F0}e^{-(e \alpha \, \underline{\beta}+\mu_F)t}+\displaystyle\int_{0}^{\infty}\gamma \lVert I(\cdot\,,\,t-x) \lVert e^{-(e \alpha \,  \underline{\beta}+\mu_F)x}dx \,.
\end{align*}
The estimate of $I(a,t)$ given in Inequality~\ref{s5} leads to the following
\begin{eqnarray}\label{s1}
N_F   &\leq & N_{F0}e^{-(e \alpha \, \underline{\beta}+\mu_F)t} \nonumber 
\\
&& \; + \dfrac{\gamma}{b}\displaystyle\int_{0}^{\infty}\displaystyle\int_{0}^{a_{\dagger}}\displaystyle\int_{0}^{a} e^{-\int_{u}^{a}\mu(r)dr}\beta(u) S^0(u)N_F(u+t-x-a)e^{-(e \alpha \, \underline{\beta}+\mu_F)x}\, du\, da\, dx \nonumber\,.
\end{eqnarray}
Let 
$$C = \limsup_{t \rightarrow \infty}N_F(t)
\,.
$$
Taking the limit superior as 
$t \rightarrow \infty$ on both sides of Inequality~\eqref{s1}, and applying Fatou's Lemma, we obtain
\begin{align*}
C  & \leq \dfrac{C\gamma}{b}\displaystyle\int_{0}^{\infty}\displaystyle\int_{0}^{a_{\dagger}}\displaystyle\int_{0}^{a} e^{-\int_{u}^{a}\mu(r)dr}\beta(u) S^0(u)   e^{-(e \alpha \, \underline{\beta}+\mu_F)x} \,du\, da\, dx \,,\\
 & \leq \dfrac{C\gamma}{b}\displaystyle\int_{0}^{a_{\dagger}}\displaystyle\int_{0}^{a} e^{-\int_{u}^{a}\mu(r)dr}\beta(u) S^0(u)   \displaystyle\int_{0}^{\infty} e^{-(e \alpha \, \underline{\beta}+\mu_F)x}\,dx\, du\, da\,.
 \end{align*}
 After integration, one gets
 \begin{align}\label{s2}
 C & \leq C \dfrac{\gamma\|S^0\|}{b(e \alpha \, \underline{\beta}+\mu_F)}\displaystyle\int_{0}^{a_{\dagger}}\displaystyle\int_{0}^{a} e^{-\int_{u}^{a}\mu(r)dr}\beta(u) \ell(u) \,du\, da\,.
\end{align}
Hence, Inequality~\ref{s2} becomes
 \[  C \leq C\dfrac{\|S^0\|(\sigma+\mu_F)}{b(e \alpha \, \underline{\beta}+\mu_F)}\mathcal{N}
 \quad \text{i.e.}\quad  C \leq C\dfrac{\mathcal{N}}{\mathcal{N}_0}
\,.
\] 
The last inequality shows that, in the case where \( \mathcal{N} < \mathcal{N}_0 \), one has  
\(\dsp 
C = \limsup_{t \rightarrow \infty} N_F(t) = 0,
\)
and and since $\dsp N_F(t)\geq 0$ it follows that
\(\dsp 
\lim_{t \rightarrow \infty} N_F(t) = 0.
\)
Now, 
from Inequality~\ref{s4}, we see that $\lim\limits_{t \rightarrow \infty}I(a,t)=0$, and consequently,  
$$\lim\limits_{t \rightarrow \infty} N_I(t) = 0 \quad \mbox{
and} \quad \lim\limits_{t \rightarrow \infty} S(a,t) = S^0(a).$$ 
\qed 
\section{Numerical analysis and simulations}
\noindent
In this section, we present the semi-implicit Euler scheme employed for the numerical study of the model introduced at the beginning of this paper. 
We show that the scheme is consistent, in the sense that it converges to the continuous solution when the age-dependent parameters and the data are sufficiently smooth. 
Finally, we provide numerical simulations to illustrate and validate the theoretical results.

\subsection{Discretization   and basic properties of the discrete solutions}

We now present the numerical approximations of the solutions to our model system~(\ref{Model diagram}). To compute these approximate solutions, we introduce a discretization of the time and age intervals. Specifically, we consider a set of points \( (t^n)_{n=0, ~  \ldots\, ,~ N} \) on the time interval \([0, T]\), and a set of points \( (a_j)_{j=0, \ldots, M} \) on the age interval \([0, a_{\dagger})\).
For simplicity, we use the same constant discretization parameter for both time and age, denoted by \( h := \Delta t = \Delta a \).

Let \( S_j^n, I_j^n, N_F^n, N_I^n \), and \( P^n \) be numerical approximations of 
\( S(a_j, t^n) \), \( I(a_j, t^n) \), \( N_F(t^n) \), \( N_I(t^n) \), and \( P(t^n) \), respectively, 
at age \( a_j = jh \) and time \( t^n = nh \). 
Since the banana-plantain plant requires six (06) months to mature after flowering, we simulate system~(\ref{Model diagram})
over the time interval \( [0, T] \), with \( T = 360 \) days and \( a_{\dagger} = 300 \) days.
We set 
\( \dsp
T = hN \quad \text{and} \quad a_{\dagger} = h M.
\)
The numerical approximation of \( P(t^n) = P(nh) \) is given by:
\begin{equation}
	\label{pn}
	P(nh ) \approx P^n := h \sum_{j=1}^{M}p^n_j\quad \mbox{where}\quad p_j^n :=   S_j^n + I_j^n 
	\,.
\end{equation}
Using the semi-implicit Euler method, the discrete unknowns satisfy the following system of equations, which corresponds to a discretization of the model system~(\ref{Model diagram}) for \(1 \leq j \leq M\) and \(1 \leq n \leq N\):
\begin{eqnarray}
 \label{eqq1aa}
	\dfrac{S_j^{n} - S_j^{n-1}}{h} + \dfrac{S_j^{n-1} - S_{j-1}^{n-1}}{h} &=& -\left( \Lambda_{j}^{n-1} + \mu_j \right) S_j^n\,,
\\
\label{eqq1ab}
	\dfrac{I{}_j^{n} - I{}_j^{n-1} }{h } +  \dfrac{I{}_j^{n-1}  - I{}_{j-1}^{n-1} }{h } &=& \Lambda^{n-1}_{j} S{}^{n}_j - \left(\mu_{j}+\frac{d_{j}N_{I}^{n-1}}{K_d + B^{n-1}}\right) I{}^{n}_j\,,
\label{eqq1ac}
\\
	\dfrac{N_F^{n}-N_F^{n-1}}{h } &=&  -\alpha \frac{C^{n-1}}{P^{n-1}}N_F^{n} +\gamma B^{n-1} - \mu_F N_F^{n} 
	\,,
\label{eqq1ad}
\\
	\dfrac{N_I^{n}-N_I^{n-1}}{h } &=&  \alpha \frac{C^{n-1}}{P^{n-1}}N_F^{n} -\mu_{I}N_I^{n} + \rho \dfrac{N_{I}^{n-1}D^{n-1}}{K_d+B^{n-1}}\left(1-\dfrac{N_I^{n}}{K}\right)
	\,,
\end{eqnarray}
where 
$$
\mu_j= \mu(jh ), \;  \beta_j= \beta(jh ), \;    d_j= d(jh ), \;\Lambda(a,t)=   
\beta(a)\mathcal{W}( P(t),N_{F}(t))\,,  $$    
 $$
  B(nh)=\int_0^{a_{\dagger}}I(a, nh)da\,,\; C(nh)=\int_0^{a_{\dagger}}\beta(a)(S(a,nh)+eI(a,nh))da \,,\; D(nh)= \int_0^{a_{\dagger}}d(a)I(a,nh)da\,, $$ 
and 
$$\Lambda^n_j = \beta(a_j)\mathcal{W}( P^n, N_{F}^n), \; 
B^{n} =h \sum_{j=1}^{M}I{}^n_j  , \; 
C^{n} =h \sum_{j=1}^{M}\beta_j(S{}^n_j+eI{}^n_j) \;  \mbox{and}
\;
D^{n} =h \sum_{j=1}^{M}d_jI{}^n_j  \,.
$$    
Reorganizing the terms of the previous equations yields the following induction scheme for the discrete unknowns:
\begin{equation}
	S{}_j^{n} = \dfrac{S{}_{j-1}^{n-1}}{1 + h  \left(\mu_{j} + \Lambda^{n-1}_{j}\right)} 
	\,,\quad 
	I{}_j^{n}  = \dfrac{I{}_{j-1}^{n-1} +  h ~ \Lambda^{n-1}_{j}S{}_j^{n}}{1 + h  \left(\mu_{j}+\frac{d_{j}N_{I}^{n-1}}{K_d+B^{n-1}}\right)} 
	\,,
\end{equation}
\begin{equation}
	N_F^{n} = \dfrac{N_F^{n-1}+ h ~ \gamma B^{n-1}}{1 + h  \left(\mu_F+ \dfrac{\alpha C^{n-1}}{P^{n-1}}\right)}   
	\,, \quad 
	N_I^{n} = \dfrac{N_I^{n-1}\left(1+ \dfrac{ h ~ \rho D^{n-1}}{K_d+B^{n-1}}\right) + h ~\dfrac{  \alpha C^{n-1}}{P^{n-1}}N_F^n}{1 + h  ~\left(\mu_{I}+ \dfrac{\rho N_I^{n-1} D^{n-1}}{K(K_d+B^{n-1})}\right)}\,,
\end{equation}
with the following initialization:   
$$
S{}^n_{0} = m , \;  S{}^0_{j} = S_0(a_{j}),\;  I{}^{n}_0 = 0, \; I{}^0_{j} = I_0 (a_{j}),\; N_F^0=N_{F0},\; N_I^0=N_{I0} \,.      
$$
This proves the existence of the numerical scheme and the non-negativity of the approximation values at each step, assuming non-negative initial and boundary conditions. For later use, we note that adding equations (\ref{eqq1aa}) and (\ref{eqq1ab}) gives:
\begin{equation}\label{eq01}
p_j^{n}  = \dfrac{  p_{j-1}^{n-1} +  h~\dfrac{d_j N_{I}^{n-1}}{K_d+B^{n-1}}S{}_j^n}{1 + h  \left(\mu_{j} + \frac{d_{j} N_{I}^{n-1}}{K_d+B^{n-1}} \right)} 
\,.
\end{equation} 
\begin{lem} Under the hypotheses of Theorem~\ref{thm1},  for all $\, 0\leq n\leq N \,$ and $\, 0\leq j\leq M\,,$
	\begin{equation}
		 0\leq S{}^n_j\,,\;I{}^n_j\,,\;   p_j^n \,  \leq \widetilde{K}_0:=\max(m, \|S_{0}\|_{\infty}+\|I_{0}\|_{\infty})
		\,.
	\end{equation}
\end{lem} 
\proof
Recall that by assumption $(A3)$, $\frac{S_{0}}{\pi}\,,\; \frac{I_{0}}{\pi} \in W^{1,1}(0,a_{\dagger}) \hookrightarrow L^{\infty}(0,a_{\dagger})$. Therefore $S_{0}, \, I_{0}$ are essentially bounded. Next, observe that for $\, 0\leq n\leq N,\,$ and $\, 0\leq j\leq M$, one has $\displaystyle  0\leq S{}^n_{0}\leq \max(m, \|S_{0}\|_{\infty})$ and $\displaystyle \, 0\leq S{}^0_{j}\leq \max(m, \|S_{0}\|_{\infty}).$ Now, by induction, assume that $\displaystyle 0\leq S{}^{n-1}_{j-1}\leq \max(m, \|S_{0}\|_{\infty})$ then we have,
$$
0\leq  S{}_j^{n} = \dfrac{S{}_{j-1}^{n-1}}{1 + h  \left(\mu_{j} + \Lambda^{n-1}_{j}\right)}\leq S{}_{j-1}^{n-1} \leq  \max(m, \|S_{0}\|_{\infty})\,. 
$$
Finally, assuming that $ \displaystyle 0\leq p^{n-1}_{j-1}\leq \widetilde{K}_0$, leads to the following inequality (see \eqref{eq01}) 
$$ 0\leq  p_j^{n}  = \dfrac{p_{j-1}^{n-1} +  h ~\frac{d_jN_{I}^{n-1}}{K_d+B^{n-1}}S{}_j^n}{1 + h  \left(\mu_{j} + \frac{d_jN_{I}^{n-1}}{K_d+B^{n-1}} \right)} \leq   \dfrac{  \widetilde{K}_0\left(1 +  h \frac{ d_jN_{I}^{n-1}}{K_d+B^{n-1}}\right)}{1 + h  \left(\mu_{j} +\frac{ d_{j}N_{I}^{n-1}}{K_d+B^{n-1}} \right)}
\leq  \widetilde{K}_0\,,
$$  
and the proof is complete.
\qed
\begin{remark}
	Note that for $\,0\leq n \leq N,$ one has
	\[
0\leq B^n\leq a_{\dagger}\widetilde{K}_0\,,\quad  0\leq C^n\leq a_{\dagger}(1+e)\|\beta\|_{\infty}\widetilde{K}_0 \quad   \mbox{and}\quad  0\leq D^n\leq a_{\dagger}\|d\|_{\infty}\widetilde{K}_0
	 \,.
	 \]  
\end{remark}

\begin{lem}
Under the hypotheses of the previous Lemma, there exists two constants $\widetilde{K}_1$ and $\widetilde{K}_2$, depending on $T$, $K_d$, $\gamma$, the bounds of the initial and boundary conditions, and those of $S^n_j$, $I^n_j$, and $p_j^n$, such that
\begin{equation}
	0 \leq P^n \leq \widetilde{K}_1 \quad \mbox{and} \quad 	0 \leq N_F^n + N_I^n \leq \widetilde{K}_2, \quad  \forall n \in \{0,1,\ldots, N\}.
\end{equation}
\end{lem}
\proof 

Combining Equations (\ref{eqq1aa}) and (\ref{eqq1ab})  gives  
$$
P^n+h^2\left(\sum_{j=1}^{M}\mu_j(S{}^n_j +I{}^n_j)+\sum_{j=1}^{M}d_j\frac{N_{I}^{n-1}}{K_d+B^{n-1}} I{}^n_j\right)= P^{n-1} +h(S{}^{n-1}_0 -S{}^{n-1}_M - I{}^{n-1}_M)
\,.
$$
Therefore,
 \(
0\leq P^n\leq    P^{n-1} +h(S{}^{n-1}_0 -S{}^{n-1}_M - I{}^{n-1}_M)
 \leq P^{n-1} +mh
\,.
\) 
By induction, it follows that
\[ 
0\leq P^{n} \leq P^0+n~m~h \leq P^0+N~m~h=P^0+m~T =:\widetilde{K}_1\,.
\] 
In the same way, set  $\displaystyle Q^n:= N_F^n + N_I^n\,. $ 
%
%
Adding Equations (\ref{eqq1ac}) and (\ref{eqq1ad}) gives
$$
Q^n\left(1+h\min(\mu_F,\mu_I)\right)+\frac{ h~\rho~D^{n-1}}{K_d+B^{n-1}}\left(\frac{N_I^{n-1}~N_I^n}{K}+N_F^{n-1}\right) = Q^{n-1}\left(1+ \frac{h~\rho~D^{n-1}}{K_d+B^{n-1}}\right) +h\gamma B^{n-1}
\,.
$$
Therefore,
\[
0\leq Q^n  \leq    Q^{n-1}\left(1+ \frac{h~\rho~D^{n-1}}{K_d+B^{n-1}}\right) +h\gamma B^{n-1}
 \leq  
   \left(
   1+h\frac{\rho a_{\dagger}\|d\|_{\infty}\widetilde{K}_0}{K_d} \right)Q^{n-1} +  h\gamma a_{\dagger}\widetilde{K}_0
\,.
\]
By induction, it follows that
\begin{eqnarray} 
0\leq Q^{n}& \leq & 
\left(
   1+h\frac{\rho a_{\dagger}\|d\|_{\infty}\widetilde{K}_0}{K_d} \right)^nQ^{0} +\frac{\left(
      1+h\frac{\rho a_{\dagger}\|d\|_{\infty}\widetilde{K}_0}{K_d} \right)^n-1}{  \rho  \|d\|_{\infty}  } \gamma  K_d
      \,,
      \no
      \\
      &\leq &
\left(
   1+h\frac{\rho a_{\dagger}\|d\|_{\infty}\widetilde{K}_0}{K_d} \right)^N~Q^{0} +\frac{\left(
      1+h\frac{\rho a_{\dagger}\|d\|_{\infty}\widetilde{K}_0}{K_d} \right)^N-1}{  \rho  \|d\|_{\infty}  } \gamma  K_d
     \,.
     \no
\end{eqnarray} 
Now, from the trivial inequality $\dsp (1+a)^{\alpha} \leq e^{a\alpha}$, it follows that  
$$ \left(
   1+h\frac{\rho a_{\dagger}\|d\|_{\infty}\widetilde{K}_0}{K_d} \right)^N \leq e^{Nh\frac{\rho a_{\dagger}\|d\|_{\infty}\widetilde{K}_0}{K_d}}
   = e^{T\frac{\rho a_{\dagger}\|d\|_{\infty}\widetilde{K}_0}{K_d}}
   \,.   
   $$
Finally, we are lead to:
\(
0\leq Q^{n}  \leq  e^{T\frac{\rho a_{\dagger}\|d\|_{\infty}\widetilde{K}_0}{K_d}}\left(Q^0+\frac{\gamma ~K_d}{\rho~\|d\|_{\infty}}   \right) -\frac{\gamma ~K_d}{\rho~\|d\|_{\infty}}:=\widetilde{K}_1  
\)
and the proof is complete.
\qed
\subsection{Consistency of the discretization}
In this section, we prove that the numerical scheme developed in the previous section 
converges to the solution of the model system~\eqref{Model diagram} as the discretization 
parameter tends to zero. 
To this end, we introduce the discrete \(L_h^{1}\) and \(L_h^{\infty}\) norms:
 for sequences  $v = \left(v^n\right)_{1\leq n \leq N}$ and 
$u = \left(u^n_j\right)_{1\leq n \leq N, \, 1\leq j\leq M}\,,\;$  we set
\begin{equation}
\|v\|_{\infty,h} := \max_{1\leq n \leq N} |v^n|\,,\quad  \|u^n\|_{1,h} := h\sum_{j=1}^M |u_j^n|\quad \mbox{and}\quad  \|u\|_{1,\infty, h} :=  h\max_{1\leq n \leq N} \sum_{j=1}^M |u_j^n|
\,.
\end{equation}
We now introduce the error vectors in the approximation denoted $\mathcal{E}$ and $\mathbf{e}$, defined by 
$$
\mathcal{E}^n_j = \big(~ S(jh , n h ) - S{}^n_j,\, I(jh , n h ) - I{}^n_j~\big)^T 
	  =: \bigl(  \mathcal{E}_S{}^{n}_j,   \mathcal{E}_I{}^{n}_j \bigr)^T
	\,,
$$
and
\[
\mathbf{e}^n
   = \bigl( N_F(nh)-N_F^{\,n},\;
            N_I(nh)-N_I^{\,n},\;
            P(nh)-P^{\,n} \bigr)^T
   =: \bigl( \mathbf{e}_F^{\,n},\; \mathbf{e}_I^{\,n},\; \mathbf{e}_P^{\,n} \bigr)^T.
\]
Set, 
$$
|\mathcal{E}^n_j| = |  \mathcal{E}_S{}^{n}_j |+ |  \mathcal{E}_I{}^{n}_j   |\quad  \mbox{and}\quad 
|\mathbf{e}^n|
   =   |\mathbf{e}_F^{\,n}|+ |\mathbf{e}_I^{\,n}|+ |\mathbf{e}_P^{\,n}| 
\,.
$$
\begin{remark}\label{rm0}
In the sequel, \(C\) denotes a generic constant whose exact value is irrelevant
for our purposes. The constant \(C\) may change from line to line and depends
only on $a_{\dagger}, \, T$, on  the bounds of \(\beta\), on the Lipschitz constant of \(\mathcal{W}\),
and on the bounds of the continuous and discrete solutions of the model system
\eqref{Model diagram}.
\end{remark}
\noindent
We have the following
\begin{theo}
\label{theo_conssitency}
In addition to the assumptions of Section~\ref{tata}, assume that the solutions of the model system~\eqref{Model diagram} and their first and second order derivatives exist and are essentially bounded on $[0,a_{\dagger})\times [0, T]$,  
then there exists a constant \(\dsp C_0^*>0\) (as in Remark~\ref{rm0}) such that,
\begin{equation}
\text{if }\quad  C_0^* h \le \frac{1}{2} \quad \text{then} \quad 
\| \mathbf{e} \|_{\infty,h}   \leq C_0^* h
\quad\text{and}\quad
\| \mathcal{E}\|_{1,\infty, h}   \leq C_0^* h \,.
\end{equation}
\end{theo}
\begin{remark}
The hypotheses of Theorem~\ref{theo_conssitency} are satisfied, for instance, if the nonlinear perturbation 
operator~$\mathcal{H}$ associated with the linear operator~$\mathcal{A}$ 
is locally Lipschitz continuous on the Sobolev space 
\(\dsp
W^{2,1}(0,a_{\dagger}) \times W^{2,1}(0,a_{\dagger}) \times \mathbb{R} \times \mathbb{R},
\)
with values in the same space. 
As in Proposition~\ref{Reg}, this holds, for example, when 
\(\dsp
m \in C^{2}_b([0,\infty)), \; 
\beta,\, d \in C^{2}_b([0,a_{\dagger})),
\)
and 
\(\dsp 
\pi\mu,\, \pi\mu',\, \pi(\mu')^2,\, \pi\mu^2,\, \pi\mu^3 \in L^1(0,a_{\dagger}).
\)
\end{remark}

\proof
We introduce the notation \(  D^2S(a,t)= \frac{1}{2}\begin{pmatrix}
1\\1
\end{pmatrix}^T\cdot H_S(a,t)\cdot \begin{pmatrix}
1\\1
\end{pmatrix}\) where \(\dsp H_S\) is the Hessian matrix of $S$.
By the Taylor's expansion formula of order two, one has,
\[
  S(jh, nh)-S((j-1)h,(n-1)h) =
  h\,\frac{\partial S}{\partial a}(jh, nh)+ h\,\frac{\partial S}{\partial t}(jh, nh)  
  - h^2\,D^2S(jh,nh)
  + h^2 \varepsilon(h).
\]
The first equation of the model system~\eqref{Model diagram} then gives 
\[
 \frac{S(jh, nh)-S((j-1)h,(n-1)h)}{h} =
   -(\Lambda(jh,nj) +\mu_j)S(jh,nh)  
 - h\,D^2S(jh,nh)+ h \varepsilon(h)\,.
\]
Therefore,
\begin{eqnarray}
	\dfrac{ \mathcal{E}_S{}^{n}_j  -    \mathcal{E}_S{}^{n-1}_{j-1}}{h } & = &      \dfrac{S(jh , nh ) -S((j-1)h , (n-1)h )}{h } 
	-\dfrac{ S{}_{j}^{n}- S{}_{j-1}^{n-1}}{h } \,,
	\no
	\\
	&= & -(\Lambda(jh,nj) +\mu_j)S(jh,nh) +  \left(\Lambda^{n-1}_{j} + \mu_{j}\right)S{}^{n}_j  - h\,D^2S(jh,nh)  + h\varepsilon(h)\,,
	\no
	\\
	&= &  -\mu_{j}  \mathcal{E}_S{}^{n}_j +\Lambda_{j}^{n-1} S_{j}^{n}  - \Lambda(jh , nh )S(jh , nh )  - h\,D^2S(jh,nh)+ h\varepsilon(h)\,,
	\no
	\\
	& = &  -\mu_{j} \mathcal{E}_S{}^{n}_j +\Lambda_{j}^{n-1}  S_{j}^{n} +\Lambda_{j}^{n-1}  S(jh , nh )
  - \Lambda_{j}^{n-1}  S(jh , nh )
  	\no
  	\\
  	&&\qquad \qquad \qquad \qquad \qquad  - \Lambda(jh , nh )S(jh , nh )  - h\,D^2S(jh,nh) + h\varepsilon(h)\,,
	\no
	\\
	& = & -(\Lambda_{j}^{n-1} +\mu_{j})  \mathcal{E}_S{}^{n}_j - S(jh , nh ) \left( \Lambda(jh , nh ) -\Lambda_{j}^{n-1}  \right)
	 - h\,D^2S(jh,nh)+ h\varepsilon(h) 
	\,.
	\no
\end{eqnarray}
Therefore, we have the identity:
\[
	 \mathcal{E}_S{}^{n}_j\left(1+h(\Lambda^{n-1}_{j}+\mu_j)\right) =  \mathcal{E}_S{}^{n-1}_{j-1}-hS(jh,nh)(\Lambda(jh,nh)-\Lambda^{n-1}_{j-1} ) - h^2\,D^2S(jh,nh)+h^2\varepsilon(h)\,,
\]
from where we obtain the following estimate
\begin{eqnarray}
	\lvert  \mathcal{E}_S{}^{n}_j\rvert &\leq& 	\lvert  \mathcal{E}_S{}^{n}_j\rvert\left(1+h(\Lambda^{n-1}_{j}+\mu_j)\right)
	\,,
	\no
	\\
	\label{18V24a}
	& \leq&   \lvert \mathcal{E}_S{}^{n-1}_{j-1}\rvert + hS(jh,nh)\lvert \Lambda(jh,nh)-\Lambda^{n-1}_{j-1}\rvert+h^2\left(|D^2S(jh,nh)|+\varepsilon(h)\right)\,.
\end{eqnarray}
The expression $ \lvert \Lambda(jh , nh )-\Lambda_{j}^{n-1}   \lvert $ is estimated as  follows 
:
\begin{eqnarray}                              
	\lvert \Lambda(jh , nh ) -\Lambda_{j}^{n-1}   \lvert &= &  \beta_{j} \left|  \mathcal{W}(P(nh ),N_F(nh )) -   \mathcal{W}(P^{n-1},N_{F}^{n-1}) )\right|
	\,,
	\no
	\\
	&\leq &   
	C \left(|N_F(nh )-N_{F}^{n-1}| + |P(nh )- P^{n-1}| \right)
	\,,
	\no
	\\
	&\leq &  
	C \big(         
	|N_F(nh )-N_F {(n-1)h}| + |N_F((n-1)h)-N_{F}^{n-1}|
	\no
	\\
	&&
	\qquad\qquad     +\; |P(nh )- P ((n-1)h)|+ |P ((n-1)h)- P^{n-1}| \big)
	\,,
	\no
	\\
	&\leq &   C \left(  \lvert\mathbf{e}_F^{n-1}\lvert+ \lvert\mathbf{e}_P^{n-1}\lvert +   h  \right)   
	\,.
	\no
\end{eqnarray} 
In the last inequality, we applied the finite-increment formula
\(\dsp 
f(b)-f(a) = (b-a) f'\bigl(a+\delta(b-a)\bigr),
\; \delta\in(0,1),
\)
together with the assumption that the exact solutions and their derivatives are uniformly bounded.
This yields the estimate (recall that \(   h\sum_{j=1}^{M}1=Mh= a_{\dagger} \)):
\begin{equation}
\label{ineq0}
	h\sum_{j=1}^{M}\lvert \Lambda(jh , nh ) -\Lambda_{j}^{n-1}   \lvert \leq   C \left(\lvert\mathbf{e}_F^{n-1}\lvert+  \lvert\mathbf{e}_P^{n-1}\lvert +   h  \right) \,.
\end{equation}
Multiplying inequality~\eqref{18V24a} by \(h\) and summing over \(j=1,\dots,M\) gives

\begin{equation}
	\label{18III24b}
 \| \mathcal{E}_S{}^{n}\|_{1,h}	  
	\leq
		  \| \mathcal{E}_S{}^{n-1}\|_{1,h}+ Ch\left( \lvert \mathbf{e}_P^{n-1}\rvert	+ h\right)\,.
\end{equation}
Note that we have used the fact that
\(\dsp  
h \sum_{j=1}^{M} \lvert \mathcal{E}_S{}^{\,n-1}_{j-1}\rvert
   \le
h \sum_{j=1}^{M} \lvert \mathcal{E}_S{}^{\,n-1}_{j}\rvert,
\)  
since \(\dsp \lvert \mathcal{E}_S{}^{\,n-1}_{0}\rvert = 0  
\).

\noindent 
We repeat the previous argument for $I^n_{j}$. We have 
\begin{eqnarray} \dfrac{I(jh, nh)- I((j-1)h, (n-1)h)}{h}  &=&  \frac{\partial I }{ \partial{t}} (jh, nh)+  \frac{\partial I }{ \partial{a}} (jh, nh)
 - h \,D^2I(jh,nh)	+ h\varepsilon(h)\,,
	\no
	\\
	&=&  \Lambda(jh,nh)S(jh,nh) -\mu_j I(jh,nh)
\no%
	\\
	&&\qquad  - h \,D^2S(jh,nh)
	 - \frac{d_{j}N_I(nh)I(jh,nh)}{K_d+B(nh)} + h\varepsilon(h)\,,
	\no
\end{eqnarray}  
and,
\begin{eqnarray}
\dfrac{ \mathcal{E}_I{}^{n}_j  -    \mathcal{E}_I{}^{n-1}_{j-1}}{h } &= &      \dfrac{I(jh , nh ) -I((j-1)h , (n-1)h )}{h } 
	-\dfrac{ I{}_{j}^{n}- I{}_{j-1}^{n-1}}{h } 
	\,,
	\no
	\\
	&= &  -\mu_{j} \mathcal{E}_I{}^{n}_j -\Lambda_{j}^{n-1} S_{j}^{n}  +\Lambda(jh , nh )S(jh , nh )
	\no
	\\
	&&
	\qquad	\qquad	\qquad  - h\,D^2S(jh,nh)  +\frac{d_{j}N_I^{n-1}I{}^{n}_j}{K_d+B^{n-1}} 
	- \frac{d_{j}N_I(nh)I(jh,nh)}{K_d+B(nh)}+ h\varepsilon(h)
	\,,
\no
\\
	& = &  -\mu_{j} \mathcal{E}_I{}^{n}_j +  \mathcal{E}_S{}^{n}_j\Lambda_{j}^{n-1}  +S(jh , nh ) \left( \Lambda(jh , nh ) -\Lambda_{j}^{n-1}  \right)
		 	\no
		 	\\
		 	&&
		 	\qquad\qquad\qquad 
	 +\frac{d_{j}N_I^{n-1}I{}^{n}_j}{K_d+B^{n-1}} 
	  - \frac{ d_{j}N_I(nh)I(jh,nh)}{K_d+B(nh)} - h\,D^2S(jh,nh)+ h\varepsilon(h)
	  \,,
\no
\\
&=&
-\mu_{j} \mathcal{E}_I{}^{n}_j +  \mathcal{E}_S{}^{n}_j\Lambda_{j}^{n-1}  +S(jh , nh ) \left( \Lambda(jh , nh ) -\Lambda_{j-1}^{n-1}  \right)    
	-  \frac{d_j N_I^{n-1}\mathcal{E}_I{}^{n}_j }{K_d+B^{n-1}} 
\no
	\\
	&& \qquad 	\qquad    -\; d_j I(jh , nh )\left(\frac{N_I(nh)}{K_d+B(nh)}-\frac{N_I^{n-1}}{K_d+B^{n-1}}\right)  
	- h\,D^2S(jh,nh)  + h\varepsilon(h)
	\,,
\no
\\
 &= & 	
	 -\mu_{j} \mathcal{E}_I{}^{n}_j +  \mathcal{E}_S{}^{n}_j\Lambda_{j}^{n-1}  +S(jh , nh ) \left( \Lambda(jh , nh ) -\Lambda_{j-1}^{n-1}  \right)    
	 	-  \frac{d_j N_I^{n-1}\mathcal{E}_I{}^{n}_j }{K_d+B^{n-1}} - h\,D^2S(jh,nh)
\no
	\\
	&&   -\; \frac{K_d d_j I(jh , nh )\left(N_I(nh)-N_I^{n-1} \right)}{(K_d+B(nh))(K_d+B_i^{n-1})}    
 -\; \frac{d_j I(jh , nh )\left(B_i^{n-1}N_I(nh) -B(nh)N_I^{n-1}\right)}{(K_d+B(nh))(K_d+B^{n-1})}    + h\varepsilon(h)
 \,,
\no
	\\
	& = &  -\mu_{j} \mathcal{E}_I{}^{n}_j +  \mathcal{E}_S{}^{n}_j\Lambda_{j}^{n-1}  +S(jh , nh ) \left( \Lambda(jh , nh ) -\Lambda_{j-1}^{n-1}  \right)    
		 	-  \frac{d_j N_I^{n-1}}{K_d+B^{n-1}}\mathcal{E}_I{}^{n}_j  
\no
	\\
	&& \qquad -\; \frac{K_dd_j I(jh , nh )(N_I(nh )-N_I{(n-1)h}) }{(K_d+B(nh))(K_d+B^{n-1})}
-  \frac{K_dd_j I(jh , nh )\mathbf{e}_I^{n-1}}{(K_d+B(nh)\|)(K_d+B^{n-1})}
\no
	\\
	&& \qquad -\;\frac{d_j I(jh , nh )\left(B^{n-1}N_I(nh) -B(nh)N_I^{n-1}\right)}{(K_d+B(nh))(K_d+B^{n-1})} - h\,D^2S(jh,nh)   + h\varepsilon(h)
\no
	\,.
\end{eqnarray}
Therefore, we have the identity:
\begin{eqnarray}
	\label{18V24b1}
	\mathcal{E}_I{}^{n}_j\left(1+h\left(\frac{d_jN_I^{n-1}}{K_d+B^{n-1}} +\mu_j\right)\right) &=& \mathcal{E}_I{}^{n-1}_{j-1} + h \mathcal{E}_S{}^{n}_j\Lambda_{j}^{n-1}   +hS(jh,nh)(\Lambda(jh,nh)-\Lambda^{n-1}_{j-1} ) 
	\no
	\\
	& - &  \frac{h~K_dd_j I(jh , nh )(N_I(nh )-N_I{(n-1)h})}{(K_d+B(nh))(K_d+B^{n-1})} 
	\label{ineq1}
\no
	\\
	& - & \frac{h~K_dd_j I(jh , nh )\mathbf{e}_I^{n-1}}{(K_d+B(nh))(K_d+B^{n-1})} + h^2\left(D^2S(jh,nh)+\varepsilon(h)\right)
	\no
	\\
	& - & \frac{h~d_j I(jh , nh )\left(B^{n-1}N_I(nh) - B(nh)N_I^{n-1}\right)}{(K_d+B(nh))(K_d+B^{n-1})}  
	\label{ineq3} 
	\,.
\end{eqnarray}
Next, 
\begin{eqnarray}
B^{n-1}N_I(nh) -B(nh)N_I^{n-1} &=&  B^{n-1}\left(N_I(nh)-N_I((n-1)h)\right)+B^{n-1}N_I((n-1)h)  - B(nh)N_I^{n-1}
\,,
\no
\\
&=&
B^{n-1}\left(N_I(nh)-N_I((n-1)h)\right)+B^{n-1}\left(N_I((n-1)h)-N_I^{n-1}\right)
\no
\\
&& \qquad\qquad\qquad\qquad +  \left(B^{n-1} - B(nh)\right) N_I^{n-1}\,.
\label{ineq0}
\end{eqnarray}
Note that 
\(\displaystyle 
B^{n-1} - B(nh) = B^{n-1} - B((n-1)h) + B((n-1)h) - B(nh)\,,
\)
with,  
\begin{eqnarray}
 B^{n-1} - B((n-1)h)  &=&  h \sum_{j=1}^{M}I{}^{n-1}_j-\int_0^{a_{\dagger}}I(a,(n-1)h)da 
 \,,
 \no
 \\
 &=&  h \sum_{j=1}^{M}I{}^{n-1}_j -\int_0^{a_{\dagger}}I(a,(n-1)h)da+h\sum_{j=1}^{M}I(jh,(n-1)h) -h\sum_{j=1}^{M}I(jh,(n-1)h)  
 \,,
 \no
 \\
  &=&
- h	\sum_{j=1}^{M}\mathcal{E}_I{}^{n-1}_j  + h\sum_{j=1}^{M}I(jh,(n-1)h)-\int_0^{a_{\dagger}}I(a,(n-1)h)da
\no
  \,.
\end{eqnarray}
By the classical error estimate for the quadrature formula, one has
\[
\left| h\sum_{j=1}^{M}I(jh,(n-1)h)-\int_0^{a_{\dagger}}I(a,(n-1)h)da\right|
\;\le\;
C h^{2}\, 
\]
and since  \(\displaystyle |B((n-1)h) - B(nh)|\leq Ch\)  we have, 
\begin{equation}\label{ineq5}
|B^{n-1} - B(nh)| \leq \lVert \mathcal{E}_I{}^{n-1} \rVert_{1,h} + Ch\,.
\end{equation}
It then follows  from~(\ref{ineq0}) that 
\begin{equation}\label{ineq2}
 \left|B^{n-1}N_I(nh) -B(nh)N_I^{n-1}\right| \leq  C\left( \lVert \mathcal{E}_I{}^{n-1} \rVert_{1,h} +\mathbf{e}_I^{n-1} +h\right)\,.
\end{equation}
From the expression  of $ \mathcal{E}_I{}^{n}_j$ given by \eqref{18V24b1} we are lead to the following estimate:
\begin{equation}
\label{ineq4}
\lVert \mathcal{E}_I{}^{n}\rVert_{1,h} \leq \lVert \mathcal{E}_I{}^{n-1} \rVert_{1,h} + Ch \left( \lVert \mathcal{E}_I{}^{n-1} \rVert_{1,h}+ \lVert \mathcal{E}_S{}^{n}   \rVert_{1,h}  +   \lvert\mathbf{e}_I^{n-1}\lvert +\lvert\mathbf{e}_F^{n-1}\lvert +\lvert\mathbf{e}_P^{n-1}\lvert +   h    \right)
\,.
\end{equation}
We now proceed with the analysis of the components of the error vector
\(\mathbf{e}^n\).
As in the previous step, the Taylor expansion of \(N_F\) up to second order
in a neighborhood of \(nh\) reads

$$
N_F(nh)-N_F((n-1)h) = h\dot{N}_F(nf)-\frac{1}{2} h^2   \ddot{N}_F(nh) + h^2\epsilon(h)\,, 
$$
and using the ODE satisfies by $N_F$ yields, 
$$
N_F(nh)-N_F((n-1)h) = -\alpha h C(nh) \mathcal{W}(P(nh),N_{F}(nh))+h\gamma B(nh)-h\mu_FN_F -h^2\left(\frac{1}{2}  \ddot{N}_F(nh)-\epsilon(h)\right)
\,.
$$
On the other hand 
$$
 N_F^{n}- N_F^{n-1} = -\alpha h C^{n-1} \mathcal{W}(P^{n-1}, N_F^{n}) +h\gamma B^{n-1} - h\mu_F N_F^{n}
 \,.
$$
Therefore
\begin{eqnarray}
 N_F(nh)-N^n_F &=& N_F((n-1)h)-N_F^{n-1}   -\alpha h \bigg(C(nh)\mathcal{W}(C(nh)P(nh),N_{F}(nh))-C^{n-1} \mathcal{W}(P^{n-1}, N_F^{n}) \bigg)
\no
\\
&&  
+\,h\gamma \left(  B(nh)-B^{n-1}\right)-h\mu_F\big(N_F(nh)-N^n_F \big)-h^2\left(\frac{1}{2}  \ddot{N}_F(nh)-\epsilon(h)\right)
\,,
\no
\end{eqnarray}
i.e.
\begin{eqnarray}
(1+h\mu_F)\mathbf{e}_F^n&=& \mathbf{e}_F^{n-1}  -\alpha h \bigg(C(nh)\mathcal{W}(P(nh),N_{F}(nh))-C^{n-1} \mathcal{W}(P^{n-1}, N_F^{n}) \bigg)
\no
\\
&&
\qquad \qquad \qquad \qquad 
+ \;h\gamma \left(  B(nh)-B^{n-1}\right) -h^2\left(\frac{1}{2}  \ddot{N}_F(nh)-\epsilon(h)\right) 
\,.
\no
\end{eqnarray}
Now from
\begin{eqnarray}
C(nh)\mathcal{W}(P(nh),N_{F}(nh))-C^{n-1} \mathcal{W}(P^{n-1}, N_F^{n}) &= &
C(nh)\big(\mathcal{W}(P(nh),N_{F}(nh))-\mathcal{W}(P^{n-1}, N_F^{n})\big)
\no
\\
 && 
+ \big(C(nh)-C^{n-1}\big) \mathcal{W}(P^{n-1}, N_F^{n})\,,
\no
\end{eqnarray}
we see that,
\begin{eqnarray}
\lvert C(nh)\mathcal{W}(P(nh),N_{F}(nh))-C^{n-1} \mathcal{W}(P^{n-1}, N_F^{n})\rvert &\leq &
C\left( \lvert P(nh)-P^{n-1}\rvert + \lvert N_{F}(nh)-N_F^{n}\rvert + \lvert C(nh)-C^{n-1}\rvert \right) \,,
\no
\\
&\leq & C\left(\lvert P(nh)-P((n-1)h)\rvert + \lvert \mathbf{e}_P^{\,n-1}\rvert+ \lvert\mathbf{e}_F^{\,n}\rvert \right)
\no
\\
&& +\; C \lvert C(nh)-C((n-1)h)\rvert +\lvert C((n-1)h)-  C^{n-1}\rvert 
\,.
\no
\end{eqnarray}
But, 
\begin{eqnarray} 
\lvert C((n-1)h)-  C^{n-1}\rvert  &=&  \lvert\int_0^{a_{\dagger}}\beta(a)(S(a,(n-1)h)+eI(a,(n-1)h))da - h \sum_{j=1}^{M}\beta_j(S{}^{n-1}_j+eI{}^{n-1}_j)\rvert 
\,,
\no
 \\
 &=&
 \lvert \int_0^{a_{\dagger}}\beta(a)(S(a,(n-1)h)+eI(a,(n-1)h))da  
 \no
 \\
 &&
 \qquad \qquad \qquad -\;  h \sum_{j=1}^{M}\beta_j \left(S(jh,(n-1)h) +e I(jh,(n-1)h)\right) \rvert     
  \no
  \\
  &&
  + h \lvert\left( \sum_{j=1}^{M}\beta_j\left((S{}^{n-1}_j-S(jh,(n-1)h)+e\left(I{}^{n-1}_j-I(jh,(n-1)h)\right)  \right)  \right)\rvert
  \,,
  \no
  \\
 & \leq&
C\left(h^2+ \lVert \mathcal{E}_S{}^{n-1}+\lVert \mathcal{E}_I{}^{n-1} \rVert_{1,h}    \right)
\,.
\no
\end{eqnarray}
Thus,
\begin{eqnarray}
\lvert C(nh)\mathcal{W}(P(nh),N_{F}(nh))-C^{n-1} \mathcal{W}(P^{n-1}, N_F^{n})\rvert 
&\leq & C\left(h+ \lvert \mathbf{e}_P^{\,n-1}\rvert+ \lvert\mathbf{e}_F^{\,n}\rvert + \lVert \mathcal{E}_S{}^{n-1}+\lVert \mathcal{E}_I{}^{n-1} \rVert_{1,h}   \right)
\,.
\no
\end{eqnarray}
Therefore we obtain the following estimate (recall,
\( |B^{n-1} - B(nh)| \leq \lVert \mathcal{E}_I{}^{n-1} \rVert_{1,h} + Ch\,,  \) see \eqref{ineq5}):
\begin{equation}\label{ineq6}
\lvert
\mathbf{e}_F^n   
\rvert
\leq 
\lvert
\mathbf{e}_F^{n-1}   
\rvert
+
Ch\left(  \lvert \mathbf{e}_P^{\,n-1}\rvert + \lVert \mathcal{E}_S{}^{n-1}\rVert+\lVert \mathcal{E}_I{}^{n-1} \rVert_{1,h}  \right)+Ch \lvert\mathbf{e}_F^{\,n}\rvert +Ch^2 
\,. 
\end{equation}
Now we continue with  the estimate for  \(\mathbf{e}_I^n  \).  
Again,  the Taylor expansion of \(N_I\) up to second order
in a neighborhood of \(nh\) gives,
$$
N_I(nh)-N_I((n-1)h) = h\dot{N}_I(nh)-\frac{1}{2} h^2 \ddot{N}_I(nh) + h^2\epsilon(h)\,, 
$$
and using the ODE satisfies by $N_I$ yields, 
\begin{eqnarray}
N_I(nh)-N_I((n-1)h) &=& \alpha h C(nh) \mathcal{W}(P(nh),N_{F}(nh)) -h\mu_IN_I(nh)
\no
\\
&&
\qquad   \;+
h\rho \frac{D(nh)N_I(nh)}{K_d+B(nh)}\left(1-\frac{N_I(nh)}{K}\right)
 -h^2\left(\frac{1}{2}  \ddot{N}_I(nh)-\epsilon(h)\right)
\,.
\no
\end{eqnarray}
On the other hand 
$$
 N_I^{n}- N_I^{n-1} =  \alpha h C^{n-1} \mathcal{W}(P^{n-1}, N_F^{n})  - h\mu_I N_I^{n} +h\rho \frac{N_I^{n-1}D^{n-1}}{K_d+B^{n-1}}
 \left(1-\frac{N_I^n}{K}\right)
 \,.
$$
Therefore,
\begin{eqnarray}
 N_I(nh)-N^n_I &=& N_I((n-1)h)-N_I^{n-1}   +\alpha h \bigg(C(nh)\mathcal{W}(P(nh),N_{F}(nh))-C^{n-1} \mathcal{W}(P^{n-1}, N_F^{n}) \bigg)
\no
\\
&&  
-h\mu_I\left(N_I(nh)-N^n_I \right) +h\rho\underbrace{\left(  \frac{D(nh)N_I(nh)}{K_d+B(nh)}\left(1-\frac{N_I(nh)}{K}\right) - \frac{N_I^{n-1}D^{n-1}}{K_d+B^{n-1}}
 \left(1-\frac{N_I^n}{K}\right)  \right)}_{:= R^n}
\no
\\
&&
-h^2\left(\frac{1}{2}  \ddot{N}_I(nh)-\epsilon(h)\right)
\,,
\no
\end{eqnarray}
i.e.
\[
(1+h\mu_I)\mathbf{e}_I^n =  \mathbf{e}_I^{n-1}  +\alpha h \big(C(nh)\mathcal{W}(P(nh),N_{F}(nh))-C^{n-1} \mathcal{W}(P^{n-1}, N_F^{n}) \big)
+h\rho R^n  -h^2\left(\frac{1}{2}  \ddot{N}_I(nh)-\epsilon(h)\right) 
\,.
\]
From the preceding steps, the remaining task is to estimate the term \(R^n\); this is done as follows:
$$
R^n= 
 \underbrace{ \frac{D(nh)N_I(nh)}{K_d+B(nh)}  - \frac{N_I^{n-1}D^{n-1}}{K_d+B^{n-1}}}_{=: R^n_1}
- \underbrace{\left(\frac{D(nh)N_I(nh)}{K_d+B(nh)} \frac{N_I(nh)}{K}  - \frac{N_I^{n-1}D^{n-1}}{K_d+B^{n-1}}
  \frac{N_I^n}{K} \right)}_{:=R^n_2}  
  \,. 
$$
We have
$$
R_1^n= \underbrace{\frac{D(nh)(N_I(nh)-N_I^{n-1})}{K_d+B(nh)}}_{=: R^n_{11}} + 
N_I^{n-1}\underbrace{ \frac{K_d(D(nh)-D^{n-1})+(D(nh)B^{n-1}-B(nh)D^{n-1})}{(K_d+B(nh))(K_d+B^{n-1})}}_{=:R^n_{12}} 
\,.
$$
Form the identity \( N_I(nh)-N^{n-1}=   N_I(nh)-N_I((n-1)h+  N_I((n-1)h-N^{n-1} \) we have 
$$
\lvert R_{11}^n\rvert \leq C\left(h+\lvert \mathbf{e}_N^{n-1}\rvert \right) 
\,.
$$
The same computation as those we did when estimating   \(  B^{n-1} - B(nh)   \) lead to the following (see \eqref{ineq5}):
\begin{equation}
\label{inqe6}
  |D(nh)-D^{n-1} \leq \lVert \mathcal{E}_I{}^{n-1} \rVert_{1,h} + Ch\,.   
\end{equation}
From the identity
\(
D(nh)B^{n-1}-B(nh)D^{n-1} =  (D(nh)-D^{n-1})B^{n-1} + D^{n-1}\left( B^{n-1}-B(nh)\right)
\) and the inequalities \eqref{ineq5} and  \eqref{ineq6}, we have
\begin{equation}
\label{inqe7}
  \lvert D(nh)B^{n-1}-B(nh)D^{n-1} \rVert \leq \lVert \mathcal{E}_I{}^{n-1} \rVert_{1,h} + Ch\,,   
\end{equation}
which shows that
\begin{equation}
\label{inqe8}
  \lvert R^n_{1} \rVert \leq   C\left( \lvert \mathbf{e}_N^{n-1}\rvert + \mathcal{E}_I{}^{n-1} \rVert_{1,h} +  h\right)\,.   
\end{equation}
We continue with the estimate of  $R^n_2$. We have 
\[
R^n_2 = \underbrace{ \frac{D(nh)N_I(nh)}{K_d+B(nh)} \frac{N_I(nh)-N_I^n}{K} }_{=:R_{21}}     +\frac{N_I^n}{K}\underbrace{\left( \frac{D(nh)N_I(nh)}{K_d+B(nh)}- \frac{N_I^{n-1}D^{n-1}}{K_d+B^{n-1}} 
   \right) }_{=:R_{22}}
\,.
\]
Note that 
\(
\lvert R_{21}\rvert \leq C\mathbf{e}^n_I
\)
. Proceeding further, one has 
$$
R_{22}=   \frac{K_d\left( D(nh)N_I(nh)-N_I^{n-1}D^{n-1}\right) +D(nh)N_I(nh)B^{n-1}-B(nh)N_I^{n-1}D^{n-1}  }{(K_d+B(nh))(K_d+B^{n-1})} 
\,. 
$$
Following the same reasoning as above, one finds that,
\[
R_{22} \leq C\left(\lvert \mathbf{e}_N^{n-1}\rvert +\lVert \mathcal{E}_I{}^{n-1} \rVert_{1,h} +  h\right)\,,   
\]
leading to 
\[
  \lvert R^n_{2} \rvert \leq   C\left( \lvert \mathbf{e}_N^{n-1}\rvert + \lVert\mathcal{E}_I{}^{n-1} \rVert_{1,h} +  h\right)\,,   
\]
and then, 
\begin{equation}
\label{inqe9}
  \lvert R^n  \rvert \leq   C\left( \lvert \mathbf{e}_N^{n-1}\rvert + \lVert\mathcal{E}_I{}^{n-1} \rVert_{1,h} +  h\right)\,.   
\end{equation}
We are therefore led to 
\begin{equation}\label{ineq10}
\lvert
\mathbf{e}_I^n   
\rvert
\leq 
\lvert
\mathbf{e}_I^{n-1}   
\rvert
+
Ch\left(  \lvert \mathbf{e}_P^{\,n-1}\rvert+ \lvert \mathbf{e}_N^{\,n-1}\rvert + \lVert \mathcal{E}_S{}^{n-1}\rVert_{1,h}+\lVert \mathcal{E}_I{}^{n-1} \rVert_{1,h}  \right)+Ch \lvert\mathbf{e}_F^{\,n}\rvert +Ch^2 
\,. 
\end{equation}
Next, we continue with the analysis of $\mathbf{e}_P^n$. We have 
\[
\mathbf{e}_P^n =  P(nh)-P^{\,n} =
  \int_0^{a_{\dagger}} p(a,nh)da -  h \sum_{j=1}^{M}p(jh,nh) +  h \sum_{j=1}^{M}(p(jh,nh)-p_j^n) 
  \,.
\]
Now, from \eqref{18III24b} and \eqref{ineq4}, we have 
\begin{eqnarray}
\lvert  p(jh,nh)-p_j^n \rvert  & \leq & 
 \lvert \mathcal{E}_S{}^{n}_j  \rvert +\lvert \mathcal{E}_I{}^{n}_j  \rvert
\,,
\no
\\
&\leq &
\| \mathcal{E}_S{}^{n-1}\|_{1,h}+\lVert \mathcal{E}_I{}^{n-1} \rVert_{1,h} + Ch \bigg( \lVert \mathcal{E}_I{}^{n-1} \rVert_{1,h}+ \lVert \mathcal{E}_S{}^{n}   \rVert_{1,h} 
\no
+   \lvert\mathbf{e}_I^{n-1}\lvert +   \lvert\mathbf{e}_F^{n-1}\lvert +\lvert\mathbf{e}_P^{n-1}\lvert +   h    \bigg)
\,,
\end{eqnarray}
from where we obtain the following
\begin{eqnarray}
\lefteqn{h\sum_{j=1}^{M}\lvert  p(jh,nh)-p_j^n \rvert  \leq }     \nonumber\\
 && C\left( \| \mathcal{E}_S{}^{n-1}\|_{1,h}+\lVert \mathcal{E}_I{}^{n-1} \rVert_{1,h}\right) + Ch \bigg( \lVert \mathcal{E}_I{}^{n-1} \rVert_{1,h}+ \lVert \mathcal{E}_S{}^{n}   \rVert_{1,h} 
+  \lvert\mathbf{e}_F^{n-1}\lvert +   \lvert\mathbf{e}_I^{n-1}\lvert +\lvert\mathbf{e}_P^{n-1}\lvert +   h    \bigg)
\,.
\no
\end{eqnarray}
Therefore, using again the classical estimate for the quadrature formula one obtains the following estimate:
\begin{equation}
\label{ineq11}
\lvert \mathbf{e}_P^n\rvert \leq 
C\left( \| \mathcal{E}_S{}^{n-1}\|_{1,h}+\lVert \mathcal{E}_I{}^{n-1} \rVert_{1,h}\right) + Ch \bigg( \lVert \mathcal{E}_I{}^{n-1} \rVert_{1,h}+ \lVert \mathcal{E}_S{}^{n}   \rVert_{1,h} 
+  \lvert\mathbf{e}_F^{n-1}\lvert +   \lvert\mathbf{e}_I^{n-1}\lvert +\lvert\mathbf{e}_P^{n-1}\lvert +   h    \bigg)
\,.
\end{equation}
Let \(C^{*}\) denote the largest of the constants \(C\) obtained in the preceding estimates (see Remark~\ref{rm0}).
Now set
 \[ U^n= \| \mathcal{E}_S{}^{n}\|_{1,h} +\| \mathcal{E}_I{}^{n}\|_{1,h} + \lvert \mathbf{e}_P^n\rvert +\lvert \mathbf{e}_F^n\rvert +\lvert \mathbf{e}_I^n\rvert \,.\]
  Adding Inequalities~\eqref{18III24b}, \eqref{ineq4}, \eqref{ineq6},
\eqref{ineq10}, and \eqref{ineq11} gives
\begin{equation}
(1-C^*h)U^n \leq (1+C^*h)U^{n-1} + C^*h^2
\no
\,,
\end{equation}
which leads to the elementary estimate:
\[
U^n \leq \left(\frac{1+C^*h}{1-C^*h}
\right)^n U^0+\frac{1}{2}h(1-C^*h)\left( \left(\frac{1+C^*h}{1-C^*h}
\right)^n-1 \right)
\,,
\]
which implies that (recall $N=T/h$)
\[
U^n \leq \left(\frac{1+C^*h}{1-C^*h}
\right)^{\frac{T}{h}} U^0+\frac{1}{2}h(1-C^*h)\left( \left(\frac{1+C^*h}{1-C^*h}
\right)^{\frac{T}{h}}-1 \right)\,.
\]
Note that 
 \[ U^0 
= \lvert \mathbf{e}_P^0\rvert= \lvert \int_0^{a_\dagger}S_0(a)+I_0(a) da -h\sum_{j=1}^{M}S_0(a_j)+I_0(a_j)\rvert
  \,.\]
Applying once more the standard error bound for the quadrature formula, one finds that \(0\leq U^0 \leq C^*h\) and finally  
\begin{equation}
\label{ineq13}
U^n \leq C^*h\left(2\left(\frac{1+C^*h}{1-C^*h}
\right)^{\frac{T}{h}}  -1 \right)  \,.
\end{equation}
Note that \(\displaystyle \lim_{x\rightarrow 0} \left(\frac{1+x}{1-x}
\right)^{\frac{C^*T}{x}} = e^{2C^*T}  \) which shows that the function \( \displaystyle x\mapsto \left(\frac{1+x}{1-x}
\right)^{\frac{C^*T}{x}}  \) is continuous  on the compact set $[0,\frac{1}{2}]$ and consequently,  by the Weierstrass's Theorem it is a bounded function. Therefore, there exists a constant \(C^*_1> 0\) such that \(\displaystyle 0\leq 2\left(\frac{1+C^*h}{1-C^*h}
\right)^{\frac{T}{h}}  -1 \leq  C^*_1\). Finally by setting \(\displaystyle C_0^*= C^*C^*_1 \) we obtain that,
\[
0\leq U^n \leq C_0^*h \,,\;\forall 0\leq n\leq N \,.
\]
Hence, our semi-implicit Euler scheme converges in the discrete \(L_h^{\infty}\) norm with order \(O(h)\).
The proof is  complete.
\qed 

\subsection{Numerical simulations}
\subsubsection{Choice of parameters }
This section is devoted to the  numerical simulations of system \eqref{Model diagram}.  The numerical values of the parameters are summarized in the following table (Table \ref{t1}). 
\begin{table}[ht]
  \caption{Parameters of the model.}
  \begin{center}
    \begin{tabular}{llll}
      \hline 
      \textbf{Symbol} & \textbf{Description} & \textbf{Value}  & References \\
      \hline
      $m(t)$ & Recruitment of healthy plants & 300 g.day$^{-1}$ &  \cite{a300} \\
      $\beta(a)$ & Infection rate &  $ \beta_{\max} e^{-\frac{\left(a-a_\text{opt}\right)^2}{2\sigma^2}}$  & Assumed \\
      $e$ & Reinfection fraction & 0.0002 & \cite{a300} \\
      $\alpha$ & Infection conversion factor & 100 nematodes.g$^{-1}$  &  \cite{a100} \\
      $\gamma$ & Production rate of free nematodes & 1000 day$^{-1}$ &   Assumed \\
      $d(a)$ & Plant consumption rate  &  $d_{\max}e^{-\eta a}$ &  Assumed \\
      $K_d$ & Plant consumption half-saturation & 60 g &  \cite{a200}\\
      $\rho$ & Plant consumption conversion factor & 400 nematodes.g$^{-1}$ &   \cite{a200} \\ 
      $K$ & Nematode carrying capacity & 1000 nematodes.g$^{-1}$  &  \cite{a300} \\
      $\mu(a)$ & Mortality rate of plants &  $ 1/(a_{\dagger}-a)^3$  &   Assumed \\
      $\mu_{F}$ & Mortality rate of free nematodes & 0.0495 day$^{-1}$  &   \cite{a200} \\
      $\mu_{I}$ & Mortality rate of infesting nematodes & 0.045 day$^{-1}$ &  \cite{a200} \\
      \hline
    \end{tabular}    
  \end{center}
  \label{t1}
\end{table}
We will be using  the following initial and boundary  conditions:
\begin{align}\label{inc}
I_0(a) = 0\,,\quad  N_{F0} = 10^4\,, \quad  N_{I0} = 0 \quad \text{and} \quad
S_0(a) =
\left\{
\begin{array}{c}
0 \quad \text{if} \quad a>0
\\
 100 \quad \text{if} \quad a= 0
\end{array}
\right.
\,. 
\end{align}
The parameters in the infection and plants consumption rates formulae are chosen as: 
\begin{align}\label{par}
 a_\text{opt} = 5.5 ,\,\sigma_p = 2.5,\,d_{\max} = 10^{-4},\,\eta = 2.5.  
\end{align} 
 We will successively assign two distinct values to $\beta_{\text{max}}$, which will allow us to present our various graphs for the cases (i) $0 \leq \mathcal{N} < 1$ and (ii) $\mathcal{N} > 1$.
\subsubsection{Disease free equilibrium}
The choice of $\beta_{\text{max}} = 10^{-6}$ results in a basic reproduction number $\mathcal{N} = 0.13 < 1$, which is found to ensure the stability of the nematode-free equilibrium point. The subsequent plots corroborate these analytical findings.
 \begin{figure}[H]	
 \begin{center}
 \includegraphics[width=0.49\linewidth]{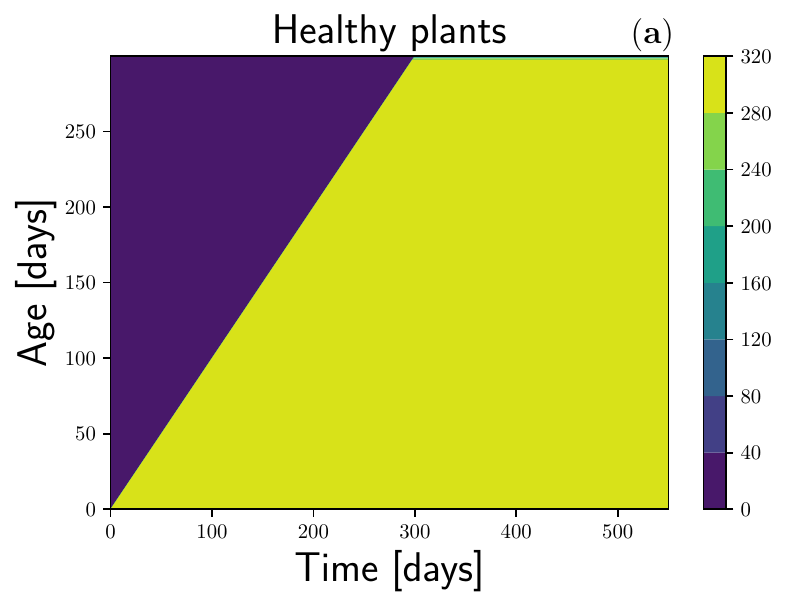}
 \includegraphics[width=0.49\linewidth]{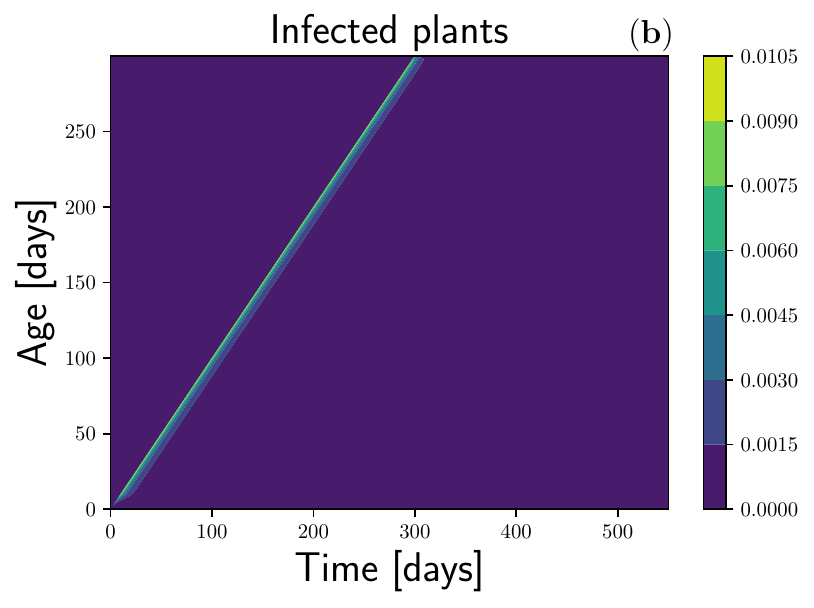}
 \includegraphics[width=0.49\linewidth]{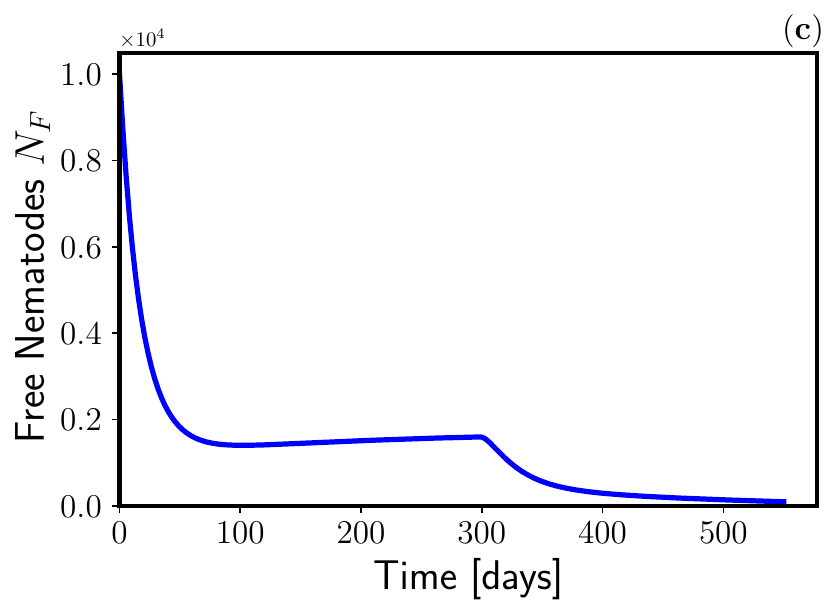}
 \includegraphics[width=0.49\linewidth]{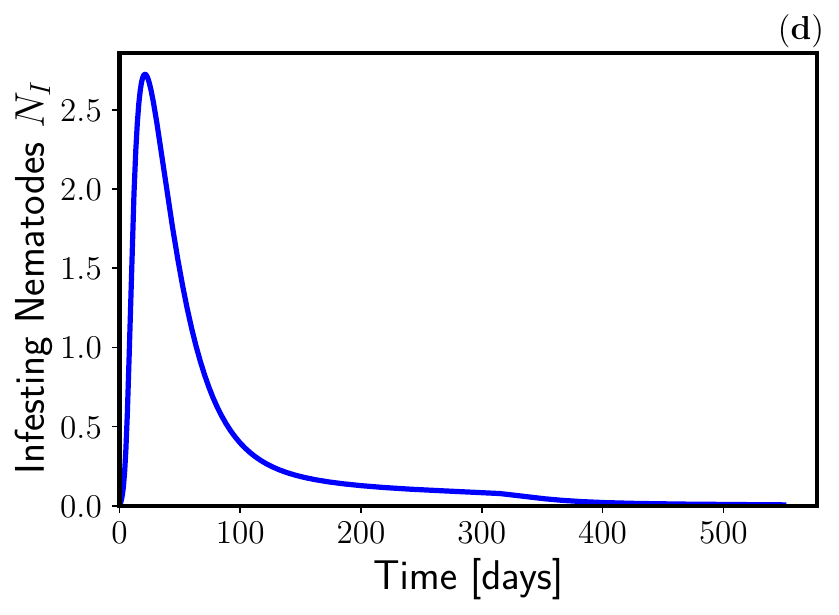}
 \caption{ Simulation in the case $0\leq \mathcal{N}<1$ of the dynamics of healthy plants $S(a,t)$, infected plants $I(a,t)$, free nematodes $N_F(t)$ and infesting nematodes $N_I(t)$ using the parameters given in Table~\ref{t1} and equations~\eqref{inc}--\eqref{par}.
 }
 \label{dfe}
 \end{center}
 \end{figure}
\vspace{-5mm}
\noindent
Subfigure \ref{dfe}(a) shows that the density of healthy plants remains dominant and tends toward a stable configuration, while Subfigure \ref{dfe}(b) illustrates a rapid decrease in the density of infected plants until it becomes negligible.
Likewise, Subfigures \ref{dfe}(c) and \ref{dfe}(d) indicate that the populations of free and infesting nematodes progressively decline over time, converging toward zero. These results confirm the local asymptotic stability of the disease-free equilibrium, characterized by the extinction of infections and the exclusive persistence of healthy plants.

\subsubsection{Endemic equilibrium}
The following plots correspond to the case $\mathcal{N} = 8.86 > 1$, obtained with the parameter choice $\beta_{\max} = 7 \times 10^{-5}$.

\begin{figure}[H]	
\begin{center}
\includegraphics[width=0.495\linewidth]{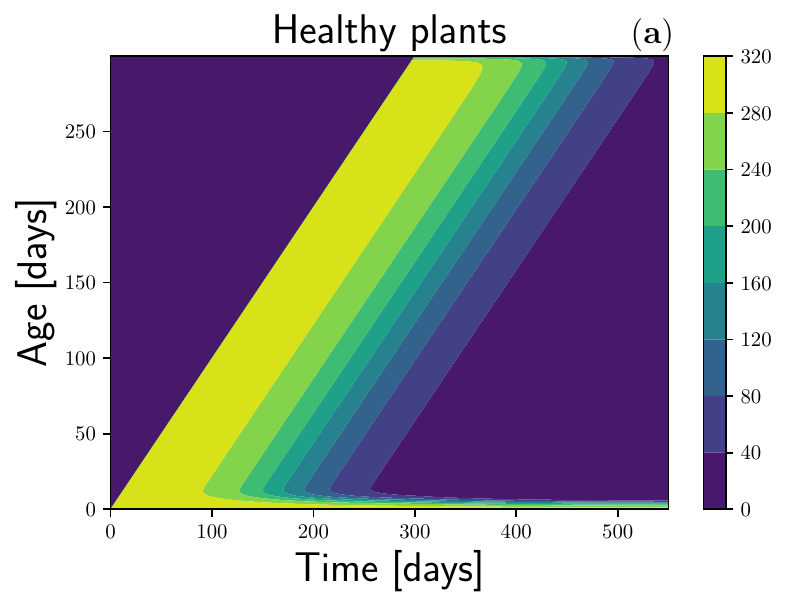}
\includegraphics[width=0.495\linewidth]{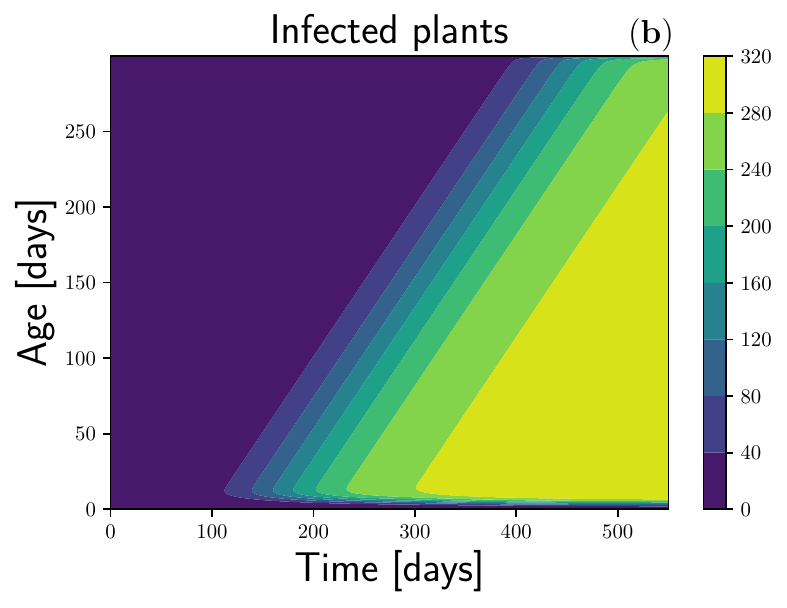}
\includegraphics[width=0.495\linewidth]{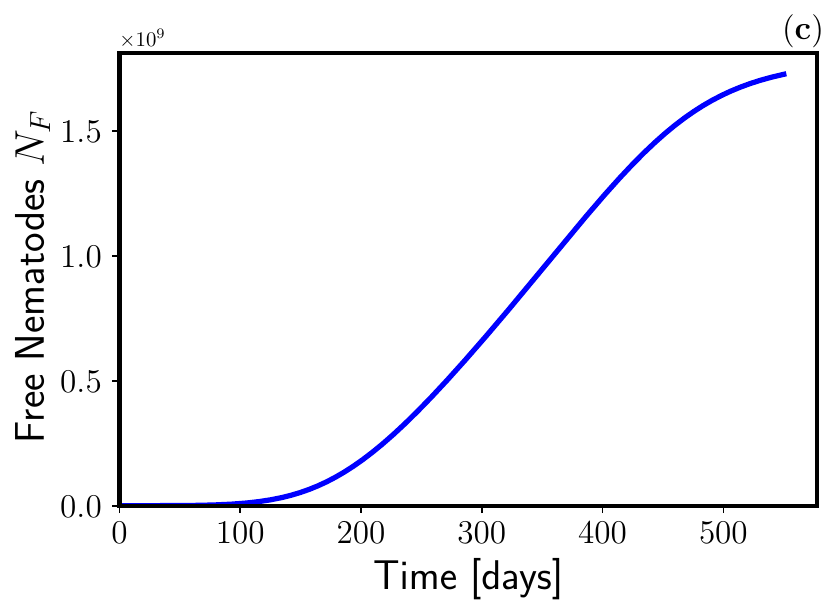}
\includegraphics[width=0.495\linewidth]{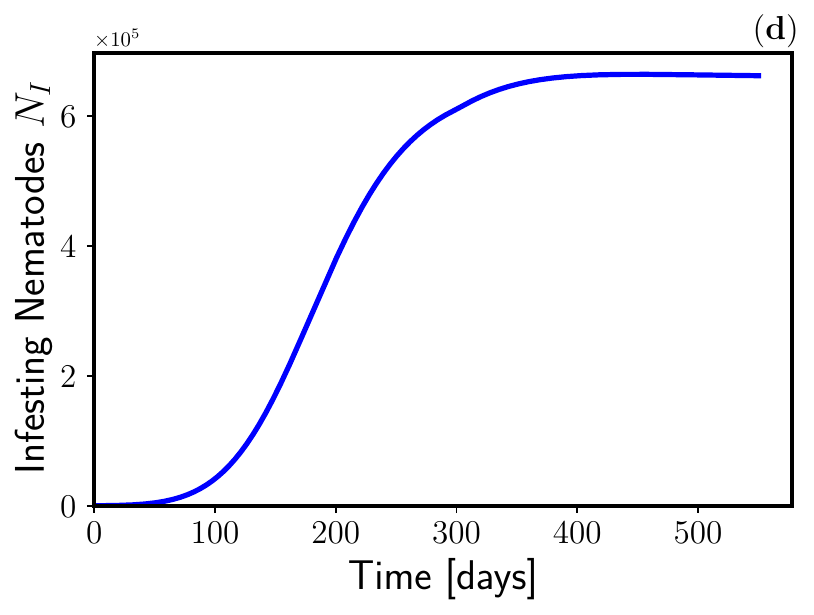} 
\caption{Simulation in the case $\mathcal{N} >1$ of the dynamics of healthy plants $S(a,t)$, infected plants $I(a,t)$, free nematodes $N_F(t)$ and Infesting nematodes $N_I(t)$ using the parameters given in Table~\ref{t1} and equations~\eqref{inc}--\eqref{par}. }
\label{eeq}
\end{center}
\end{figure}
Subfigures~\ref{eeq}(a) and~\ref{eeq}(b) illustrate the dynamics of healthy and infected banana plants as functions of age and time.
In Subfigure~\ref{eeq}(a), the density of healthy plants is initially high (yellow region), corresponding to young and abundant individuals. Over time, this density gradually decreases after approximately day~100 and remains low beyond day~300, indicating a reduction in the healthy population.
Conversely, Subfigure~\ref{eeq}(b) shows that the density of infected plants is initially almost zero (purple region), then increases markedly around day~100, reaching high levels (green and yellow regions) after day~300 for intermediate ages.
These results demonstrate how infestation progressively reduces the proportion of healthy plants while increasing that of infected plants.

 Similarly, Subfigures~\ref{eeq}(c) and~\ref{eeq}(d) show the dynamics of free nematodes in the soil and infesting nematodes within the plants, respectively.
These nematodes, which are responsible for root infestation, exhibit a rapid growth phase following an initial period of slow increase, and then reach a stable level, following a similar pattern.
Their proliferation directly contributes to the decline in healthy plants and the increase in infected plants observed in the previous figures. 

\begin{figure}[h]	
\begin{center}
\includegraphics[width=0.495\linewidth]{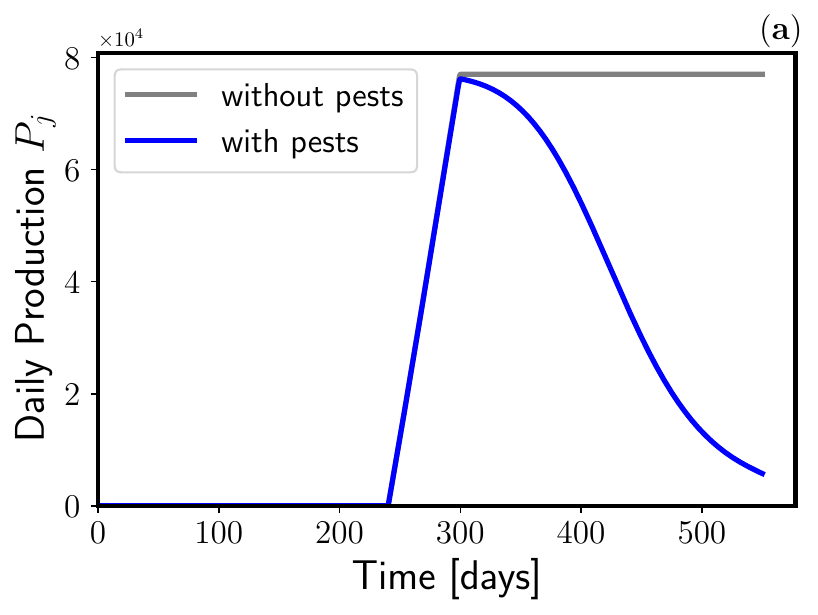}
\includegraphics[width=0.495\linewidth]{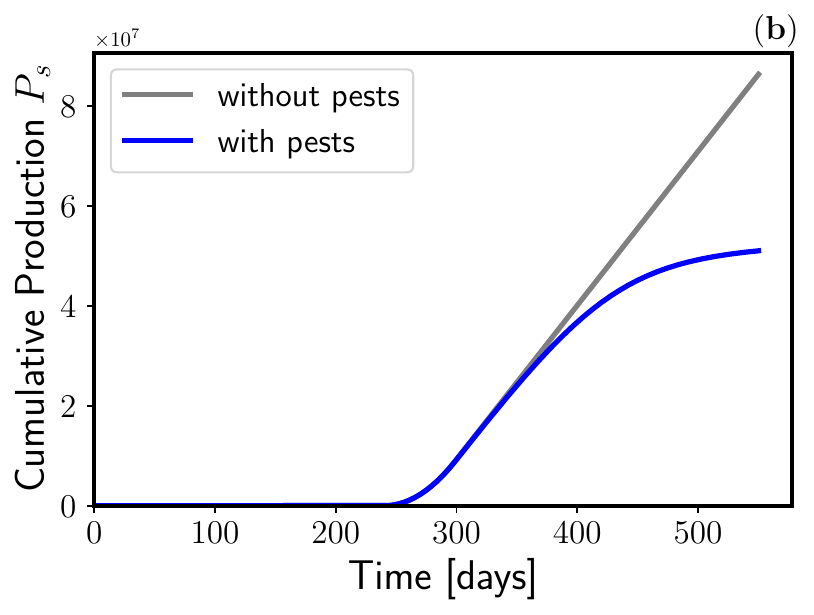}
\caption{Simulation of the daily and final cumulative production without pests (gray curves) 
and with pests (blue curves), using the parameters given in Table~\ref{t1} 
and equations~\eqref{inc}--\eqref{par}.
}
\label{pro-wthc}
\end{center}
\end{figure}

Finally, Subfigures~\ref{pro-wthc}(a) and~\ref{pro-wthc}(b) illustrate the impact of infestation on production.
In Subfigure~\ref{pro-wthc}(b), cumulative production increases steadily in the absence of parasites, whereas in the presence of infestation it quickly diverges from the ideal trajectory and stabilizes at a much lower level, indicating a substantial yield loss.
In Subfigure~\ref{pro-wthc}(b), daily production is initially high but gradually declines under the effect of infestation, becoming nearly zero by the end of the cropping season.

\subsubsection{Impulsive control application}

We aim to reduce the population of free nematodes in the soil by using nematicides such as the bionematicide Paecilomyces lilacinus strain 251 (non-toxic) \cite{a50} or the Sisal residue \cite{eco_control}.
The control variable will be represented by an impulsive function $u(t)$ which  acts directly on the mortality of free nematodes, incorporating itself into their dynamics via the term $-uN_F$ in the third equation of the model  system \eqref{Model diagram}. Hence, the controlled model is given as our model system where the third equation is substituted by   
\begin{equation}
\dot{N}_F(t)  = -\alpha(\lVert\beta(\cdot) S(\cdot\,,\,t)\lVert + e \lVert\beta(\cdot)I(\cdot\,,\,t))\lVert) \mathcal{W}(P,N_{F})+ \gamma\lVert I(\cdot\,,\,t) \lVert -( \mu_{F} + u(t))N_F(t) \,.
\end{equation}
The impulsive biological/ecological control consists of spraying nematicides two day per month, starting at the beginning of the season. The analysis aims to illustrate how such a regular, two-day-per-month intervention can effectively reduce the nematode population and protect production without resorting to a complex optimization process (which, despite its difficulties, offers its own advantages).
We consider that the application of nematicides occurs every $16$ days (two times per month) over a total period of $550$ days. Consequently, the   mathematical expression of $u$ is given as
\begin{equation}
u(t) \;=\; u_{\max}\sum_{k=0}^{\kappa} 
\mathbf{1}_{[\,16\,k,\;16\,k+1\,)}(t), \quad \text{with} \quad \kappa = \left\lfloor  550/16  \right\rfloor\,.
 \end{equation}

\begin{figure}[H]	
\begin{center}
\includegraphics[width=0.495\linewidth]{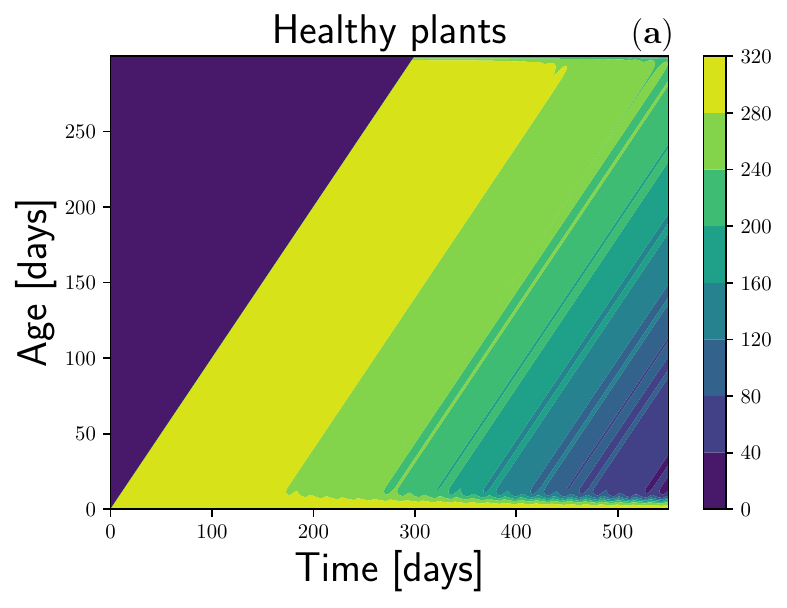}
\includegraphics[width=0.495\linewidth]{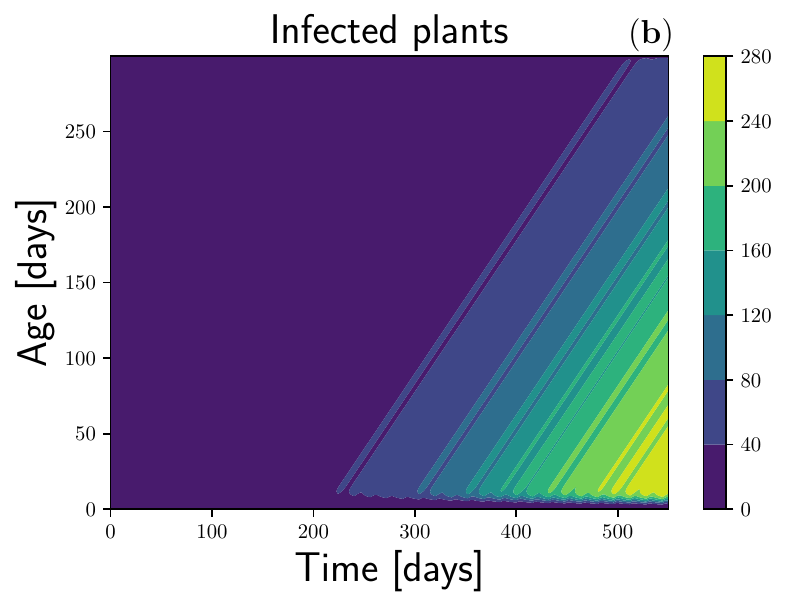}
\includegraphics[width=0.495\linewidth]{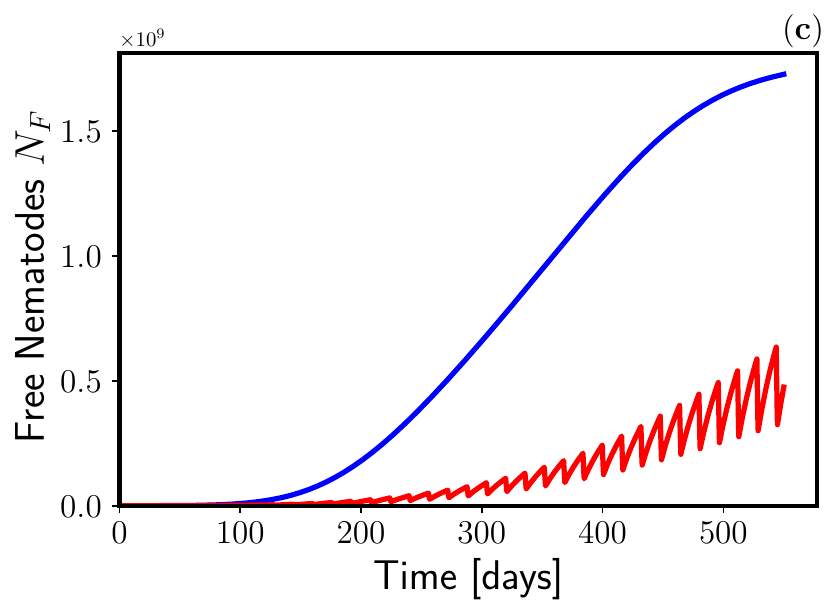}
\includegraphics[width=0.495\linewidth]{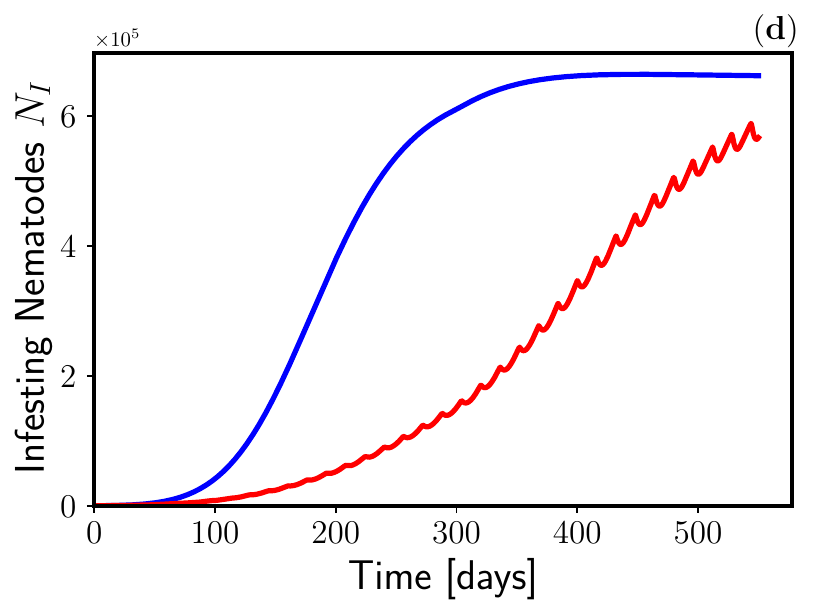}
\caption{Simulation of the dynamics of healthy plants $S(a,t)$ and infected plants $I(a,t)$ under control, along with a comparison of the dynamics of free and infesting nematodes without control (blue curves) and with control (red curves), using the parameters provided in Table~\ref{t1} and equations \eqref{inc}-\eqref{par}.}
\label{ste-wc}
\end{center}
\end{figure}

Figure~\ref{ste-wc} shows that the application of the specified control exerts an overall beneficial effect on the system dynamics. In the presence of control, the density of healthy plants rises sharply, the nematode population is drastically reduced, and the progression of infected plants slows down (as observed in \cite{a50} and\cite{eco_control}). These effects translate directly into a substantial improvement in production, which approaches the ideal, pest-free scenario (gray curves). 

\begin{figure}[H]	
\begin{center}
\includegraphics[width=0.495\linewidth]{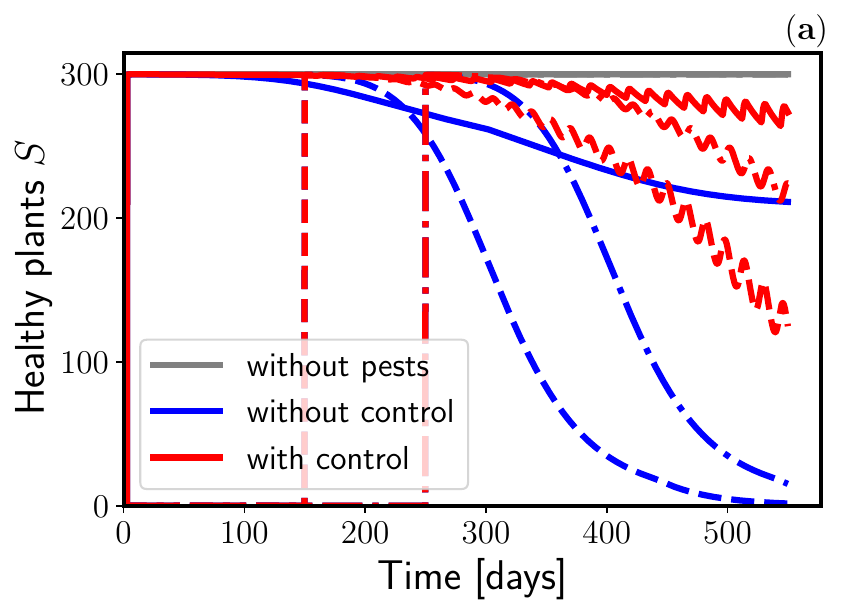}
\includegraphics[width=0.495\linewidth]{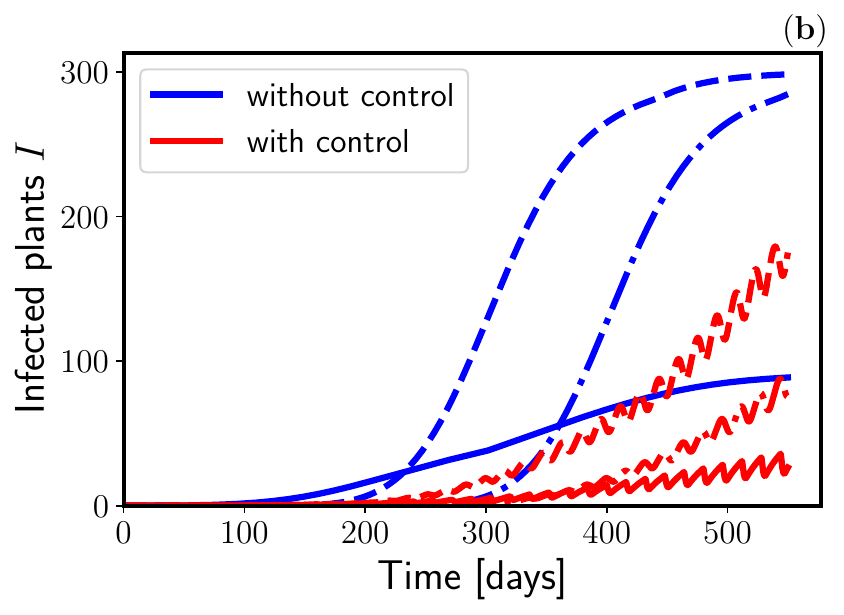}
\includegraphics[width=0.495\linewidth]{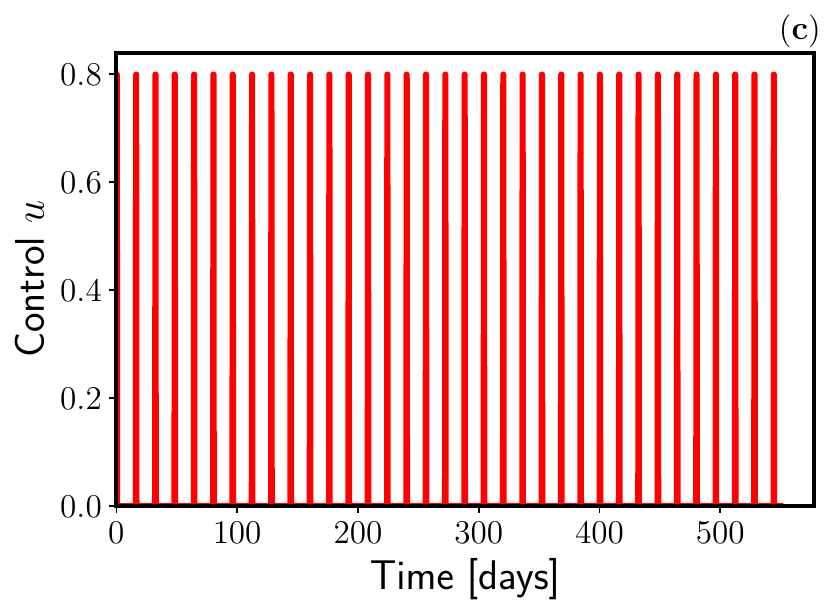}
\caption{Comparison of the dynamics of healthy and infected plants for selected age values without pests (gray curves), without control (blue curves) and with control (red curves), using parameters given in Table~\ref{t1} and equations \eqref{inc}-\eqref{par}.}
\label{slec-ages}
\end{center}
\end{figure}

In Figures~\ref{slec-ages} we represent the dynamics of healthy ans infected plants for few selected age values: age $12$ (solid curves), age $150$ (dashed curves) and age $250$ (thin mixed curves). This choice clearly illustrates the effect of control according to certain age classes, while avoiding graphic overload.
However, the beneficial effect of control varies depending on the age considered. For age $12$ (solid curves), the impact is very marked, with strong preservation of healthy plants and a clear limitation of infected plants. For age $150$ (dashed curves), the improvement remains significant but more moderate, and for age $250$ (thin mixed curves), the differences between the scenarios with and without control are less pronounced, although the advantage of control still persists.
\begin{figure}[H]	
\begin{center}
\includegraphics[width=0.495\linewidth]{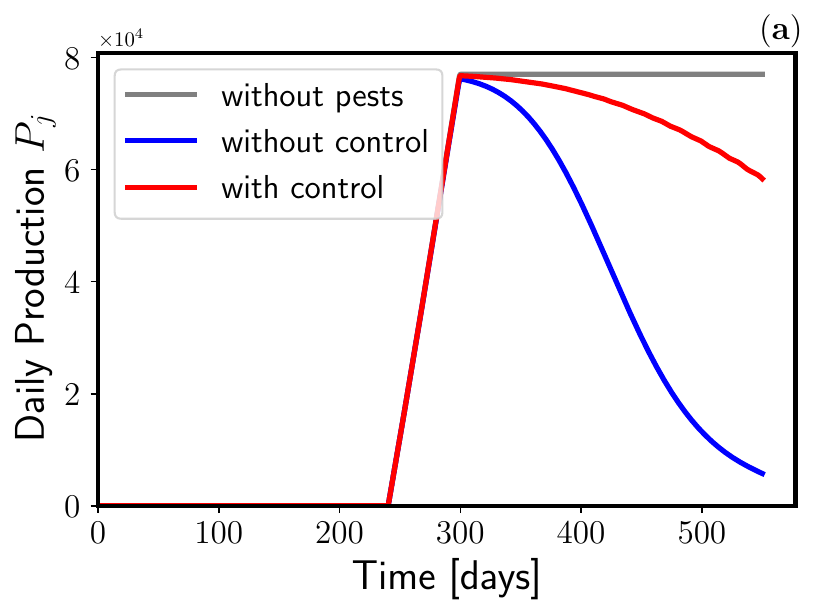}
\includegraphics[width=0.495\linewidth]{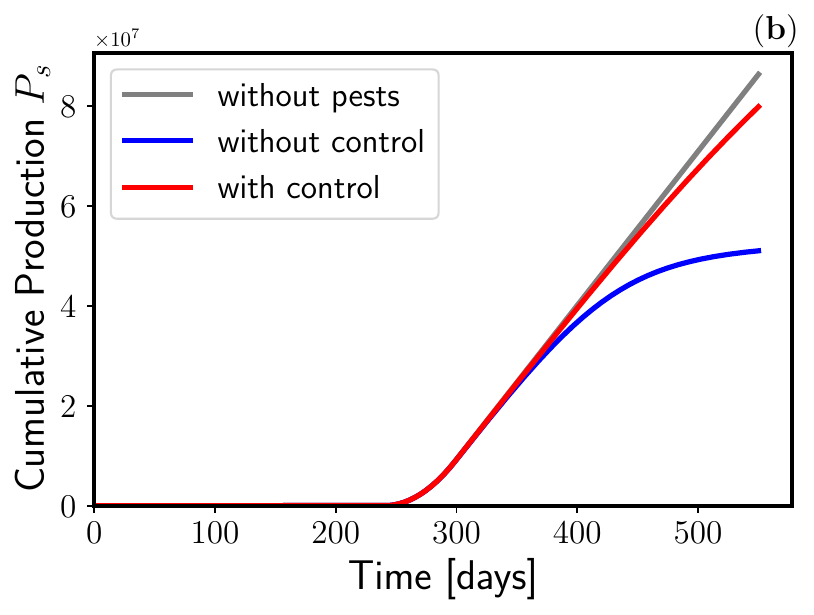}
\caption{Comparison of the dynamics of daily production and cumulative production without pest (gray curves), without control (blue curves) and with control (red curves), using parameters given in Table~\ref{t1} and equations \eqref{inc}-\eqref{par}.}
\label{pro-wc}
\end{center}
\end{figure}

The simulation results for production losses are estimated for a plot of $600$
 square meters with a spacing of $2.5$ meters between trees (recall that the initial conditions are $S_0(a) = 0$ for all $a \ne 0$, except $S_0(0) = 100$) and are presented in Table~\ref{tab2}:
the  crop loss rate is 40.91\% in the absence of any control measures. This is consistent with values reported in the literature, where losses of up to 60\% have been observed in many countries around the world~\cite{a4}. The difference may be attributed to the plant density and pest pressure considered in the model, as local environmental conditions were not included.
\begin{table}[h!]
\centering
\caption{Final cumulative production and associated losses/gains}

\begin{tabular}{|l|c|c|}
\hline
\textbf{Indicator} & \textbf{Value} & \textbf{Unit / \%} \\
\hline
Production without pests & $8.63 \times 10^7$ & kg \\
\hline 
Production without control & $5.09 \times 10^7$ & kg \\
\hline 
Production with control & $7.98 \times 10^7$ & kg \\
\hline
\hline 
Loss due to pests & 40.91 & \% \\
\hline 
Loss with control & 7.52 & \% \\
\hline 
Relative gain from control & 56.51 & \% \\
\hline
\end{tabular} 
\label{tab2}
\end{table}
In summary, the results of simulations (without and with control ) show that 
(i)  the nematode proliferation leads to a drastic reduction in healthy plants in favor of infected plants, resulting in a significant loss of yield, both in terms of cumulative and daily production as shown Table \ref{tab2}.
 and (ii) the control plays an essential role in maintaining a favorable balance in the system, as it not only increases the proportion of healthy plants but also reduces the dynamics of infected plants and nematodes. 

\section{Conclusion and discussion}

We have formulated and analyzed an age-structured mathematical model describing the infestation dynamics of \textit{Radopholus similis} in a banana-plantain  plantation, taking into account the age of the host plant. 
Our analysis showed that the model is ecologically well-posed; that is, it admits a unique, positive, and bounded solution. We identified the trivial equilibrium point and analyzed its stability properties. We found that the infection-free steady state (IFS) is locally asymptotically stable (L.A.S.) when $\mathcal{N} <1$, and unstable when  $\mathcal{N} \geq 1$. Moreover, it turns out that under suitable conditions on the recruitment and mortality rates, the IFS is globally asymptotically stable (G.A.S.).
Furthermore, we proposed a consistent numerical scheme, the validity of which was demonstrated both theoretically and through simulations. The numerical results indicate that the proposed model successfully captures the real-world dynamics of \textit{R. similis} infestation in banana-plantain  plantations, as reflected by the predicted yield loss rates, which are in close agreement with field observations. We  have implemented an impulsive control strategy, which confirms that the use of nematicides - whether chemical or biological - helps to mitigate the devastating effects of nematodes and enhances crop yield.

The expression of the basic reproduction number 
\[
\mathcal{N} = \frac{\gamma}{\sigma+\mu_F} 
\int_{0}^{a_{\dagger}} \int_{0}^{a} 
\ell(\zeta)\, e^{-\int_{\zeta}^{a}\mu(r)\,dr}\, d\zeta\, da
=
\frac{\gamma}{\sigma+\mu_F} 
\int_{0}^{a_{\dagger}} \ell(\zeta)
\left( \int_{\zeta}^{a_{\dagger}} e^{-\int_{\zeta}^{a} \mu(r)\,dr}\, da \right) d\zeta\,,
\]
provides a biologically meaningful threshold that determines the local stability of the infection-free steady state. 
Each component of this expression has a clear ecological interpretation. 
The parameter $\gamma$ represents the migration or infection rate of free nematodes in the soil, while $\mu_F$ denotes their natural mortality rate, and $\sigma$ accounts for additional losses such as trapping or environmental removal. 
Hence, the prefactor $\dfrac{\gamma}{\sigma+\mu_F}$ measures the average number of successful contacts made by a free nematode before it dies or is lost from the system. 
The double integral quantifies the total contribution of all root age classes to the production or maintenance of nematode populations. 
The function $\ell(\zeta)$ measures the reproductive potential or attractiveness of roots of age $\zeta$, while the exponential term 
$e^{-\int_{\zeta}^{a}\mu(r)\,dr}$ represents the probability that a root segment born at age $\zeta$ remains alive until age $a$, where $\mu(r)$ is the age-dependent natural mortality rate of the roots. 
Consequently, the integral accumulates the weighted reproductive contribution over all root ages, thereby linking nematode dynamics to the plant's age structure.

Mathematically, $\mathcal{N}$ is always positive and varies monotonically with biologically relevant parameters. 
It increases with the migration rate $\gamma$ and the root productivity function $\ell(\zeta)$, and decreases with higher mortality rates of free nematodes ($\mu_F$) or with stronger environmental losses ($\sigma$). 
An increase in the root mortality rate $\mu(a)$ also reduces $\mathcal{N}$, reflecting that short-lived roots provide a smaller window for nematode reproduction.
The threshold condition $\mathcal{N} < 1$ implies that, on average, each free nematode produces less than one successor before being lost, leading to the extinction of the infestation. 
Conversely, when $\mathcal{N} > 1$, the nematode population can invade and persist in the plantation. 
This threshold has clear biological and management implications: practices that reduce nematode migration ($\gamma$), increase mortality of free stages ($\mu_F$), or shorten the lifespan of susceptible roots (increase $\mu(a)$) will contribute to maintaining $\mathcal{N}$ below unity.

Overall, the structure of $\mathcal{N}$ is consistent with ecological intuition and measurable in practice, as its parameters can be estimated through experimental or field data. 
The age-structured formulation therefore provides a realistic and flexible framework for describing the population dynamics of \textit{Radopholus similis} in banana or plantain plantations, and for designing control strategies based on biologically meaningful quantities.

It should be noted that the determination and analysis of the endemic equilibrium has not been carried out in this study. 
However, numerical simulations suggest the existence of such an equilibrium (the parameters used for numerical simulations are such that $\mathcal{N}=8.86 >1 $).
We plan to address this in a forthcoming paper, which will also include the investigation of optimal control strategies, 
guided by the insights provided by the basic reproduction number $\mathcal{N}$.

{\tiny {\tiny }} \clearpage
 \addcontentsline{toc}{chapter}{Bibliography}

 \bibliographystyle{plain}
 \bibliographystyle{plainurl} 
 
\bibliography{data} 
\end{document}